\documentclass{article}

\usepackage{microtype}
\usepackage{graphicx}
\usepackage{subfigure}
\usepackage{booktabs} 
\usepackage{amssymb}
\usepackage{amsmath}
\usepackage{times}
\usepackage{enumerate}
\usepackage{algorithm}
\usepackage{algorithmic}
\usepackage{multirow}
\usepackage{bm}
\usepackage{mathtools}
\usepackage{amsthm}
\usepackage[pagebackref=true,breaklinks=true,letterpaper=true,colorlinks,bookmarks=false]{hyperref}
\usepackage{setspace}
\usepackage{algorithm}
\usepackage{algorithmic}
\usepackage{epsfig}
\usepackage{epstopdf}
\usepackage{bm}
\usepackage{url}
\usepackage{natbib}
\include{standardHeader}

\newtheorem{theorem}{Theorem}

\newtheorem{lemma}{Lemma}
\newtheorem{corollary}{Corollary}

\newtheorem{assumption}{Assumption}

\newcommand{\x}{\mathbf{x}}
\newcommand{\y}{\mathbf{y}}

\newcommand{\1}{\mathbf{1}}
\newcommand{\0}{\mathbf{0}}

\renewcommand{\v}{\mathbf{v}}

\renewcommand{\b}{b}

\newcommand{\A}{A}

\newcommand{\M}{\mathbf{M}}
\newcommand{\W}{W}

\newcommand{\I}{I}
\newcommand{\U}{U}

\newcommand{\R}{\mathbb{R}}
\newcommand{\N}{\mathcal{N}}
\newcommand{\blambda}{\mathbf{\lambda}}

\newcommand{\<}{\left\langle}
\renewcommand{\>}{\right\rangle}

\DeclareMathOperator*{\argmin}{argmin}

\usepackage[accepted]{icml2020}

\icmltitlerunning{Revisiting EXTRA for Smooth Distributed Optimization}

\begin{document}

\onecolumn
\icmltitle{Revisiting EXTRA for Smooth Distributed Optimization}

\begin{icmlauthorlist}
\icmlauthor{Huan Li}{to,goo}
\icmlauthor{Zhouchen Lin}{goo}
\end{icmlauthorlist}

\icmlaffiliation{to}{Institute of Robotics and Automatic Information Systems, Nankai University, Tianjin, China (lihuan\_ss@126.com).\\}
\icmlaffiliation{goo}{Key Lab. of Machine Perception (MOE), School of EECS, Peking University, Beijing, China
    (zlin@pku.edu.cn).\\}
\icmlcorrespondingauthor{Zhouchen Lin}{zlin@pku.edu.cn}





\vskip 0.3in



\printAffiliationsAndNotice{}  

\begin{abstract}
  EXTRA is a popular method for dencentralized distributed optimization and has broad applications. This paper revisits EXTRA. First, we give a sharp complexity analysis for EXTRA with the improved $O((\frac{L}{\mu}+\frac{1}{1-\sigma_2(\W)})\log\frac{1}{\epsilon(1-\sigma_2(\W))})$ communication and computation complexities for $\mu$-strongly convex and $L$-smooth problems, where $\sigma_2(\W)$ is the second largest singular value of the weight matrix $\W$. When the strong convexity is absent, we prove the $O((\frac{L}{\epsilon}+\frac{1}{1-\sigma_2(\W)})\log\frac{1}{1-\sigma_2(\W)})$ complexities. Then, we use the Catalyst framework to accelerate EXTRA and obtain the $O(\sqrt{\frac{L}{\mu(1-\sigma_2(\W))}}\log\frac{ L}{\mu(1-\sigma_2(\W))}\log\frac{1}{\epsilon})$ communication and computation complexities for strongly convex and smooth problems and the $O(\sqrt{\frac{L}{\epsilon(1-\sigma_2(\W))}}\log\frac{1}{\epsilon(1-\sigma_2(\W))})$ complexities for nonstrongly convex ones. Our communication complexities of the accelerated EXTRA are only worse by the factors of $(\log\frac{L}{\mu(1-\sigma_2(\W))})$ and $(\log\frac{1}{\epsilon(1-\sigma_2(\W))})$ from the lower complexity bounds for strongly convex and nonstrongly convex problems, respectively.
\end{abstract}

\section{Introduction}
In this paper, we consider the following convex problem
\begin{eqnarray}
\min_{x\in\R^n} F(x)=\frac{1}{m}\sum_{i=1}^m f_i(x)\label{problem}
\end{eqnarray}
in the decentralized distributed environment, where $m$ agents form an undirected communication network and collaboratively solve the above problem. Each agent $i$ privately holds a local objective function $f_i(x)$ and can exchange information only with its immediate neighbors. We only consider the network that does not have a centralized agent. Distributed computation has broad applications, ranging from machine learning \cite{Dekel2012,Forero2012,Agarwal2011,Recht2011}, to sensor networks \cite{Duarte2014}, to flow and power control problems \cite{Duarte2014,Gan2013}.
\subsection{Literature Review}
Distributed optimization has gained significant attention in engineering applications for a long time \cite{Bertsekas1983,Tsitsiklis1986}. The distributed subgradient method was first proposed in \cite{Nedic-2009} with the convergence and convergence rate analysis for the general network topology and further extended to the asynchronous variant in \cite{nedic2011asynchronous}, the stochastic variant in \cite{Ram-2010}, and a study with fixed step-size in \cite{Yuan-2016}. In \cite{Jakovetic-2014,chen2012b}, the accelerated distributed gradient method in the sense of Nesterov has been proposed, and the authors of \cite{li-2018-pm} gave a different explanation with sharper analysis, which builds upon the accelerated penalty method. Although the optimal computation complexity and near optimal communication complexity were proved in \cite{li-2018-pm}, the accelerated distributed gradient method employs multiple consensus after each gradient computation and thus places more burdens in the communication-limited environment.

A different class of distributed approaches with efficient communication is based on the Lagrangian dual and they work in the dual space. Classical algorithms include dual ascent \cite{Terelius-2011,dasent,Uribe-2017}, ADMM \cite{Iutzeler-2016,makhdoumi-2017,Aybat2018}, and the primal-dual method \cite{Lam-2017,scaman-2018,hong-2017,jakovetic-2017}. Specifically, accelerated dual ascent \cite{dasent} and the primal-dual method \cite{scaman-2018} attain the optimal communication complexities for smooth and nonsmooth problems, respectively. However, the dual-based methods require the evaluation of the Fenchel conjugate or the proximal mapping and thus have a larger computation cost per iteration.

EXTRA \cite{Shi-2015} and the gradient tracking based method \cite{aug-dgm,qu2017} (also called DIGing in \cite{shi2017}) can be seen as a trade-off between communications and computations, which need equal numbers of communications and gradient computations at each iteration. As a comparison, the accelerated distributed gradient method needs more communications, while the dual-based methods require more computations. EXTRA uses the differences of gradients and guarantees the convergence to the exact optimal solution with constant step-size. The proximal-gradient variant was studied in \cite{Shi-2015-2}. Recently, researchers have established the equivalence between the primal-dual method and EXTRA \cite{hong-2017,Mokhtari-2016,jakovetic-2017}. Specifically, the authors of \cite{hong-2017} study the nonconvex problem and the authors of \cite{Mokhtari-2016} focus on the stochastic optimization, while the authors of \cite{jakovetic-2017} give a unified framework for EXTRA and the gradient tracking based method. The gradient tracking based method shares some similar features to EXTRA, e.g., using the differences of gradients and constant step-size. The accelerated version of the gradient tracking based method was studied in \cite{qu2017-2}.

In this paper, we revisit EXTRA and give a sharper complexity analysis for the original EXTRA. Then, we propose an accelerated EXTRA, which answers the open problem proposed in [Section V]\cite{Shi-2015-2} on how to improve the rate of EXTRA with certain acceleration techniques.
\subsection{Notation and Assumption}
Denote $x_{(i)}\in\R^n$ to be the local copy of the variable $x$ for agent $i$ and $x_{(1:m)}$ to be the set of vectors consisting of $x_{(1)},...,x_{(m)}$. We introduce the aggregate objective function $f(\x)$ of the local variables with its argument $\x\in\R^{m\times n}$ and gradient $\nabla f(\x)\in\R^{m\times n}$ as
\begin{equation}\label{define_f}
f(\x)=\sum_{i=1}^m f_i(x_{(i)}),\qquad\x=\left(
  \begin{array}{c}
    x_{(1)}^T\\
    \vdots\\
    x_{(m)}^T
  \end{array}
\right),\qquad\nabla f(\x)=\left(
  \begin{array}{c}
    \nabla f_1(x_{(1)})^T\\
    \vdots\\
    \nabla f_m(x_{(m)})^T
  \end{array}
\right).
\end{equation}
For a given matrix, we use $\|\cdot\|_F$ and $\|\cdot\|_2$ to denote its Frobenius norm and spectral norm, respectively. We denote $\|\cdot\|$ as the $l_2$ Euclidean norm for a vector. Denote $\I\in\R^{m\times m}$ as the identity matrix and $\1=(1,1,...,1)^T\in\R^m$ as the vector with all ones. For any matrix $\x$, we denote its average across the rows as
\begin{eqnarray}
\alpha(\x)=\frac{1}{m}\sum_{i=1}^mx_{(i)}.\label{def_avg_x}
\end{eqnarray}
Define two operators measuring the consensus violation. The first one is
\begin{eqnarray}
\Pi=\I-\frac{1}{m}\1\1^T\in\R^{m\times m}\label{define_pi}
\end{eqnarray}
and $\|\Pi\x\|_F$ measures the distance between $x_{(i)}$ and $\alpha(\x)$ for all $i$. The second one follows \cite{Shi-2015},
\begin{eqnarray}
\U=\sqrt{\frac{\I-\W}{2}}\in\R^{m\times m}.\label{define_u}
\end{eqnarray}
Let $\N_i$ be the neighbors of agent $i$ and $\mbox{Span}(\U)$ be the linear span of all the columns of $\U$.

We make the following assumptions for the local objectives.
\begin{assumption}\label{assumption_f}
\begin{enumerate}
\item Each $f_i(x)$ is $\mu$-strongly convex: $f_i(y)\geq f_i(x)+\<\nabla f_i(x),y-x\>+\frac{\mu}{2}\|y-x\|^2$. Especially, $\mu$ can be zero, and we say $f_i(x)$ is convex in this case.
\item Each $f_i(x)$ is $L$-smooth: $f_i(y)\leq f_i(x)+\<\nabla f_i(x),y-x\>+\frac{L}{2}\|y-x\|^2$.
\end{enumerate}
\end{assumption}
Then, $F(x)$ and $f(\x)$ are also $\mu$-strongly convex and $L$-smooth. Assume that the set of minimizers of problem (\ref{problem}) is nonempty. Denote $x^*$ as one minimizer, and let $\x^*=\1(x^*)^T$.

We make the following assumptions for the weight matrix $\W\in\R^{m\times m}$ associated to the network.
\begin{assumption}\label{assumption_w}
\begin{enumerate}
\item $\W_{i,j}\neq 0$ if and only if agents $i$ and $j$ are neighbors or $i=j$. Otherwise, $\W_{i,j}=0$.
\item $\W=\W^T$, $\I\succeq\W\succeq -\I$, and $\W\1=\1$.
\item $\sigma_2(\W)<1$, where $\sigma_2(\W)$ is the second largest singular value of $\W$.
\end{enumerate}
\end{assumption}
Part 2 of Assumption \ref{assumption_w} implies that the singular values of $\W$ lie in $[0,1]$ and its largest one $\sigma_1(\W)$ equals 1. Moreover, part 3 can be deduced by part 2 and the assumption that the network is connected. Examples satisfying Assumption \ref{assumption_w} can be found in \cite{Shi-2015}.

When minimizing a convex function, the performance of the first-order methods is affected by the smoothness constant $L$ and the strong convexity constant $\mu$, as well as the target accuracy $\epsilon$. When we solve the problem over a network, the connectivity of the network also directly affects the performance. Typically, $\frac{1}{1-\sigma_2(\W)}$ is a good indication of the network connectivity \cite{Jakovetic-2014,dasent} and it is often related to $m$ [Proposition 5]\cite{nedic-2018}. For example, for any connected and undirected graph, $\frac{1}{1-\sigma_2(\W)}\leq m^2$ \cite{nedic-2018}. In this paper, we study the complexity of EXTRA with explicit dependence on $L$, $\mu$, $1-\sigma_2(\W)$, and $\epsilon$.

Denote $x_{(1:m)}^0$ to be the initializers. Assume that $\|x_{(i)}^0-x^*\|^2\leq R_1$, $\|x^*\|^2\leq R_1$, and $\|\nabla f_i(x^*)\|^2\leq R_2$ for all $i=1,\cdots,m$. Then we can simply have
\begin{eqnarray}
\begin{aligned}\label{fact}
\|\x^0-\x^*\|_F^2\leq mR_1,\quad\|\x^*\|_F^2\leq mR_1\quad \mbox{and}\quad\|\nabla f(\x^*)\|_F^2\leq mR_2.
\end{aligned}
\end{eqnarray}
In this paper, we only regard $R_1$ and $R_2$ as the constants which can be dropped in our complexities.

\subsection{Proposed Algorithm}
\begin{algorithm}[t]
   \caption{EXTRA}
   \label{extra}
\begin{algorithmic}
   \STATE Input $F(x)$, $K$, $x_{(1:m)}^0$, $v_{(1:m)}^0$
   \FOR{$k=0,1,2,\cdots,K$}
   \STATE $x_{(i)}^{k+1}=x_{(i)}^k-\alpha\left(\nabla f_i(x_{(i)}^k)+v_{(i)}^k+\frac{\beta}{2}\left(x_{(i)}^k-\sum_{j\in\N_i}\W_{i,j}x_{(j)}^k\right)\right)\forall i.$
   \STATE $v_{(i)}^{k+1}=v_{(i)}^k+\frac{\beta}{2}\left(x_{(i)}^{k+1}-\sum_{j\in\N_i}\W_{i,j}x_{(j)}^{k+1}\right)\forall i.$
   \ENDFOR
   \STATE Output $x_{(1:m)}^{K+1}$ and $v_{(1:m)}^{K+1}$.
\end{algorithmic}
\end{algorithm}

Before presenting the proposed algorithm, we first rewrite EXTRA in the primal-dual framework in Algorithm \ref{extra}. When we set $\alpha=\frac{1}{\beta}$, Algorithm \ref{extra} reduces to the original EXTRA. \footnote{Initialize $\v^0=\0$ and define $\widetilde\W=\frac{\I+\W}{2}$. The second step of Algorithm \ref{extra} leads to $\v^k=\beta\sum_{t=1}^k(\widetilde\W-\W)\x^k$. Plugging it into the first step and letting $\alpha=1/\beta$ leads to equation (3.5) in \cite{Shi-2015}.} In this paper, we specify $\alpha=\frac{1}{2(L+\beta)}$ and $\beta=L$ for the strongly convex problems to give a faster convergence rate than the original EXTRA, which is crucial to obtain the near optimal communication complexities after acceleration.

We use the Catalyst framework \cite{catalyst} to accelerate Algorithm \ref{extra}. It has double loops and is described in Algorithm \ref{acc-extra}. The inner loop calls Algorithm \ref{extra} to approximately minimize a well-chosen auxiliary function of $G^k(x)$ for $T_k$ iterations with warm-start. $T_k$ and $\tau$ are given for two cases:
\begin{enumerate}
\item When each $f_i(x)$ is strongly convex with $\mu>0$, then $\tau=L(1-\sigma_2(\W))-\mu>0$ and $T_k=O(\frac{1}{1-\sigma_2(\W)}\log\frac{L}{\mu(1-\sigma_2(\W))})$, which is a constant.
\item When each $f_i(x)$ is convex with $\mu=0$, then $\tau=L(1-\sigma_2(\W))$ and $T_k=O(\frac{1}{1-\sigma_2(\W)}\log \frac{k}{1-\sigma_2(\W)})$, which is nearly a constant.
\end{enumerate}
Although Algorithm \ref{acc-extra} employs the double loop, it places almost no more burdens than the original EXTRA. A good property of Algorithms \ref{extra} and \ref{acc-extra} in practice is that they need equal numbers of gradient computations and communications at each iterations.

\begin{algorithm}[t]
   \caption{Accelerated EXTRA}
   \label{acc-extra}
\begin{algorithmic}
   \STATE Initialize $x_{(i)}^0=y_{(i)}^0$, $v_{(i)}^0=\0$, $q=\frac{\mu}{\mu+\tau}$; set $\theta_k=\sqrt{q},\forall k$ if $\mu>0$; otherwise, set $\theta_0=1$ and update $\theta_{k+1}\in (0,1)$ by solving the equation $\theta_{k+1}^2=(1-\theta_{k+1})\theta_k^2$.
   \FOR{$k=0,1,2,\cdots,K$}
   \STATE Define $g_i^k(x)= f_i(x)+\frac{\tau}{2}\|x-y_{(i)}^k\|^2$ and $G^k(x)=\frac{1}{m}\sum_{i=1}^m g_i^k(x)$.
   \STATE $(x_{(1:m)}^{k+1},v_{(1:m)}^{k+1})=\mbox{EXTRA}(G^k(x),T_k,x_{(1:m)}^k,v_{(1:m)}^k)$.
   \STATE $y_{(i)}^{k+1}=x_{(i)}^{k+1}+\frac{\theta_k(1-\theta_k)}{\theta_k^2+\theta_{k+1}}(x_{(i)}^{k+1}-x_{(i)}^k) \forall i$.
   \ENDFOR
\end{algorithmic}
\end{algorithm}

\subsection{Complexities}
We study the communication and computation complexities of EXTRA and its accelerated version in this paper. They are presented as the numbers of communications and computations to find an $\epsilon$-optimal solution $x$ such that $F(x)-F(x^*)\leq\epsilon$. We follow \cite{li-2018-pm} to define one communication to be the operation that all the agents receive information from their neighbors once, i.e., $\sum_{j\in\N_i}\W_{ij}x_{(j)}$ for all $i=1,2,...,m$. One computation is defined to be the gradient evaluations of all the agents once, i.e., $\nabla f_i(x_{(i)})$ for all $i$. Note that the gradients are evaluated in parallel on each nodes.

To find an $\epsilon$-optimal solution, Algorithm \ref{extra} needs $O((\frac{L}{\mu}+\frac{1}{1-\sigma_2(\W)})\log\frac{1}{\epsilon(1-\sigma_2(\W))})$ and $O((\frac{L}{\epsilon}+\frac{1}{1-\sigma_2(\W)})\log\frac{1}{1-\sigma_2(\W)})$ iterations for strongly convex and nonstrongly convex problems, respectively. The computation and communication
complexities are identical for EXTRA, which equal the number of iterations. For Algorithm \ref{acc-extra}, we establish the $O(\sqrt{\frac{L}{\mu(1-\sigma_2(\W))}}\log\frac{L}{\mu(1-\sigma_2(\W))}\log\frac{1}{\epsilon})$ complexity for strongly convex problems and the $O(\sqrt{\frac{L}{\epsilon(1-\sigma_2(\W))}}\log\frac{1}{\epsilon(1-\sigma_2(\W))})$ one for nonstrongly convex problems.

Our first contribution is to give a sharp analysis for EXTRA with improved complexity. The complexity of the original EXTRA is at least $O(\frac{L^2}{\mu^2(1-\sigma_2(\W))}\log\frac{1}{\epsilon})$\footnote{The authors of \cite{Shi-2015} did not give an explicit complexity. We try to simplify equation (3.38) in \cite{Shi-2015} and find it to be at least $O(\frac{L^2}{\mu^2(1-\sigma_2(\W))}\log\frac{1}{\epsilon})$. The true complexity may be larger than $O\left(\frac{L^2}{\mu^2(1-\sigma_2(\W))}\log\frac{1}{\epsilon}\right)$.} for strongly convex problems. For nonstrongly convex ones, although the $O(\frac{1}{\epsilon})$ complexity was studied in \cite{Shi-2015}, no explicit dependence on $1-\sigma_2(\W)$ was given.\footnote{The authors of \cite{Shi-2015} proved the $O(\frac{1}{K})$ rate in the sense of $\frac{1}{K}\sum_{k=1}^K\|\U\blambda^k+\nabla f(\x^k)\|_{\widetilde\W}^2\leq O(\frac{1}{K})$ and $\frac{1}{K}\sum_{k=1}^K\|\U\x^k\|_F^2\leq O(\frac{1}{K})$, where $\widetilde\W=\frac{\I+\W}{2}$. They omitted the dependence on $1-\sigma_2(\W)$ in their analysis.} It is remarkable that the sum of $\frac{L}{\mu}$ (or $\frac{L}{\epsilon}$) and $\frac{1}{1-\sigma_2(\W)}$, rather than their product, dominates our complexities. When we only consider the gradient computation time and if $\frac{1}{1-\sigma_2(\W)}$ is smaller than $\frac{L}{\mu}$ (or $\frac{L}{\epsilon}$), we can see that due to parallelization, EXTRA takes almost only $\frac{1}{m}$ of computation time compared to nondistributed gradient descent, whose complexities are $O(\frac{mL}{\mu}\log\frac{1}{\epsilon})$ and $O(\frac{mL}{\epsilon})$ for strongly convex and nonstrongly convex problems, respectively \cite{nesterov2013introductory}. Thus, EXTRA achieves a linear speed up if we ignore the logarithm factor.



Our second contribution is to give an accelerated EXTRA with the near optimal communication complexity and a competitive computation complexity. In Table \ref{table1}, we summarize the comparisons to the state-of-the-art decentralized optimization algorithms, namely, the accelerated dual ascent method and the accelerated penalty method with consensus. We also present the complexities of the nonaccelerated EXTRA and the lower complexity bounds. Our communication complexities of the accelerated EXTRA match the lower bounds, except the extra factors of $\left(\log\frac{L}{\mu(1-\sigma_2(\W))}\right)$ and $\left(\log\frac{1}{\epsilon(1-\sigma_2(\W))}\right)$ for strongly convex and nonstrongly convex problems, respectively. When high precision is required, i.e., $\frac{1}{\epsilon}>\frac{L}{\mu(1-\sigma_2(\W))}$ for strongly convex problems and $\frac{1}{\epsilon}>\frac{1}{1-\sigma_2(\W)}$ for nonstrongly convex problems, our communication complexities are competitive to the state-of-the-art ones in \cite{dasent,Uribe-2017,li-2018-pm}.
On the other hand, our computation complexities are better than those of \cite{Uribe-2017} for applications with large $\frac{L}{\epsilon}$ and $\frac{L}{\mu}$ and moderate $\log\frac{1}{1-\sigma_2(\W)}$, but worse than those of \cite{li-2018-pm}.\footnotemark[5]\footnotetext[4]{The dependence on $\frac{1}{1-\sigma_2(\W)}$ is omitted in \cite{Shi-2015}.}\footnotetext[5]{Although the authors of \cite{dasent} also gives the $O\left(\sqrt{\frac{L}{\mu}}\log\frac{1}{\epsilon}\right)$ computation complexity, they defines one computation to be the cost of solving an `$\argmin$' subproblem. We cite \cite{Uribe-2017} in Table \ref{table1}, where the authors study the computation complexity with the total number of gradient computations, which is a more reasonable measurement.} Our result is a significant complement to the existing work in the sense that EXTRA and its accelerated version have equal numbers of communications and computations, while the accelerated dual ascent has more computations than communications and the accelerated penalty method needs more communications than computations.

\begin{table*}
\begin{center}
\scriptsize
\begin{tabular}{|c|c|c|}
\hline
\multicolumn{3}{c}{Non-strongly convex case}\\
\hline
 Methods & Complexity of gradient computation & Complexity of communication\\
\hline
\cite{Shi-2015}'s result for EXTRA\footnotemark[4] & $O\left(\frac{1}{\epsilon}\right)$\hspace*{0.1cm} \cite{Shi-2015}& $O\left(\frac{1}{\epsilon}\right)$\hspace*{0.1cm} \cite{Shi-2015}\\
Our result for EXTRA & $O\left(\left(\frac{L}{\epsilon}+\frac{1}{1-\sigma_2(\W)}\right)\log\frac{1}{1-\sigma_2(\W)}\right)$ & $O\left(\left(\frac{L}{\epsilon}+\frac{1}{1-\sigma_2(\W)}\right)\log\frac{1}{1-\sigma_2(\W)}\right)$\\
Accelerated Dual Ascent& $O\left(\frac{L}{\epsilon\sqrt{1-\sigma_2(\W)}}\log^2\frac{1}{\epsilon}\right)$ \cite{Uribe-2017} & $O\left(\sqrt{\frac{L}{\epsilon(1-\sigma_2(\W))}}\log\frac{1}{\epsilon}\right)$ \cite{Uribe-2017}\\
\hspace*{-0.4cm}Accelerated Penalty Method \hspace*{-0.4cm}& $O\left({\sqrt{\frac{L}\epsilon}}\right)$ \cite{li-2018-pm}& $O\left(\sqrt{\frac{L}{\epsilon(1-\sigma_2(\W))}}\log\frac{1}{\epsilon}\right)$ \cite{li-2018-pm}\\
Our Accelerated EXTRA & \hspace*{-0.2cm}$O\left(\sqrt{\frac{L}{\epsilon(1-\sigma_2(\W))}}\log\frac{1}{\epsilon(1-\sigma_2(\W))}\right)$ \hspace*{-0.2cm}& \hspace*{-0.2cm} $O\left(\sqrt{\frac{L}{\epsilon(1-\sigma_2(\W))}}\log\frac{1}{\epsilon(1-\sigma_2(\W))}\right)$\hspace*{-0.2cm}\\
\hline
Lower Bound & $O\left({\sqrt{\frac{L}\epsilon}}\right)$ \cite{nesterov2013introductory} & $O\left(\sqrt{\frac{L}{\epsilon(1-\sigma_2(\W))}}\right)$ \cite{scaman-2019}\\
\hline
\multicolumn{3}{c}{Strongly convex case}\\
\hline
 Methods & Complexity of gradient computation & Complexity of communication\\
\hline
\cite{Shi-2015}'s result for EXTRA & at least $O\left(\frac{L^2}{\mu^2(1-\sigma_2(\W))}\log\frac{1}{\epsilon}\right)$ \cite{Shi-2015} & at least $ O\left(\frac{L^2}{\mu^2(1-\sigma_2(\W))}\log\frac{1}{\epsilon}\right)$ \cite{Shi-2015}\\
Our result for EXTRA & $O\left(\left(\frac{L}{\mu}+\frac{1}{1-\sigma_2(\W)}\right)\log\frac{1}{\epsilon(1-\sigma_2(\W))}\right)$ & $O\left(\left(\frac{L}{\mu}+\frac{1}{1-\sigma_2(\W)}\right)\log\frac{1}{\epsilon(1-\sigma_2(\W))}\right)$\\
Accelerated Dual Ascent& $O\left(\frac{L}{\mu\sqrt{1-\sigma_2(\W)}}\log^2\frac{1}{\epsilon}\right)$ \cite{Uribe-2017}& $O\left(\sqrt{\frac{L}{\mu(1-\sigma_2(\W))}}\log\frac{1}{\epsilon}\right)$ \cite{dasent,Uribe-2017}\\
\hspace*{-0.4cm}Accelerated Penalty Method\hspace*{-0.4cm} & $O\left(\sqrt{\frac{L}{\mu}}\log\frac{1}{\epsilon}\right)$ \cite{li-2018-pm}& $O\left(\sqrt{\frac{L}{\mu(1-\sigma_2(\W))}}\log^2\frac{1}{\epsilon}\right)$ \cite{li-2018-pm}\\
Our Accelerated EXTRA & \hspace*{-0.2cm}$O\left(\sqrt{\frac{L}{\mu(1-\sigma_2(\W))}}\log\frac{L}{\mu(1-\sigma_2(\W))}\log\frac{1}{\epsilon}\right)$\hspace*{-0.2cm} & \hspace*{-0.2cm} $O\left(\sqrt{\frac{L}{\mu(1-\sigma_2(\W))}}\log\frac{L}{\mu(1-\sigma_2(\W))}\log\frac{1}{\epsilon}\right)$\hspace*{-0.2cm}\\
\hline
Lower Bound & $O\left(\sqrt{\frac{L}{\mu}}\log\frac{1}{\epsilon}\right)$ \cite{nesterov2013introductory}& $O\left(\sqrt{\frac{L}{\mu(1-\sigma_2(\W))}}\log\frac{1}{\epsilon}\right)$ \cite{dasent}\\
\hline
\end{tabular}
\end{center}
\caption{Complexity comparisons between the accelerated dual ascent, accelerated penalty method with consensus, EXTRA and accelerated EXTRA for smooth convex problems.}\label{table1}
\end{table*}

\subsection{Paper Organization}
The rest of the paper is organized as follows. Section \ref{sec:algorithm} gives a sharper analysis on the original EXTRA, and Section \ref{sec:algorithm2} develops the accelerated EXTRA. Section \ref{sec:analysis} proves the complexities, and Section \ref{sec:exp} gives some numerical experiments. Finally, we conclude in Section \ref{sec:conclude}.

\section{Enhanced Results on EXTRA}\label{sec:algorithm}
We give a sharper analysis on EXTRA in this section. Specifically, section \ref{sec:extra} studies the strongly convex problems and section \ref{sec:extra2} studies the nonstrongly convex ones, respectively.
\subsection{Sharper Analysis on EXTRA for Strongly Convex Problems}\label{sec:extra}
We first describe EXTRA in the primal-dual form. From Assumption \ref{assumption_w} and the definition in (\ref{define_u}), we know that $\x$ is consensus if and only if $\U\x=\0$. Thus, we can reformulate problem (\ref{problem}) as the following linearly constrained problem:
\begin{eqnarray}
\min_{\x\in\R^{m\times n}} f(\x)\quad \mbox{s.t.}\quad \U\x=\0.\label{constrained_problem}
\end{eqnarray}
Introduce the augmented Lagrangian function
\begin{eqnarray}
L(\x,\blambda)=f(\x)+\<\blambda,\U\x\>+\frac{\beta}{2}\|\U\x\|_F^2.\notag
\end{eqnarray}
Problem (\ref{constrained_problem}) can be solved by the classical primal-dual method \cite{hong-2017,jakovetic-2017}. Specifically, it uses the Gauss$-$Seidel-like order to compute the saddle point of the augmented Lagrangian function and consists of the following iterations:
\begin{subequations}
\begin{align}
&\x^{k+1}=\x^k-\frac{1}{2(L+\beta)}(\nabla f(\x^k)+\U\blambda^k+\beta\U^2\x^k),\label{alm_s1}\\
&\blambda^{k+1}=\blambda^k+\beta\U\x^{k+1},\label{alm_s2}
\end{align}
\end{subequations}
where we specify the step-size in the primal step as $\frac{1}{2(L+\beta)}$. Step (\ref{alm_s2}) involves the operation of $\U\x$, which is uncomputable in the distributed environment. We introduce the auxiliary variable
\begin{equation}
\v^k=\U\blambda^k.\notag
\end{equation}
Multiplying both sides of (\ref{alm_s2}) by $\U$ it leads to
\begin{equation}
\v^{k+1}=\v^k+\beta\U^2\x^{k+1}.\notag
\end{equation}
From the definition of $\U$ in (\ref{define_u}), we have Algorithm \ref{extra}. Now, we establish the convergence of Algorithm \ref{extra}. Define
\begin{eqnarray}\label{define_rho}
\rho_k=(L+\beta)\|\x^k-\x^*\|_F^2+\frac{1}{2\beta}\|\blambda^k-\blambda^*\|_F^2,
\end{eqnarray}
where $(\x^*,\blambda^*)$ is a KKT point of the saddle point problem $\min_{\x}\max_{\blambda} f(\x)+\<\blambda,\U\x\>$. We prove the exponentially diminishing of $\rho_k$ in the following theorem. Especially, we choose a smaller $\beta$, i.e., a larger step-size $\alpha$, in the primal step than that in \cite{Shi-2015} to obtain a faster convergence rate. More precisely, the original EXTRA uses the step-size of $O\left(\frac{\mu}{L^2}\right)$ [Remark 4]\cite{Shi-2015} and an open problem was proposed in \cite{Shi-2015} on how to prove linear convergence under the larger step-size of $O\left(\frac{1}{L}\right)$. Our analysis addresses this open problem. We leave the proof in Section \ref{sec:analysis1} and describe the crucial tricks there.
\begin{theorem}\label{theorem_alm}
Suppose that Assumptions \ref{assumption_f} and \ref{assumption_w} hold with $\mu>0$. Let $\v^0\in\mbox{Span}(\U^2)$, $\alpha=\frac{1}{2(L+\beta)}$ and $\beta=L$. Then, for Algorithm \ref{extra}, we have
\begin{eqnarray}
\begin{aligned}\notag
\rho_{k+1}\leq\left(1-\delta\right)\rho_k,
\end{aligned}
\end{eqnarray}
where $\delta=\frac{1}{39(\frac{L}{\mu}+\frac{1}{1-\sigma_2(\W)})}$.
\end{theorem}

Based on Theorem \ref{theorem_alm}, we can give the following corollary, which proves that Algorithm \ref{extra} needs $O((\frac{L}{\mu}+\frac{1}{1-\sigma_2(\W)})\log\frac{L}{\epsilon(1-\sigma_2(\W))})$ iterations to find an $\epsilon$-optimal solution. Recall that $\alpha(\x)$ is the average of $x_{(1)},...,x_{(m)}$ defined in (\ref{def_avg_x}).

\begin{corollary}\label{sc_complexity_lemma}
Under the assumptions of Theorem \ref{theorem_alm} and letting $\v^0=\0$, Algorithm \ref{extra} needs
\begin{eqnarray}
\begin{aligned}\notag
O\left( \left(\frac{L}{\mu}+\frac{1}{1-\sigma_2(\W)}\right)\log\frac{LR_1+R_2/L}{\epsilon(1-\sigma_2(\W))} \right)
\end{aligned}
\end{eqnarray}
iterations to achieve an $\epsilon$-optimal solution $\x$ such that
\begin{eqnarray}
\begin{aligned}\notag
F(\alpha(\x))-F(x^*)\leq \epsilon\quad\mbox{and}\quad\frac{1}{m}\sum_{i=1}^m \left\|x_{(i)}-\alpha(\x)\right\|^2\leq \epsilon^2.
\end{aligned}
\end{eqnarray}
\end{corollary}

\subsection{Sharper Analysis on EXTRA for Non-strongly Convex Problems}\label{sec:extra2}
We study EXTRA for nonstrongly convex problems in this section. Specifically, we study the original EXTRA in Section \ref{sec:nsc1} and the regularized EXTRA in Section \ref{sec:nsc2}, respectively.
\subsubsection{Complexity for the Original EXTRA}\label{sec:nsc1}
The $O\left(\frac{1}{K}\right)$ convergence rate of EXTRA was well studied in \cite{Shi-2015,Shi-2015-2}. However, the authors of \cite{Shi-2015,Shi-2015-2} did not establish the explicit dependence on $1-\sigma_2(\W)$. In this section, we study the original EXTRA and give the $O(\frac{L}{K\sqrt{1-\sigma_2(\W)}})$ convergence rate in the following lemma.
\begin{lemma}\label{nsc_lemma}
Suppose that Assumptions \ref{assumption_f} and \ref{assumption_w} hold with $\mu=0$. Let $\alpha=\frac{1}{2(L+\beta)}$ and $\beta=\frac{L}{\sqrt{1-\sigma_2(\W)}}$, and define $\hat\x^K=\frac{1}{K}\sum_{k=1}^K\x^k$. Assume that $K\geq\frac{1}{\sqrt{1-\sigma_2(\W)}}$. Then, for Algorithm \ref{extra}, we have
\begin{eqnarray}
\begin{aligned}\notag
&F(\alpha(\hat\x^K))-F(x^*)\leq\frac{34}{K\sqrt{1-\sigma_2(\W)}}\left(LR_1+\frac{R_2}{L}\right),\\
&\frac{1}{m}\sum_{i=1}^m \left\|\hat x_{(i)}^K-\alpha(\hat\x^K)\right\|^2\leq \frac{16}{K^2(1-\sigma_2(\W))}\left(R_1+\frac{R_2}{L^2}\right).
\end{aligned}
\end{eqnarray}
\end{lemma}

We assume $K\geq\frac{1}{\sqrt{1-\sigma_2(\W)}}$ in Lemma \ref{nsc_lemma}, and it is a reasonable assumption. Take the linear network as an example, where all the agents connect in a line. For this special network, we know that $\frac{1}{1-\sigma_2(\W)}=m^2$ \cite{nedic-2018}. Algorithm \ref{extra} needs at least $m$ iterations to exchange messages between the two farthest nodes in the network. Thus, any convergent method needs at least $\frac{1}{\sqrt{1-\sigma_2(\W)}}$ iterations.

In Section \ref{sec:extra}, we establish the $O((\frac{L}{\mu}+\frac{1}{1-\sigma_2(\W)})\log\frac{L}{\epsilon(1-\sigma_2(\W))})$ complexity for strongly convex problems. Naturally, one may expect the $O(\frac{L}{\epsilon}+\frac{1}{1-\sigma_2(\W)})$ complexity for nonstrongly convex ones. However, Lemma \ref{nsc_lemma} only proves the $O(\frac{L}{\epsilon\sqrt{1-\sigma_2(\W)}})$ complexity. We describe the technical challenges in Section \ref{sec:analysis3}. It is currently unclear how to establish the faster rate for the original EXTRA, and we leave it as an open problem. In the following section, we improve the complexity via solving a regularized problem.

\subsubsection{Complexity for the Regularized EXTRA}\label{sec:nsc2}
When the complexity for the strongly convex problems is well studied, the regularization technique is a common way to solve the nonstrongly convex ones \cite{Allen-Zhu-2016-nips}. Namely, we add a small strongly convex regularizer to the objective and solve the regularized problem instead. Define the regularized version of $F(x)$ as
\begin{eqnarray}
\begin{aligned}\label{problem2}
F_{\epsilon}(x)=\frac{1}{m}\sum_{i=1}^mf_i(x)+\frac{\epsilon}{2}\|x\|^2
\end{aligned}
\end{eqnarray}
and denote $x_{\epsilon}^*=\argmin_{x} F_{\epsilon}(x)$. It can be easily checked that the precision between problems (\ref{problem}) and (\ref{problem2}) satisfies
\begin{eqnarray}
\begin{aligned}\label{regularized_relation}
F(x)-F(x^*)\leq F_{\epsilon}(x)-F_{\epsilon}(x_{\epsilon}^*)+\frac{\epsilon}{2}\|x^*\|^2.
\end{aligned}
\end{eqnarray}
Thus, to attain an $\epsilon$-optimal solution of problem (\ref{problem}), we only need to find an $\epsilon$-optimal solution of problem (\ref{problem2}). Denote $L_{\epsilon}=L+\epsilon$. Define
\begin{eqnarray}
\begin{aligned}\notag
f_{\epsilon}(\x)=f(\x)+\frac{\epsilon}{2}\|\x\|_F^2
\end{aligned}
\end{eqnarray}
and that is $L_{\epsilon}$-smooth and $\epsilon$-strongly convex. Problem (\ref{problem2}) can be reformulated as the following constrained problem:
\begin{eqnarray}
\begin{aligned}\label{reg_problem}
\min_{\x} f_{\epsilon}(\x)\quad \mbox{s.t.}\quad \U\x=0.
\end{aligned}
\end{eqnarray}
Denote $(\x_{\epsilon}^*,\blambda_{\epsilon}^*)$ to be a pair of KKT points of problem (\ref{reg_problem}). We use Algorithm \ref{extra} to solve problem (\ref{reg_problem}), and Corollary \ref{sc_complexity_lemma} needs $O((\frac{L}{\epsilon}+\frac{1}{1-\sigma_2(\W)})\log\frac{L}{\epsilon(1-\sigma_2(\W))})$ iterations to find an $\epsilon$-optimal solution of problem (\ref{problem2}), which is also an $\epsilon$-optimal solution of problem (\ref{problem}).

When $\epsilon\leq 1-\sigma_2(\W)$, the complexity of the above regularized EXTRA is dominated by $O(\frac{L}{\epsilon}\log\frac{L}{\epsilon})$. We want to further reduce the complexity by the $(\log\frac{L}{\epsilon})$ factor. As discussed in Section \ref{sec:analysis3}, the main reason for the slow rate of the original EXTRA discussed in Section \ref{sec:nsc1} is that $(L+\beta)\|\x^0-\x^*\|_F^2+\frac{1}{2\beta}\|\blambda^0-\blambda^*\|_F^2$ has the same order of magnitude as $O(1)$, rather than $O(1-\sigma_2(\W))$. Our motivation is that we may find a good enough initializer in a reasonable time such that $(L+\beta)\|\x^0-\x^*\|_F^2+\frac{1}{2\beta}\|\blambda^0-\blambda^*\|_F^2$ is of the order $O(1-\sigma_2(\W))$. With this inspiration, our algorithm consists of two stages. In the first stage, we run Algorithm \ref{extra} for $K_0$ iterations to solve problem (\ref{problem2}) and use its output $(\x^{K_0},\blambda^{K_0})$ as the initializer of the second stage. In the second stage, we run Algorithm \ref{extra} on problem (\ref{problem2}) again for $K$ iterations and output the averaged solution $\hat\x^K$. Although we analyze the method in two stages, we implement it in a single loop and only average over the last $K$ iterations.

The complexity of our two-stage regularized EXTRA is described in the next lemma. We see that the complexity is improved from $O((\frac{L}{\epsilon}+\frac{1}{1-\sigma_2(\W)})\log\frac{L}{\epsilon(1-\sigma_2(\W))})$ to $O((\frac{L}{\epsilon}+\frac{1}{1-\sigma_2(\W)})\log\frac{1}{1-\sigma_2(\W)})$ via the two-stage strategy.
\begin{lemma}\label{nsc_lemma2}
Suppose that Assumptions \ref{assumption_f} and \ref{assumption_w} hold with $\mu=0$. Let $\v^0\in\mbox{Span}(\U^2)$, $\alpha=\frac{1}{2(L_{\epsilon}+\beta)}$, and $\beta=L_{\epsilon}$. Run Algorithm \ref{extra} on problem (\ref{problem2}). Then, we only need
\begin{eqnarray}
\begin{aligned}\notag
O\left(\left(\frac{L}{\epsilon}+\frac{1}{1-\sigma_2(\W)}\right)\log\frac{1}{1-\sigma_2(\W)}\right)
\end{aligned}
\end{eqnarray}
iterations for the first stage and $K=O(\frac{LR_1+R_2/L}{\epsilon})$ iterations for the second stage such that
\begin{eqnarray}
\begin{aligned}\notag
F_{\epsilon}(\alpha(\hat\x^K))-F_{\epsilon}(x_{\epsilon}^*)\leq \epsilon\quad\mbox{and}\quad\frac{1}{m}\sum_{i=1}^m \left\|\hat x_{(i)}^K-\alpha(\hat\x^K)\right\|^2\leq \epsilon^2,
\end{aligned}
\end{eqnarray}
where $\hat\x^K=\frac{1}{K}\sum_{k=1}^K\x^k$ in the second stage.
\end{lemma}

\section{Accelerated EXTRA}\label{sec:algorithm2}
We first review Catalyst and then use it to accelerate EXTRA.
\subsection{Catalyst}\label{sec:catalyst}
Catalyst \cite{catalyst} is a general scheme for accelerating gradient-based optimization methods in the sense of Nesterov. It builds upon the inexact accelerated proximal point algorithm, which consists of the following iterations:
\begin{subequations}
\begin{align}
&x^{k+1}\approx \argmin_{x\in\R^n} F(x)+\frac{\tau}{2}\|x-y^k\|^2,\label{catalyst_s1}\\
&y^{k+1}=x^{k+1}+\frac{\theta_k(1-\theta_k)}{\theta_k^2+\theta_{k+1}}(x^{k+1}-x^k)\label{catalyst_s2},
\end{align}
\end{subequations}
where $\theta_k$ is defined in Algorithm \ref{acc-extra}. Catalyst employs double loop and approximately solves a sequence of well-chosen auxiliary problems in step (\ref{catalyst_s1}) in the inner loop. The following theorem describes the convergence rate for the outer loop.
\begin{theorem}\cite{schnidt-2011-nips,catalyst}\label{theorem_catalyst}
Suppose that $F(x)$ is convex and the following criterion holds for all $k\leq K$ with $\varepsilon_k\leq \frac{1}{k^{4+2\xi}}$:
\begin{eqnarray}
\begin{aligned}\label{inexact_proximal}
F(x^{k+1})+\frac{\tau}{2}\|x^{k+1}-y^k\|^2\leq \min_{x}\left(F(x)+\frac{\tau}{2}\|x-y^k\|^2\right)+\varepsilon_{k},
\end{aligned}
\end{eqnarray}
where $\xi$ can be any small positive constant. Then, Catalyst generates iterates $(x^k)_{k=0}^{K+1}$ such that
\begin{eqnarray}
\begin{aligned}\notag
F(x^{K+1})-F(x^*)\leq\frac{1}{(K+2)^2}\left(6\tau\|x^0-x^*\|^2+\frac{48}{\xi^2}+\frac{12}{1+2\xi}\right).
\end{aligned}
\end{eqnarray}
Suppose that $F(x)$ is $\mu$-strongly convex and (\ref{inexact_proximal}) holds for all $k\leq K$ with the precision of $\varepsilon_k\leq \frac{2(F(x^0)-F(x^*))}{9}(1-\rho)^{k+1}$, where $\rho<\sqrt{q}$ and $q=\frac{\mu}{\mu+\tau}$. Then, Catalyst generates iterates $(x_k)_{k=0}^{K+1}$ such that
\begin{eqnarray}
\begin{aligned}\notag
F(x^{K+1})-F(x^*)\leq \frac{8}{(\sqrt{q}-\rho)^2}(1-\rho)^{K+2}(F(x^0)-F(x^*)).
\end{aligned}
\end{eqnarray}
\end{theorem}

Briefly, Catalyst uses some linearly convergent method to solve the subproblem in step (\ref{catalyst_s1}) with warm-start, balances the outer loop and inner loop and attains the near optimal global complexities.

\subsection{Accelerating EXTRA via Catalyst}\label{sec:acc_extra}
We first establish the relation between Algorithm \ref{acc-extra} and Catalyst. Recall the definition of $G^k(x)$ in Algorithm \ref{acc-extra}:
\begin{eqnarray}
\begin{aligned}\notag
G^k(x)=\frac{1}{m}\sum_{i=1}^m g_i^k(x),\quad \mbox{where}\quad g_i^k(x)=f_i(x)+\frac{\tau}{2}\|x-y_{(i)}^k\|^2,
\end{aligned}
\end{eqnarray}
which is $(L+\tau)$-smooth and $(\mu+\tau)$-strongly convex. Denote $L_g=L+\tau$ and $\mu_g=\mu+\tau$ for simplicity. We can easily check that
\begin{eqnarray}
\begin{aligned}\notag
G^k(x)=&\frac{1}{m}\sum_{i=1}^m f_i(x)+\frac{\tau}{2}\sum_{i=1}^m\frac{1}{m}\|x-y_{(i)}^k\|^2\\
=&\frac{1}{m}\sum_{i=1}^m f_i(x)+\frac{\tau}{2}\|x-\alpha(\y^k)\|^2-\frac{\tau}{2}\|\alpha(\y^k)\|^2+\frac{\tau}{2}\sum_{i=1}^m\frac{1}{m}\|y_{(i)}^k\|^2.
\end{aligned}
\end{eqnarray}
Recall that $\alpha(\y^k)$ is the average of $y_{(1)}^k,...,y_{(m)}^k$ defined in (\ref{def_avg_x}). In Algorithm \ref{acc-extra}, we call EXTRA to minimize $G^k(x)$ approximately, i.e., to minimize $F(x)+\frac{\tau}{2}\|x-\alpha(\y^k)\|^2$. Thus, Algorithm \ref{acc-extra} can be interpreted as
\begin{eqnarray}
\begin{aligned}\notag
&\alpha(\x^{k+1})\approx\argmin_{x} F(x)+\frac{\tau}{2}\|x-\alpha(\y^k)\|^2,\\
&\alpha(\y^{k+1})=\alpha(\x^{k+1})+\frac{\theta_k(1-\theta_k)}{\theta_k^2+\theta_{k+1}}(\alpha(\x^{k+1})-\alpha(\x^k)),
\end{aligned}
\end{eqnarray}
and it belongs to the Catalyst framework. Thus, we only need to ensure (\ref{inexact_proximal}), i.e., $G^k(\alpha(\x^{k+1}))$ $\leq \min_x G^k(x)+\varepsilon_k$ for all $k$. Catalyst requires the liner convergence in the form of
\begin{equation}\notag
G^k(z^t)-\min_x G^k(x)\leq (1-\delta)^t\left(G^k(z^0)-\min_x G^k(x)\right)
\end{equation}
when solving the subproblem in step (\ref{catalyst_s1}), which is not satisfied for Algorithm \ref{extra} due to the existence of terms $\|\blambda^k-\blambda^*\|_F^2$ and $\|\blambda^{k+1}-\blambda^*\|_F^2$ in Theorem \ref{theorem_alm}. Thus, the conclusion in \cite{catalyst} cannot be directly applied to Algorithm \ref{acc-extra}. By analyzing the inner loop carefully, we can have the following theorem, which establishes that a suitable constant setup of $T_k$ is sufficient to ensure (\ref{inexact_proximal}) in the strongly convex case and thus allows us to use the Catalyst framework for distributed optimization, where $T_k$ is the number of inner iterations when calling Algorithm \ref{extra}.

\begin{theorem}\label{acc_extra_the1}
Suppose that Assumptions \ref{assumption_f} and \ref{assumption_w} hold with $\mu>0$. We only need to set $T_k=O((\frac{L+\tau}{\mu+\tau}+\frac{1}{1-\sigma_2(\W)})\log\frac{L+\tau}{\mu(1-\sigma_2(\W))})$ in Algorithm \ref{acc-extra} such that $G^k(\alpha(\x^{k+1}))\leq \min_x G^k(x)+\varepsilon_k$ holds for all $k$, where $\varepsilon_k$ is defined in Theorem \ref{theorem_catalyst}.
\end{theorem}

Based on Theorems \ref{theorem_catalyst} and \ref{acc_extra_the1}, we can establish the global complexity via finding the optimal balance between the inner loop and outer loop. Specifically, the total number of inner iterations is
\begin{eqnarray}
\begin{aligned}\notag
\sum_{k=0}^{\sqrt{1+\frac{\tau}{\mu}}\log\frac{1}{\epsilon}}T_k=\sqrt{1+\frac{\tau}{\mu}}\left(\frac{L+\tau}{\mu+\tau}+\frac{1}{1-\sigma_2(\W)}\right)\log\frac{L+\tau}{\mu(1-\sigma_2(\W))}\log\frac{1}{\epsilon}.
\end{aligned}
\end{eqnarray}
We obtain the minimal value with the optimal setting of $\tau$, which is described in the following corollary. On the other hand, when we set $\tau\approx0$, it approximates the original EXTRA.
\begin{corollary}\label{acc_extra_cor1}
Under the settings of Theorem \ref{acc_extra_the1} and letting $\tau=L(1-\sigma_2(\W))-\mu$, Algorithm \ref{acc-extra} needs
\begin{eqnarray}
\begin{aligned}\notag
O\left(\sqrt{\frac{L}{\mu(1-\sigma_2(\W))}}\log\frac{L}{\mu(1-\sigma_2(\W))}\log\frac{1}{\epsilon}\right)
\end{aligned}
\end{eqnarray}
total inner iterations to achieve an $\epsilon$-optimal solution such that
\begin{eqnarray}
\begin{aligned}\notag
F(\alpha(\x))-F(x^*)\leq \epsilon\quad\mbox{and}\quad\frac{1}{m}\sum_{i=1}^m \left\|x_{(i)}-\alpha(\x)\right\|^2\leq \epsilon^2.
\end{aligned}
\end{eqnarray}
\end{corollary}

When the strong convexity is absent, we have the following conclusions, which is the counterpart of Theorem \ref{acc_extra_the1} and Corollary \ref{acc_extra_cor1}.
\begin{theorem}\label{acc_extra_the2}
Suppose that Assumptions \ref{assumption_f} and \ref{assumption_w} hold with $\mu=0$. We only need to set $T_k=O((\frac{L+\tau}{\tau}+\frac{1}{1-\sigma_2(\W)})\log \frac{k}{1-\sigma_2(\W)})$ in Algorithm \ref{acc-extra} such that $G^k(\alpha(\x^{k+1}))\leq \min_x G^k(x)+\varepsilon_k$ holds for all $k$.
\end{theorem}
\begin{corollary}\label{acc_extra_cor2}
Under the settings of Theorem \ref{acc_extra_the2} and letting $\tau=L(1-\sigma_2(\W))$, Algorithm \ref{acc-extra} needs
\begin{eqnarray}
\begin{aligned}\notag
O\left(\sqrt{\frac{L}{\epsilon(1-\sigma_2(\W))}}\log\frac{1}{\epsilon(1-\sigma_2(\W))}\right)
\end{aligned}
\end{eqnarray}
total inner iterations to achieve an $\epsilon$-optimal solution such that
\begin{eqnarray}
\begin{aligned}\notag
F(\alpha(\x))-F(x^*)\leq \epsilon\quad\mbox{and}\quad\frac{1}{m}\sum_{i=1}^m \left\|x_{(i)}-\alpha(\x)\right\|^2\leq \epsilon^2.
\end{aligned}
\end{eqnarray}
\end{corollary}

The accelerated EXTRA needs to know $\frac{1}{1-\sigma_2(\W)}$ in advance to set $T_k$. Generally speaking, $\frac{1}{1-\sigma_2(\W)}$ relates to the global connectivity of the network. The authors of [Proposition 5]\cite{nedic-2018} give the estimation of $\frac{1}{1-\sigma_2(\W)}$ by $m$ for many frequently used networks, e.g., the 2-D graph, the geometric graph, the expander graph and the Erd\H{o}s$-$R\'{e}nyi random graph. See \cite{nedic-2018} for the details.
\section{Proof of Theorems}\label{sec:analysis}

In this section, we give the proofs of the theorems, corollaries, and lemmas in Sections \ref{sec:algorithm} and \ref{sec:algorithm2}. We first present several supporting lemmas, which will be used in our analysis.

\begin{lemma}\label{lemma04}
Assume that Assumption \ref{assumption_w} holds and then we have that $\|\Pi\x\|_F\leq\sqrt{\frac{2}{1-\sigma_2(\W)}}\|\U\x\|_F$.
\end{lemma} 

The proof is similar to that of [Lemma 5]\cite{li-2018-pm} and we omit the details.
\begin{lemma}\label{U_ineq}
Suppose that $\x^*$ is the optimal solution of problem (\ref{constrained_problem}). There exists $\blambda^*\in\mbox{Span}(\U)$ such that $(\x^*,\blambda^*)$ is a KKT point of the saddle point problem $\min_{\x}\max_{\blambda} f(\x)+\<\blambda,\U\x\>$. For $\blambda^*$ and any $\blambda\in\mbox{Span}(\U)$, we have $\|\blambda^*\|_F\leq\frac{\sqrt{2}\|\nabla f(\x^*)\|_F}{\sqrt{1-\sigma_2(\W)}}$ and $\|\U(\blambda-\blambda^*)\|_F^2\geq \frac{1-\sigma_2(\W)}{2}\|\blambda-\blambda^*\|_F^2$.
\end{lemma}

The existence of $\blambda^*\in\mbox{Span}(\U)$ was proved in [Lemma 3.1]\cite{Shi-2015} and $\|\blambda^*\|_F\leq\frac{\|\nabla f(\x^*)\|_F}{\widetilde\sigma_{\min}(\U)}$ was proved in [Theorem 2]\cite{Lam-2017}, where $\widetilde\sigma_{\min}(\U)$ denotes the smallest nonzero singular value of $\U$ and it is equal to $\sqrt{\frac{1-\sigma_2(\W)}{2}}$. The last inequality can be obtained from a similar induction to the proof of Lemma \ref{lemma04}, and we omit the details. From Lemma \ref{U_ineq}, we can see that when we study the dependence on $1-\sigma_2(\W)$, we should deal with $\|\blambda^*\|_F$ carefully. $\|\blambda^*\|_F$ cannot be regarded as a constant that can be dropped in the complexities.

\begin{lemma}
Assume that $f(\x)$ is $\mu$-strongly convex and $L$-smooth. Then, we have
\begin{eqnarray}
\begin{aligned}\label{sc_ineq}
\frac{\mu}{2}\|\x-\x^*\|_F^2\leq f(\x)-f(\x^*)+\<\blambda^*,\U\x\>\leq \frac{L}{2}\|\x-\x^*\|_F^2.
\end{aligned}
\end{eqnarray}
Assume that $f(\x)$ is convex and $L$-smooth. Then, we have
\begin{eqnarray}
\begin{aligned}\label{sc_ineq2}
\frac{1}{2L}\|\nabla f(\x)+\U\blambda^*\|_F^2\leq f(\x)-f(\x^*)+\<\blambda^*,\U\x\>.
\end{aligned}
\end{eqnarray}
\end{lemma}
\emph{Proof:} We can easily see that $f(\x)+\<\blambda^*,\U\x\>$ is $\mu$-strongly convex and $L$-smooth in $\x$. Since $\x^*=\argmin_{\x} f(\x)+\<\blambda^*,\U\x\>$ and $\U\x^*=\0$, we have (\ref{sc_ineq}). From $\nabla f(\x^*)+\U\blambda^*=\0$ and the smoothness of $f(\x)+\<\blambda^*,\U\x\>$ [Theorem 2.1.5]\cite{nesterov2013introductory}, we have (\ref{sc_ineq2}).
\hfill$\Box$

\begin{lemma}\label{lemma_cont1}
Suppose that Assumptions \ref{assumption_f} and \ref{assumption_w} hold with $\mu=0$. Assume that $f(\x)-f(\x^*)+\<\blambda^*,\U\x\>\leq\epsilon_1$ and $\|\U\x\|_F\leq\epsilon_2$. Then, we have
\begin{eqnarray}
\begin{aligned}\notag
F(\alpha(\x))-F(x^*)\leq \frac{1}{m}\left(\epsilon_1+\frac{3\|\nabla f(\x^*)\|_F+2L\|\x-\x^*\|_F}{\sqrt{1-\sigma_2(\W)}}\epsilon_2+\frac{L}{1-\sigma_2(\W)}\epsilon_2^2\right).
\end{aligned}
\end{eqnarray}
\end{lemma}
\emph{Proof:} Recall that $\frac{1}{m}\1\1^T\x=\1(\alpha(\x))^T$ from (\ref{def_avg_x}) and $\x^*=\1 (x^*)^T$. From the definitions of $F(x)$ and $f(\x)$ in (\ref{problem}) and (\ref{define_f}), respectively, we have $F(\alpha(\x))=\frac{1}{m}f\left(\frac{1}{m}\1\1^T\x\right)$ and $F(x^*)=\frac{1}{m}f(\x^*)$. Thus, we only need to bound $f(\frac{1}{m}\1\1^T\x)-f(\x^*)$.
\begin{eqnarray}
\begin{aligned}\notag
&f\left(\frac{1}{m}\1\1^T\x\right)-f(\x^*)\\
&=f\left(\frac{1}{m}\1\1^T\x\right)-f(\x)+f(\x)-f(\x^*)\\
&\overset{a}\leq \<\nabla f(\x),\frac{1}{m}\1\1^T\x-\x\>+\frac{L}{2}\|\Pi\x\|_F^2+f(\x)-f(\x^*)\\
&\overset{b}\leq \left(\|\nabla f(\x^*)\|_F+L\|\x-\x^*\|_F\right)\|\Pi\x\|_F+\frac{L}{2}\|\Pi\x\|_F^2+f(\x)-f(\x^*)\\
&\overset{c}\leq \left(\|\nabla f(\x^*)\|_F+L\|\x-\x^*\|_F\right)\sqrt{\frac{2}{1-\sigma_2(\W)}}\|\U\x\|_F+\frac{L}{1-\sigma_2(\W)}\|\U\x\|_F^2
\end{aligned}
\end{eqnarray}
\begin{eqnarray}
\begin{aligned}\notag
&+f(\x)-f(\x^*)+\<\blambda^*,\U\x\>+\|\blambda^*\|_F\|\U\x\|_F\\
\overset{d}\leq& \left(\|\nabla f(\x^*)\|_F+L\|\x-\x^*\|_F\right)\sqrt{\frac{2}{1-\sigma_2(\W)}}\|\U\x\|_F+\frac{L}{1-\sigma_2(\W)}\|\U\x\|_F^2\\
&+f(\x)-f(\x^*)+\<\blambda^*,\U\x\>+\frac{\sqrt{2}\|\nabla f(\x^*)\|_F}{\sqrt{1-\sigma_2(\W)}}\|\U\x\|_F,
\end{aligned}
\end{eqnarray}
where we use the smoothness of $f(\x)$ and (\ref{define_pi}) in $\overset{a}\leq$ and $\overset{b}\leq$, Lemma \ref{lemma04} and $-\<\blambda^*,\U\x\>\leq \|\blambda^*\|_F\|\U\x\|_F$ in $\overset{c}\leq$ and Lemma \ref{U_ineq} in $\overset{d}\leq$.
\hfill$\Box$

The following lemma is the well-known coerciveness property of the proximal operator.
\begin{lemma}[Lemma 22]\cite{catalyst}\label{lemma05}
Given a convex function $F(x)$ and a positive constant $\tau$, define $p(y)=\argmin_{x} F(x)+\frac{\tau}{2}\|x-y\|^2$. For any $y$ and $y'$, the following inequality holds,
\begin{eqnarray}
\begin{aligned}\notag
\|y-y'\|\geq \|p(y)-p(y')\|.
\end{aligned}
\end{eqnarray}
\end{lemma}

Finally, we study the regularized problem (\ref{reg_problem}).
\begin{lemma}\label{x_star_bound}
Suppose that Assumptions \ref{assumption_f} and \ref{assumption_w} hold with $\mu=0$. Then, we have $\|\x^*-\x_{\epsilon}^*\|_F\leq \|\x^*\|_F$ and $\|\x_{\epsilon}^*\|_F\leq 2\|\x^*\|_F$.
\end{lemma}

The proof is similar to that of [Claim 3.4]\cite{Allen-Zhu-2016-nips}. We omit the details.

\subsection{Proofs of Theorem \ref{theorem_alm} and Corollary \ref{sc_complexity_lemma}}\label{sec:analysis1}
Now, we prove Theorem \ref{theorem_alm}, which is based on the following lemma. It gives a progress in one iteration of Algorithm \ref{extra}. Some techniques in this proof have already appeared in \cite{Shi-2015}, and we present the proof for the sake of completeness.
\begin{lemma}
Suppose that Assumptions \ref{assumption_f} and \ref{assumption_w} hold with $\mu=0$. Then, for procedure (\ref{alm_s1})$-$(\ref{alm_s2}), we have
\begin{eqnarray}
\begin{aligned}\label{alm_ineq}
&f(\x^{k+1})-f(\x^*)+\<\blambda^*,\U\x^{k+1}\>\\
&\leq (L+\beta)\|\x^k-\x^*\|_F^2-(L+\beta)\|\x^{k+1}-\x^*\|_F^2\\
&\quad+\frac{1}{2\beta}\|\blambda^k-\blambda^*\|_F^2-\frac{1}{2\beta}\|\blambda^{k+1}-\blambda^*\|_F^2-\frac{\beta+L}{2}\|\x^{k+1}-\x^k\|_F^2.
\end{aligned}
\end{eqnarray}
\end{lemma}
\emph{Proof:} From the $L$-smoothness and convexity of $f(\x)$, we have
\begin{eqnarray}
\begin{aligned}\notag
f(\x^{k+1})\leq &f(\x^k)+\<\nabla f(\x^k),\x^{k+1}-\x^k\>+\frac{L}{2}\|\x^{k+1}-\x^k\|_F^2\\
= &f(\x^k)+\<\nabla f(\x^k),\x^*-\x^k\>+\<\nabla f(\x^k),\x^{k+1}-\x^*\>+\frac{L}{2}\|\x^{k+1}-\x^k\|_F^2\\
\leq &f(\x^*)+\<\nabla f(\x^k),\x^{k+1}-\x^*\>+\frac{L}{2}\|\x^{k+1}-\x^k\|_F^2.
\end{aligned}
\end{eqnarray}
Plugging (\ref{alm_s1}) into the above inequality, adding $\<\blambda^*,\U\x^{k+1}\>$ to both sides, and rearranging the terms, we have
\begin{eqnarray}
\begin{aligned}\notag
&f(\x^{k+1})-f(\x^*)+\<\blambda^*,\U\x^{k+1}\>\\
&\leq -\<2(L+\beta)(\x^{k+1}-\x^k)+\U\blambda^k+\beta\U^2\x^k,\x^{k+1}-\x^*\>\\
&\quad+\frac{L}{2}\|\x^{k+1}-\x^k\|_F^2+\<\blambda^*,\U\x^{k+1}\>\\
&\overset{a}= -2(L+\beta)\<\x^{k+1}-\x^k,\x^{k+1}-\x^*\>-\frac{1}{\beta}\<\blambda^k-\blambda^*,\blambda^{k+1}-\blambda^k\>\\
&\quad-\beta\<\U\x^k,\U\x^{k+1}\>+\frac{L}{2}\|\x^{k+1}-\x^k\|_F^2,
\end{aligned}
\end{eqnarray}
where we use $\U\x^*=0$ and (\ref{alm_s2}) in $\overset{a}=$. Using the identity of $2\<a,b\>=\|a\|^2+\|b\|^2-\|a-b\|^2$, we have
\begin{eqnarray}
\begin{aligned}\notag
&f(\x^{k+1})-f(\x^*)+\<\blambda^*,\U\x^{k+1}\>\\
&\leq (L+\beta)\|\x^k-\x^*\|_F^2-(L+\beta)\|\x^{k+1}-\x^*\|_F^2\\
&\quad+\frac{1}{2\beta}\|\blambda^k-\blambda^*\|_F^2-\frac{1}{2\beta}\|\blambda^{k+1}-\blambda^*\|_F^2+\frac{1}{2\beta}\|\blambda^{k+1}-\blambda^k\|_F^2\\
&\quad-\frac{\beta}{2}\|\U\x^k\|_F^2-\frac{\beta}{2}\|\U\x^{k+1}\|_F^2+\frac{\beta}{2}\|\U\x^{k+1}-\U\x^k\|_F^2-\left(\frac{L}{2}+\beta\right)\|\x^{k+1}-\x^k\|_F^2\\
&\overset{b}\leq (L+\beta)\|\x^k-\x^*\|_F^2-(L+\beta)\|\x^{k+1}-\x^*\|_F^2\\
&\quad+\frac{1}{2\beta}\|\blambda^k-\blambda^*\|_F^2-\frac{1}{2\beta}\|\blambda^{k+1}-\blambda^*\|_F^2-\frac{\beta+L}{2}\|\x^{k+1}-\x^k\|_F^2,
\end{aligned}
\end{eqnarray}
where we use (\ref{alm_s2}) and $\|\U\|_2^2\leq 1$ in $\overset{b}\leq$.\hfill$\Box$

A crucial property in (\ref{alm_ineq}) is that we keep the term $-\frac{\beta+L}{2}\|\x^{k+1}-\x^k\|_F^2$, which will be used in the following proof to eliminate the term $\left(\frac{1}{\nu}-1\right)\frac{9(L+\beta)^2}{2L}\|\x^{k+1}-\x^k\|_F^2$ to attain (\ref{sm2_ineq}). In the following proof of Theorem \ref{theorem_alm}, we use the strong convexity and smoothness of $f(\x)$ to obtain two inequalities, i.e., (\ref{sm2_ineq}) and (\ref{sc_ineq33}). A convex combination leads to (\ref{cont4}). The key thing here is to design the parameters carefully. Otherwise, we may only obtain a suboptimal result with a worse dependence on $\frac{L}{\mu}$ and $\frac{1}{1-\sigma_2(\W)}$.

\noindent\noindent\emph{Proof of Theorem \ref{theorem_alm}:}
We use (\ref{sc_ineq2}) to upper bound $\|\blambda^{k+1}-\blambda^*\|_F^2$. From procedure (\ref{alm_s1})$-$(\ref{alm_s2}), we have
\begin{eqnarray}
\begin{aligned}\notag
&2(L+\beta)\left(\x^{k+1}-\x^k\right)+\nabla  f(\x^k)+\U\blambda^{k+1}+\beta\U^2(\x^k-\x^{k+1})=0.
\end{aligned}
\end{eqnarray}
Thus, we obtain
\begin{eqnarray}
\begin{aligned}\label{sc_sm_inequ}
&\frac{1}{2L}\|\nabla f(\x^{k+1})+\U\blambda^*\|_F^2\\
&=\frac{1}{2L}\left\| 2(L+\beta)\left(\x^{k+1}-\x^k\right)+\beta\U^2(\x^k-\x^{k+1})+\nabla f(\x^k)-\nabla f(\x^{k+1})+\U(\blambda^{k+1}-\blambda^*) \right\|_F^2\\
&\overset{c}\geq\frac{1-\nu}{2L}\|\U(\blambda^{k+1}-\blambda^*)\|_F^2\\
&\quad-\frac{1/\nu-1}{2L}\left\| 2(L+\beta)\left(\x^{k+1}-\x^k\right)+\beta\U^2(\x^k-\x^{k+1})+\nabla f(\x^k)-\nabla f(\x^{k+1})\right\|_F^2\\
&\overset{d}\geq\frac{(1-\nu)(1-\sigma_2(\W))}{4L}\|\blambda^{k+1}-\blambda^*\|_F^2-\left(\frac{1}{\nu}-1\right)\frac{9(L+\beta)^2}{2L}\|\x^{k+1}-\x^k\|_F^2,
\end{aligned}
\end{eqnarray} 
where we use $\|a+b\|^2\geq (1-\nu)\|a\|^2-(1/\nu-1)\|b\|^2$ for some $\nu\in(0,1)$ in $\overset{c}\geq$, Lemma \ref{U_ineq} and the smoothness of $f(\x)$ in $\overset{d}\geq$. Lemma \ref{U_ineq} requires $\blambda^k\in\mbox{Span}(\U)$. From the initialization and (\ref{alm_s2}), we know it holds for all $k$.

Letting $\nu=\frac{9(\beta+L)}{9(\beta+L)+L}$, then $\left(\frac{1}{\nu}-1\right)\frac{9(L+\beta)^2}{2L}=\frac{L+\beta}{2}$. Plugging $\nu$ into the above inequality and using (\ref{sc_ineq2}) and (\ref{alm_ineq}), we have
\begin{eqnarray}
\begin{aligned}\label{sm2_ineq}
\hspace*{1cm}\frac{1-\sigma_2(\W)}{36(\beta+L)+4L}\|\blambda^{k+1}-\blambda^*\|_F^2\leq& (L+\beta)\|\x^k-\x^*\|_F^2-(L+\beta)\|\x^{k+1}-\x^*\|_F^2\\
\hspace*{1cm}&+\frac{1}{2\beta}\|\blambda^k-\blambda^*\|_F^2-\frac{1}{2\beta}\|\blambda^{k+1}-\blambda^*\|_F^2.
\end{aligned}
\end{eqnarray}
From (\ref{sc_ineq}) and (\ref{alm_ineq}), we also have
\begin{eqnarray}
\begin{aligned}\label{sc_ineq33}
\frac{\mu}{2}\|\x^{k+1}-\x^*\|_F^2\leq& (L+\beta)\|\x^k-\x^*\|_F^2-(L+\beta)\|\x^{k+1}-\x^*\|_F^2\\
&+\frac{1}{2\beta}\|\blambda^k-\blambda^*\|_F^2-\frac{1}{2\beta}\|\blambda^{k+1}-\blambda^*\|_F^2.
\end{aligned}
\end{eqnarray}
Multiplying (\ref{sm2_ineq}) by $\eta$, multiplying (\ref{sc_ineq33}) by $1-\eta$, adding them together and rearranging the terms, we have
\begin{eqnarray}
\begin{aligned}\label{cont4}
&\left(L+\beta+\frac{(1-\eta)\mu}{2}\right)\|\x^{k+1}-\x^*\|_F^2+\left(\frac{1}{2\beta}+\frac{\eta(1-\sigma_2(\W))}{36(\beta+L)+4L}\right)\|\blambda^{k+1}-\blambda^*\|_F^2\\
&\leq (L+\beta)\|\x^k-\x^*\|_F^2+\frac{1}{2\beta}\|\blambda^k-\blambda^*\|_F^2.
\end{aligned}
\end{eqnarray}
Letting $\frac{(1-\eta)\mu}{2(L+\beta)}=\frac{\beta \eta(1-\sigma_2(\W))}{18(\beta+L)+2L}$, we have $\eta=\frac{\frac{\mu}{2(L+\beta)}}{\frac{\mu}{2(L+\beta)}+\frac{\beta(1-\sigma_2(\W))}{18(\beta+L)+2L}}$. Plugging it into (\ref{cont4}) and recalling the definition of $\rho_k$ in (\ref{define_rho}), it leads to
\begin{eqnarray}
\begin{aligned}\notag
\left( 1+\frac{\mu\beta(1-\sigma_2(\W))}{\mu(18(\beta+L)+2L)+2(L+\beta)\beta(1-\sigma_2(\W))} \right)\rho_{k+1}\leq \rho_k.
\end{aligned}
\end{eqnarray}
We can easily check that
\begin{eqnarray}
\begin{aligned}\notag
&\frac{\mu\beta(1-\sigma_2(\W))}{\mu(18(\beta+L)+2L)+2(L+\beta)\beta(1-\sigma_2(\W))}\\
&=\frac{\mu(1-\sigma_2(\W))}{ \frac{20L\mu}{\beta}+2\beta(1-\sigma_2(\W))+2L(1-\sigma_2(\W))+18\mu }\\
&\overset{e}\geq\frac{1}{38}\frac{\mu(1-\sigma_2(\W))}{ L(1-\sigma_2(\W))+\mu }
\end{aligned}
\end{eqnarray}
by letting $\beta=L$ in $\overset{e}\geq$. Thus, we have the conclusion.\hfill$\Box$

Finally, we prove Corollary \ref{sc_complexity_lemma}.

\emph{Proof:} From Theorem \ref{theorem_alm}, $\blambda^0=\0$, (\ref{fact}), $\beta=L$ and Lemma \ref{U_ineq}, we have
\begin{eqnarray}
\begin{aligned}\label{cont005}
(L+\beta)\|\x^k-\x^*\|_F^2+\frac{1}{2\beta}\|\blambda^k-\blambda^*\|_F^2\leq& \left(1-\delta\right)^k\left((L+\beta)\|\x^0-\x^*\|_F^2+\frac{1}{2\beta}\|\blambda^*\|_F^2\right)\\
\leq& \left(1-\delta\right)^k\frac{2m(LR_1+R_2/L)}{1-\sigma_2(\W)}.
\end{aligned}
\end{eqnarray}
On the other hand, from $\x^*=\1(x^*)^T$, the definition of $\alpha(\x)$ in (\ref{def_avg_x}), the convexity of $\|\cdot\|^2$, and the smoothness of $F(x)$, we have
\begin{eqnarray}
\begin{aligned}\label{sc_cont2}
\qquad\frac{1}{m}\|\x^k-\x^*\|_F^2=\frac{1}{m}\sum_{i=1}^m\|x_{(i)}^k-x^*\|^2\geq \|\alpha(\x^k)-x^*\|^2\geq \frac{2}{L}\left(F(\alpha(\x^k))-F(x^*)\right).
\end{aligned}
\end{eqnarray}
So we have
\begin{eqnarray}
\begin{aligned}\notag
F(\alpha(\x^k))-F(x^*)\leq& \left(1-\delta\right)^k\frac{LR_1+R_2/L}{1-\sigma_2(\W)}.
\end{aligned}
\end{eqnarray}
On the other hand, since $\frac{1}{m}\sum_{i=1}^m \left\|x_{(i)}-\alpha(\x)\right\|^2=\frac{1}{m}\|\Pi\x\|_F^2$, we only need to bound $\|\Pi\x\|_F^2$:
\begin{eqnarray}
\begin{aligned}\label{sc_cont1}
\|\Pi\x^k\|_F^2\overset{a}\leq&\frac{2}{1-\sigma_2(\W)}\|\U\x^k\|_F^2\overset{b}=\frac{2}{(1-\sigma_2(\W))\beta^2}\|\blambda^k-\blambda^{k-1}\|_F^2\\
\leq& \frac{4}{(1-\sigma_2(\W))\beta^2}\left(\|\blambda^k-\blambda^*\|_F^2+\|\blambda^{k-1}-\blambda^*\|_F^2\right)\\
\overset{c}\leq& \frac{16}{(1-\sigma_2(\W))\beta}\left((L+\beta)\|\x^{k-1}-\x^*\|_F^2+\frac{1}{2\beta}\|\blambda^{k-1}-\blambda^*\|_F^2\right)\\
\overset{d}\leq& \left(1-\delta\right)^{k-1}\frac{32m(R_1+R_2/L^2)}{(1-\sigma_2(\W))^2},
\end{aligned}
\end{eqnarray}
where we use Lemma \ref{lemma04} in $\overset{a}\leq$, (\ref{alm_s2}) in $\overset{b}=$, $\|\blambda^k-\blambda^*\|_F^2\leq 2\beta\rho_k\leq 2\beta\rho_{k-1}$ and $\|\blambda^{k-1}-\blambda^*\|_F^2\leq 2\beta\rho_{k-1}$ in $\overset{c}\leq$, (\ref{cont005}) and $\beta=L$ in $\overset{d}\leq$. The proof is complete.
\hfill$\Box$

\subsection{Proofs of Lemmas \ref{nsc_lemma} and \ref{nsc_lemma2}}\label{sec:analysis3}
Lemma \ref{nsc_lemma} only proves the $O(\frac{L}{K\sqrt{1-\sigma_2(\W)}})$ convergence rate, rather than $O(\frac{L}{K})$. In fact, from Lemma \ref{lemma_cont1}, to prove the $O(\frac{1}{K})$ convergence rate, we should establish $\|\U\x\|_F\leq O(\frac{\sqrt{m(1-\sigma_2(\W))}}{K})$. However, from (\ref{nsc_cont5}), we know $\|\U\x\|_F$ has only the same order of magnitude as $O(\frac{\sqrt{m}}{K}\sqrt{R_1+\frac{R_2}{\beta^2(1-\sigma_2(\W))}})$. We find that $\beta=\frac{L}{\sqrt{1-\sigma_2(\W)}}$ is the best choice to balance the terms in (\ref{nsc_cont9}).

\noindent\emph{Proof of Lemma \ref{nsc_lemma}:} Summing (\ref{alm_ineq}) over $k=0,1,...,K-1$, dividing both sides by $K$, using the convexity of $f(\x)$, and using the definition of $\hat\x^K$, we have
\begin{eqnarray}
\begin{aligned}\label{nsc_cont3}
f(\hat\x^K)-f(\x^*)+\<\blambda^*,\U\hat\x^K\>\leq \frac{1}{K}\left((L+\beta)\|\x^0-\x^*\|_F^2+\frac{1}{2\beta}\|\blambda^0-\blambda^*\|_F^2\right).
\end{aligned}
\end{eqnarray}
On the other hand, since $f(\x)-f(\x^*)+\<\blambda^*,\U\x\>\geq 0$ for all $\x$ from (\ref{sc_ineq2}), we also have
\begin{eqnarray}
\begin{aligned}\label{nsc_cont4}
&\hspace*{0.8cm}(L+\beta)\|\x^k-\x^*\|_F^2\leq (L+\beta)\|\x^0-\x^*\|_F^2+\frac{1}{2\beta}\|\blambda^0-\blambda^*\|_F^2\quad\forall k=1,...,K
\end{aligned}
\end{eqnarray}
and
\begin{eqnarray}
\begin{aligned}\notag
&\frac{1}{2\beta}\|\blambda^K-\blambda^*\|_F^2\leq (L+\beta)\|\x^0-\x^*\|_F^2+\frac{1}{2\beta}\|\blambda^0-\blambda^*\|_F^2
\end{aligned}
\end{eqnarray}
from (\ref{alm_ineq}). Using (\ref{alm_s2}) and the definition of $\hat\x^K$, we further have
\begin{eqnarray}
\begin{aligned}\label{nsc_cont5}
\|\U\hat\x^K\|_F^2=&\frac{1}{\beta^2K^2}\|\blambda^K-\blambda^0\|_F^2\\
\leq&\frac{2}{\beta^2K^2}\|\blambda^K-\blambda^*\|_F^2+\frac{2}{\beta^2K^2}\|\blambda^0-\blambda^*\|_F^2\\
\leq& \frac{4}{K^2}\left(\frac{L+\beta}{\beta}\|\x^0-\x^*\|_F^2+\frac{1}{\beta^2}\|\blambda^0-\blambda^*\|_F^2\right).
\end{aligned}
\end{eqnarray}
Summing (\ref{nsc_cont4}) over $k=1,2,\cdots,K$, dividing both sides by $K$, using the convexity of $\|\cdot\|_F^2$, and using the definition of $\hat\x^K$, we have
\begin{eqnarray}
\begin{aligned}\label{nsc_cont6}
\qquad &(L+\beta)\|\hat\x^K-\x^*\|_F^2\leq (L+\beta)\|\x^0-\x^*\|_F^2+\frac{1}{2\beta}\|\blambda^0-\blambda^*\|_F^2\quad\forall k=1,...,K.
\end{aligned}
\end{eqnarray}
From (\ref{nsc_cont3}), (\ref{nsc_cont5}), (\ref{nsc_cont6}), and Lemma \ref{lemma_cont1}, we have
\begin{eqnarray}
\begin{aligned}\label{nsc_cont9}
F(\alpha(\hat\x^K))\hspace*{-0.05cm}-\hspace*{-0.05cm}F(x^*)\hspace*{-0.05cm}\leq\hspace*{-0.05cm}& \frac{1}{mK}\hspace*{-0.05cm}\left(\hspace*{-0.05cm}\left(\hspace*{-0.05cm}1\hspace*{-0.05cm}+\hspace*{-0.05cm}\frac{8L}{K\beta(1\hspace*{-0.05cm}-\hspace*{-0.05cm}\sigma_2(\W))}\hspace*{-0.05cm}\right)\hspace*{-0.05cm}\left(\hspace*{-0.05cm}(L\hspace*{-0.05cm}+\hspace*{-0.05cm}\beta)\|\x^0\hspace*{-0.05cm}-\hspace*{-0.05cm}\x^*\|_F^2\hspace*{-0.05cm}+\hspace*{-0.05cm}\frac{1}{2\beta}\|\blambda^*\|_F^2\hspace*{-0.05cm}\right)\right.\\
&\hspace*{1cm}\left.+\frac{6\|\nabla f(\x^*)\|_F}{\sqrt{1-\sigma_2(\W)}}\left(\sqrt{\frac{L+\beta}{\beta}\|\x^0-\x^*\|_F^2+\frac{1}{\beta^2}\|\blambda^*\|_F^2}\right)\right.\\
&\hspace*{1cm}\left.+\frac{4L}{\sqrt{1\hspace*{-0.05cm}-\hspace*{-0.05cm}\sigma_2(\W)}}\sqrt{\frac{L\hspace*{-0.05cm}+\hspace*{-0.05cm}\beta}{\beta}}\hspace*{-0.05cm}\left(\hspace*{-0.05cm}\|\x^0\hspace*{-0.05cm}-\hspace*{-0.05cm}\x^*\|_F^2\hspace*{-0.05cm}+\hspace*{-0.05cm}\frac{1}{\beta(L\hspace*{-0.05cm}+\hspace*{-0.05cm}\beta)}\|\blambda^*\|_F^2\hspace*{-0.05cm}\right)\hspace*{-0.05cm}\right).
\end{aligned}
\end{eqnarray}
Plugging $\|\lambda^*\|_F^2\leq\frac{2\|\nabla f(\x^*)\|_F^2}{1-\sigma_2(\W)}$, (\ref{fact}), and the setting of $\beta$ into the above inequality, after some simple computations, we have the first conclusion. Similarly, from (\ref{nsc_cont5}) and Lemma \ref{lemma04}, we have the second conclusion.
\hfill$\Box$

\noindent\emph{Proof of Lemma \ref{nsc_lemma2}:} For the first stage, from a modification of (\ref{cont005}) on problem (\ref{reg_problem}), we know that Algorithm \ref{extra} needs
\begin{eqnarray}
\begin{aligned}\notag
K_0=O\left(\left(\frac{L_{\epsilon}}{\epsilon}+\frac{1}{1-\sigma_2(\W)}\right)\log\frac{1}{1-\sigma_2(\W)}\right)
\end{aligned}
\end{eqnarray}
iterations such that
\begin{eqnarray}
\begin{aligned}\label{cont007}
(L_{\epsilon}+\beta)\|\x^{K_0}-\x_{\epsilon}^*\|_F^2+\frac{1}{2\beta}\|\blambda^{K_0}-\blambda_{\epsilon}^*\|_F^2\leq m(1-\sigma_2(\W))\left(L_{\epsilon}R_1+R_2/L_{\epsilon}\right).
\end{aligned}
\end{eqnarray}
Let $(\x^{K_0},\blambda^{K_0})$ be the initialization of the second stage. From a modification of (\ref{nsc_cont9}) on problem (\ref{problem2}), we have
\begin{eqnarray}
\begin{aligned}\notag
F_{\epsilon}(\alpha(\hat\x^K))-F_{\epsilon}(x_{\epsilon}^*)\leq& \frac{1}{mK}\left(\left(1+\frac{8L_{\epsilon}}{K\beta(1-\sigma_2(\W))}\right)m(1-\sigma_2(\W))(L_{\epsilon}R_1+R_2/L_{\epsilon})\right.\\
&\hspace*{1cm}\left.+\frac{6\|\nabla f_{\epsilon}(\x_{\epsilon}^*)\|_F}{\sqrt{1-\sigma_2(\W)}}\sqrt{\frac{2m(1-\sigma_2(\W))(L_{\epsilon}R_1+R_2/L_{\epsilon})}{\beta}}\right.\\
&\hspace*{1cm}\left.+\frac{4L_{\epsilon}}{\sqrt{1\hspace*{-0.06cm}-\hspace*{-0.06cm}\sigma_2(\W)}}\hspace*{-0.06cm}\sqrt{\hspace*{-0.06cm}\frac{L_{\epsilon}\hspace*{-0.06cm}+\hspace*{-0.06cm}\beta}{\beta}}\frac{2m(1\hspace*{-0.06cm}-\hspace*{-0.06cm}\sigma_2(\W)\hspace*{-0.03cm})(L_{\epsilon}R_1\hspace*{-0.06cm}+\hspace*{-0.06cm}R_2/L_{\epsilon})}{L_{\epsilon}+\beta}\hspace*{-0.06cm}\right)\hspace*{-0.06cm}.
\end{aligned}
\end{eqnarray}
From the definition of $f_{\epsilon}(\x)$, the smoothness of $f(\x)$, Lemma \ref{x_star_bound} and (\ref{fact}), we have $\|\nabla f_{\epsilon}(\x_{\epsilon}^*)\|_F$ $\leq \|\nabla f(\x^*)\|_F\hspace*{-0.03cm}+\hspace*{-0.03cm}L\|\x_{\epsilon}^*\hspace*{-0.03cm}-\hspace*{-0.03cm}\x^*\|_F\hspace*{-0.03cm}+\hspace*{-0.03cm}\epsilon\|\x_{\epsilon}^*\|_F\hspace*{-0.03cm}\leq\hspace*{-0.03cm}\sqrt{mR_2}\hspace*{-0.03cm}+\hspace*{-0.03cm}2L_{\epsilon}\sqrt{mR_1}\hspace*{-0.03cm}\leq\hspace*{-0.03cm}\sqrt{8mL_{\epsilon}(L_{\epsilon}R_1\hspace*{-0.03cm}+\hspace*{-0.03cm}R_2/L_{\epsilon})}$. From $\beta=L_{\epsilon}$ and after some simple calculations, we have
\begin{eqnarray}
\begin{aligned}\notag
F_{\epsilon}(\alpha(\hat\x^K))-F_{\epsilon}(x_{\epsilon}^*)\leq\frac{41(L_{\epsilon}R_1+R_2/L_{\epsilon})}{K}.
\end{aligned}
\end{eqnarray}
On the other hand, from Lemma \ref{lemma04}, (\ref{nsc_cont5}), (\ref{cont007}), and $\beta=L_{\epsilon}$, we have
\begin{eqnarray}
\begin{aligned}\notag
\|\Pi\hat\x^K\|_F^2\leq\frac{1}{1-\sigma_2(\W)}\|\U\hat\x^K\|_F^2\leq \frac{8m(R_1+R_2/L_{\epsilon}^2)}{K^2}.
\end{aligned}
\end{eqnarray}
Thus, the second stage needs $K=O(\frac{L_{\epsilon}R_1+R_2/L_{\epsilon}}{\epsilon})$ iterations such that $F_{\epsilon}(\alpha(\hat\x^K))-F_{\epsilon}(x_{\epsilon}^*)\leq\epsilon$ and $\frac{1}{m}\sum_{i=1}^m \left\|\hat x_{(i)}^K-\alpha(\hat\x^K)\right\|^2\leq \epsilon^2$.
\hfill$\Box$

\subsection{Proofs of Theorems \ref{acc_extra_the1} and \ref{acc_extra_the2}}\label{sec:analysis2}
We consider the strongly convex problems in Section \ref{sec:sc} and the nonstrongly convex ones in Section \ref{sec:ns}, respectively.
\subsubsection{Strongly Convex Case}\label{sec:sc}
In this section, we prove Theorem \ref{acc_extra_the1}. Define
\begin{eqnarray}
\begin{aligned}\notag
x^{k,*}=\argmin_{x} G^k(x)=\argmin_{x} F(x)+\frac{\tau}{2}\|x-\alpha(\y^k)\|^2
\end{aligned}
\end{eqnarray}
and denote $(\x^{k,*},\blambda^{k,*})$ to be a KKT point of saddle point problem $\min_{\x}\max_{\blambda}g^k(\x)+\<\blambda,\U\x\>$, where $g^k(\x)\equiv f(\x)+\frac{\tau}{2}\|\x-\y^k\|_F^2$. Then, we know $\x^{k,*}=\1(x^{k,*})^T$. Let $(\x^{k,t},\U\blambda^{k,t})_{t=0}^{T_k+1}$ be the iterates generated by Algorithm \ref{extra} at the $k$th iteration of Algorithm \ref{acc-extra}. Then, $\x^{k,0}=\x^k$ and $\x^{k,T_k+1}=\x^{k+1}$. Define
\begin{eqnarray}
\begin{aligned}\notag
\rho_{k,t}=\left(L_g+\beta_g\right)\|\x^{k,t}-\x^{k,*}\|_F^2+\frac{1}{2\beta_g}\|\blambda^{k,t}-\blambda^{k,*}\|_F^2,
\end{aligned}
\end{eqnarray}
where we set $\beta_g=L_g$. Similar to (\ref{sc_cont2}), we have
\begin{eqnarray}
\begin{aligned}\notag
G^k(\alpha(\x^{k+1}))-G^k(x^{k,*})=G^k(\alpha(\x^{k,T_k+1}))-G^k(x^{k,*})\leq \frac{1}{2m}\rho_{k,T_k}.
\end{aligned}
\end{eqnarray}
Thus, we only need to prove $\rho_{k,T_k}\leq 2m\varepsilon_k$. Moreover, we prove a sharper result of $\rho_{k,T_k}\leq 2m(1-\sigma_2(\W))\varepsilon_k$ by induction in the following lemma. The reason is that we want to prove $\|\Pi\x^{K+1}\|_F^2\leq O\left(m\varepsilon_K\right)$ and thus we need to eliminate $1-\sigma_2(\W)$ in (\ref{cont2}).
\begin{lemma}\label{lemma_cont2}
Suppose that Assumptions \ref{assumption_f} and \ref{assumption_w} hold with $\mu>0$. If $\rho_{r,T_r}\leq 2m(1-\sigma_2(\W))\varepsilon_r$ holds for all $r\leq k-1$ and we initialize $\x^{k,0}=\x^{k-1,T_{k-1}+1}$ and $\blambda^{k,0}=\blambda^{k-1,T_{k-1}+1}$, then we only need $T_k=O((\frac{L_g}{\mu_g}+\frac{1}{1-\sigma_2(\W)})\log\frac{L_g}{\mu(1-\sigma_2(\W))})$ such that $\rho_{k,T_k}\leq 2m(1-\sigma_2(\W))\varepsilon_k$.
\end{lemma}
\emph{Proof:}
From Theorem \ref{theorem_alm} and (\ref{sc_cont1}), we have
\begin{eqnarray}
&&\rho_{k,T_k}\leq\left(1-\delta_g\right)^{T_k}\rho_{k,0},\label{cont1}\\
&&\|\Pi\x^{k+1}\|_F^2=\|\Pi\x^{k,T_k+1}\|_F^2\leq \frac{16}{\beta_g(1-\sigma_2(\W))}\rho_{k,T_k}\label{cont2},
\end{eqnarray}
where $\delta_g=\frac{1}{39(\frac{L_g}{\mu_g}+\frac{1}{1-\sigma_2(\W)})}$. From the initialization and Theorem \ref{theorem_alm}, we have
\begin{eqnarray}
\begin{aligned}\label{cont001}
\rho_{k,0}=&\left(L_g+\beta_g\right)\|\x^{k-1,T_{k-1}+1}-\x^{k,*}\|_F^2+\frac{1}{2\beta_g}\|\blambda^{k-1,T_{k-1}+1}-\blambda^{k,*}\|_F^2\\
\leq& 2\left(L_g+\beta_g\right)\|\x^{k-1,T_{k-1}+1}-\x^{k-1,*}\|_F^2+\frac{1}{\beta_g}\|\blambda^{k-1,T_{k-1}+1}-\blambda^{k-1,*}\|_F^2\\
&+2\left(L_g+\beta_g\right)\|\x^{k,*}-\x^{k-1,*}\|_F^2+\frac{1}{\beta_g}\|\blambda^{k,*}-\blambda^{k-1,*}\|_F^2\\
\leq&2\rho_{k-1,T_{k-1}}+2\left(L_g+\beta_g\right)\|\x^{k,*}-\x^{k-1,*}\|_F^2+\frac{1}{\beta_g}\|\blambda^{k,*}-\blambda^{k-1,*}\|_F^2.
\end{aligned}
\end{eqnarray}
From the fact that $\x^{k,*}=\1 (x^{k,*})^T$, we have
\begin{eqnarray}
\begin{aligned}\label{cont002}
\|\x^{k,*}-\x^{k-1,*}\|_F^2=&m \|x^{k,*}-x^{k-1,*}\|^2\overset{a}\leq m\|\alpha(\y^k)-\alpha(\y^{k-1})\|^2\\
\overset{b}\leq& \sum_{i=1}^m\|y_{(i)}^k-y_{(i)}^{k-1}\|^2= \|\y^k-\y^{k-1}\|_F^2,
\end{aligned}
\end{eqnarray}
where $\overset{a}\leq$ uses Lemma \ref{lemma05} and $\overset{b}\leq$ uses the definition of $\alpha(\y)$ and the convexity of $\|\cdot\|^2$. From Lemma \ref{U_ineq}, we know
\begin{eqnarray}
\|\blambda^{k,*}-\blambda^{k-1,*}\|_F^2\leq\frac{2}{1-\sigma_2(\W)}\|\U\blambda^{k,*}-\U\blambda^{k-1,*}\|_F^2. \label{cont003}
\end{eqnarray}
Recall that $(\x^{k,*},\blambda^{k,*})$ is a KKT point of $\min_{\x}\max_{\blambda}g^k(\x)+\<\blambda,\U\x\>$ and $g^k(\x)= f(\x)+\frac{\tau}{2}\|\x-\y^k\|_F^2$. From the KKT condition, we have $\U\blambda^{k,*}+\nabla g^k(\x^{k,*})=0$. Thus, we have
\begin{eqnarray}
\begin{aligned}\label{cont004}
&\|\U\blambda^{k,*}-\U\blambda^{k-1,*}\|_F^2\\
&=\left\|\nabla f(\x^{k,*})+\tau(\x^{k,*}-\y^k)-\nabla f(\x^{k-1,*})-\tau(\x^{k-1,*}-\y^{k-1})\right\|_F^2\\
&\overset{c}\leq 2\left(L+\tau\right)^2\|\x^{k,*}-\x^{k-1,*}\|_F^2+2\tau^2\|\y^k-\y^{k-1}\|_F^2\\
&\overset{d}\leq 4L_g^2\|\y^k-\y^{k-1}\|_F^2,
\end{aligned}
\end{eqnarray}
where $\overset{c}\leq$ uses the $L$-smoothness of $f(\x)$ and $\overset{d}\leq$ uses (\ref{cont002}) and $L_g=L+\tau$. Combining (\ref{cont001}), (\ref{cont002}), (\ref{cont003}), and (\ref{cont004}) and using $\beta_g=L_g$, we have
\begin{eqnarray}
\begin{aligned}\label{cont3}
\rho_{k,0}\leq  2\rho_{k-1,T_{k-1}}+\left(4L_g+\frac{8L_g}{1-\sigma_2(\W)}\right)\|\y^k-\y^{k-1}\|_F^2.
\end{aligned}
\end{eqnarray}
From a similar induction to the proof of [Proposition 12]\cite{catalyst} and the relations in Algorithm \ref{acc-extra}, we have
\begin{eqnarray}
\begin{aligned}\label{cont6}
\|\y^k-\y^{k-1}\|_F^2\leq&2\|\y^k-\x^*\|_F^2+2\|\y^{k-1}-\x^*\|_F^2\\
\leq&4(1+\vartheta_k)^2\|\x^k-\x^*\|_F^2+4\vartheta_k^2\|\x^{k-1}-\x^*\|_F^2\\
&+4(1+\vartheta_{k-1})^2\|\x^{k-1}-\x^*\|_F^2+4\vartheta_{k-1}^2\|\x^{k-2}-\x^*\|_F^2\\
\leq&40\max\{\|\x^k-\x^*\|_F^2,\|\x^{k-1}-\x^*\|_F^2,\|\x^{k-2}-\x^*\|_F^2\},
\end{aligned}
\end{eqnarray}
where we denote $\vartheta_k=\frac{\theta_{k-1}(1-\theta_{k-1})}{\theta_{k-1}^2+\theta_k}$ and use $\vartheta_k\leq 1$ for all $k$. The latter can be obtained by $\vartheta_k=\frac{\sqrt{q}-q}{\sqrt{q}+q}\leq 1$ for $\mu>0$ and $\vartheta_k=\frac{\theta_{k-1}(1-\theta_{k-1})}{\theta_{k-1}^2/\theta_k}\leq\frac{\theta_k}{\theta_{k-1}}\leq 1$ for $\mu=0$.

Since $\rho_{r,T_r}\leq 2m\varepsilon_r$ for all $r\leq k-1$, i.e., $G^r(\alpha(\x^{r+1}))-G^r(x^{r,*})\leq\varepsilon_r$, from Theorem \ref{theorem_catalyst} we know the following conclusion holds for all $r\leq k-1$:
\begin{eqnarray}
\begin{aligned}\label{cont008}
F(\alpha(\x^{r+1}))-F(x^*)\leq \frac{36}{(\sqrt{q}-\rho)^2}\varepsilon_{r+1},
\end{aligned}
\end{eqnarray}
where we use the definition of $\varepsilon_r$ in Theorem \ref{theorem_catalyst}, Thus, we have
\begin{eqnarray}
\begin{aligned}\label{cont7}
\|\x^k-\x^*\|_F^2\overset{e}=&\|\1(\alpha(\x^k))^T+\Pi\x^k-\1(x^*)^T\|_F^2\\
\leq& 2m\|\alpha(\x^k)-x^*\|^2+2\|\Pi\x^k\|_F^2\\
\overset{f}\leq& \frac{4m}{\mu}(F(\alpha(\x^k))-F(x^*))+\frac{ 32\rho_{k-1,T_{k-1}}}{\beta_g(1-\sigma_2(\W))}\\
\overset{g}\leq& \frac{144m\varepsilon_{k}}{\mu(\sqrt{q}-\rho)^2}+\frac{64m\varepsilon_{k-1}}{\beta_g}
\end{aligned}
\end{eqnarray}
where we use the definitions of $\Pi\x$ and $\alpha(\x)$ in $\overset{e}=$, the $\mu$-strong convexity of $F(x)$ and (\ref{cont2}) in $\overset{f}\leq$, (\ref{cont008}) and the induction condition of $\rho_{k-1,T_{k-1}}\leq 2m(1-\sigma_2(\W))\varepsilon_{k-1}$ in $\overset{g}\leq$.

Combining (\ref{cont1}), (\ref{cont3}), (\ref{cont6}), (\ref{cont7}), and using $\rho_{k-1,T_{k-1}}\leq 2m(1-\sigma_2(\W))\varepsilon_{k-1}$, we have
\begin{eqnarray}
\begin{aligned}\notag
\rho_{k,T_k}\leq&
(1\hspace*{-0.03cm}-\hspace*{-0.03cm}\delta_g)^{T_k}\varepsilon_k\hspace*{-0.03cm}\left(\hspace*{-0.03cm}\frac{4m}{1\hspace*{-0.03cm}-\hspace*{-0.03cm}\rho}\hspace*{-0.03cm}+\hspace*{-0.03cm}\left(\hspace*{-0.03cm}4L_g\hspace*{-0.03cm}+\hspace*{-0.03cm}\frac{8L_g}{1\hspace*{-0.03cm}-\hspace*{-0.03cm}\sigma_2(\W)}\hspace*{-0.03cm}\right)\hspace*{-0.03cm}\frac{40}{(1\hspace*{-0.03cm}-\hspace*{-0.03cm}\rho)^3}\hspace*{-0.03cm}\left(\hspace*{-0.03cm}\frac{144m\varepsilon_{k}}{\mu(\sqrt{q}\hspace*{-0.03cm}-\hspace*{-0.03cm}\rho)^2}\hspace*{-0.03cm}+\hspace*{-0.03cm}\frac{64m\varepsilon_{k-1}}{\beta_g}\hspace*{-0.03cm}\right)\hspace*{-0.03cm}\right)\\
\overset{h}\leq&(1-\delta_g)^{T_k}\frac{99844mL_g}{\mu (1-\sigma_2(\W))(1-\rho)^3(\sqrt{q}-\rho)^2}\varepsilon_k\equiv(1-\delta_g)^{T_k}C_1\varepsilon_k,
\end{aligned}
\end{eqnarray}
where we use $\sqrt{q}-\rho<1$, $\varepsilon_k\leq\varepsilon_{k-1}$ and $\beta_g\geq\mu$ in $\overset{h}\leq$.

Thus, to attain $\rho_{k,T_k}\leq 2m(1-\sigma_2(\W))\varepsilon_k$, we only need $(1-\delta_g)^{T_k}C_1\leq 2m(1-\sigma_2(\W))$, i.e., $T_k=O(\frac{1}{\delta_g}\log \frac{C_1}{2m(1-\sigma_2(\W))})=O((\frac{L_g}{\mu_g}+\frac{1}{1-\sigma_2(\W)})\log\frac{L_g}{\mu(1-\sigma_2(\W))})$.\hfill$\Box$

Based on the above lemma and Theorem \ref{acc_extra_the1}, we can prove Corollary \ref{acc_extra_cor1}.

\noindent\emph{Proof of Corollary \ref{acc_extra_cor1}:} From (\ref{cont2}) and Lemma \ref{lemma_cont2}, we have
\begin{eqnarray}
\begin{aligned}\notag
&\|\Pi\x^{K+1}\|_F^2\leq \frac{32m\varepsilon_K}{\beta_g}\overset{b}\leq \frac{32m}{\beta_g}\frac{2(F(x^0)-F(x^*))}{9}(1-\rho)^{K+1},
\end{aligned}
\end{eqnarray}
where $\overset{b}\leq$ uses the definition of $\varepsilon_k$ in Theorem \ref{theorem_catalyst}. On the other hand, from Theorem \ref{theorem_catalyst}, we have
\begin{eqnarray}
\begin{aligned}\notag
\qquad F(\alpha(\x^{K+1}))-F(x^*)\leq\frac{8}{(\sqrt{q}-\rho)^2}(1-\rho)^{K+2}(F(x^0)-F(x^*)).
\end{aligned}
\end{eqnarray}
To make $\|\Pi\x^{K+1}\|_F^2\leq O(m\epsilon^2)$ and $F(\alpha(\x^{K+1}))-F(x^*)\leq O(\epsilon)$, we only need to run Algorithm \ref{acc-extra} for $K=O(\sqrt{1+\frac{\tau}{\mu}}\log\frac{1}{\epsilon})$ outer iterations such that
\begin{eqnarray}
\begin{aligned}\notag
(F(x^0)-F(x^*))(1-\rho)^{K+1}\leq \epsilon^2.
\end{aligned}
\end{eqnarray}
Thus, the total number of inner iterations is
\begin{eqnarray}
\begin{aligned}\notag
\sum_{k=0}^{K}T_k=&\sqrt{1+\frac{\tau}{\mu}}\left(\log\frac{1}{\epsilon}\right)\left(\frac{L+\tau}{\mu+\tau}+\frac{1}{1-\sigma_2(\W)}\right)\log\frac{L+\tau}{\mu(1-\sigma_2(\W))}\\
\leq& 3\sqrt{\frac{L}{\mu(1-\sigma_2(\W))}}\log\frac{2L}{\mu(1-\sigma_2(\W))}\log\frac{1}{\epsilon}
\end{aligned}
\end{eqnarray}
by letting $\tau=L(1-\sigma_2(\W))-\mu$.

\subsubsection{Non-strongly Convex Case}\label{sec:ns}
When the strong convexity is absent, we can have the following lemma, which further leads to Theorem \ref{acc_extra_the2}. Similar to Lemma \ref{lemma_cont2}, we prove a sharper result of $\rho_{k,T_k}\leq 2m\varepsilon_k(1-\sigma_2(\W))^{3+\xi}$.
\begin{lemma}\label{lemma_cont3}
Suppose that $F(x)$ is convex. If $\rho_{r,T_r}\leq 2m\varepsilon_r(1-\sigma_2(\W))^{3+\xi}$ holds for all $r\leq k-1$ and we initialize $\x^{k,0}=\x^{k-1,T_{k-1}+1}$ and $\blambda^{k,0}=\blambda^{k-1,T_{k-1}+1}$, then we only need $T_k=O((\frac{L_g}{\mu_g}+\frac{1}{1-\sigma_2(\W)})\log\frac{k}{1-\sigma_2(\W)})$ such that $\rho_{k,T_k}\leq 2m\varepsilon_k(1-\sigma_2(\W))^{3+\xi}$.
\end{lemma}

The proof is similar to that of [Proposition 12]\cite{catalyst}, and we omit the details. Simply, when the strong convexity is absent, (\ref{cont3}) and (\ref{cont6}) also hold. But we need to bound $\|\x^k-\x^*\|_F^2$ in a different way. From Theorem \ref{theorem_catalyst}, the sequence $F(\alpha(\x^k))$ is bounded by a constant. By the bounded level set assumption, there exists $C>0$ such that $\|\alpha(\x^k)-x^*\|\leq C$. From (\ref{cont2}), we have $\|\Pi\x^k\|_F^2\leq \frac{16}{\beta_g(1-\sigma_2(\W))}\rho_{k-1,T_{k-1}}$. Thus, we have
\begin{eqnarray}
\begin{aligned}\notag
\|\x^k-\x^*\|_F^2=&\|\1(\alpha(\x^k))^T+\Pi\x^k-\1(x^*)^T\|_F^2\\
\leq& 2m\|\alpha(\x^k)-x^*\|^2+2\|\Pi\x^k\|_F^2\\
\leq& 2mC^2+\frac{32m\varepsilon_{k-1}(1-\sigma_2(\W))^{2+\xi}}{\beta_g}.
\end{aligned}
\end{eqnarray}
Thus, $\rho_{k,0}$ is bounded by constant $C_2=4m+40(4L_g+\frac{8L_g}{1-\sigma_2(\W)})(2mC^2+32m/\beta_g)$ and we only need $(1-\delta_g)^{T_k}C_2\leq 2m\varepsilon_k(1-\sigma_2(\W))^{3+\xi}$, i.e., $T_k=O(\frac{1}{\delta_g}\log\frac{C_2}{2m\varepsilon_k(1-\sigma_2(\W))^{3+\xi}})$ $=O((\frac{L_g}{\mu_g}+\frac{1}{1-\sigma_2(\W)})\log\frac{k}{1-\sigma_2(\W)})$.

Now, we come to Corollary \ref{acc_extra_cor2}. From Theorem \ref{theorem_catalyst}, to find an $\epsilon$-optimal solution such that $F(\alpha(\x^{K+1}))-F(x^*)\leq\epsilon$, we need $K=O(\sqrt{\frac{\tau R_1}{\epsilon}})$ outer iterations. On the other hand, from (\ref{cont2}) and Lemma \ref{lemma_cont3}, we have
\begin{eqnarray}
\begin{aligned}\notag
&\|\Pi\x^{K+1}\|_F^2\leq \frac{32m(1-\sigma_2(\W))^{2+\xi}\varepsilon_K}{\beta_g}
\overset{a}\leq \frac{32m(1-\sigma_2(\W))^{2+\xi}}{\beta_g K^{4+2\xi}}\overset{b}\leq\frac{32m\epsilon^2}{\beta_gL^{2+\xi}R_1^{2+\xi}},
\end{aligned}
\end{eqnarray}
where $\overset{a}\leq$ uses the definition $\varepsilon_k$ in Theorem \ref{theorem_catalyst}, $\overset{b}\leq$ uses $K=O(\sqrt{\frac{\tau R_1}{\epsilon}})$ and $\tau=L(1-\sigma_2(\W))$. Thus, the settings of $T_k$ and $K$ lead to $\|\Pi\x^{K+1}\|_F^2\leq O(m\epsilon^2)$. The total number of inner iterations is
\begin{eqnarray}
\begin{aligned}\notag
\sum_{k=0}^{\sqrt{\frac{\tau}{\epsilon}}}T_k= \sqrt{\frac{\tau}{\epsilon}}\left(\frac{L+\tau}{\tau}+\frac{1}{1-\sigma_2(\W)}\right)\log\frac{1}{\epsilon(1-\sigma_2(\W))}.
\end{aligned}
\end{eqnarray}
The setting of $\tau=L(1-\sigma_2(\W))$ leads to the minimal value of $\sqrt{\frac{L}{\epsilon(1-\sigma_2(\W))}}\log\frac{1}{\epsilon(1-\sigma_2(\W))}$.

\section{Numerical Experiments}\label{sec:exp}
Consider the decentralized least squares problem:
\begin{eqnarray}
\min_{x\in\R^{n}} \sum_{i=1}^m f_i(x)\quad \mbox{with}\quad f_i(x)\equiv\frac{1}{2}\|\A_i^T x-\b_i\|^2+\frac{\mu}{2}\|x\|^2,\label{problem_exp}
\end{eqnarray}
where each agent $\{1,...,m\}$ holds its own local function $f_i(x)$. $\A_i\in\R^{n\times s}$ is generated from the uniform distribution with each entry in $[0,1]$, and each column of $\A_i$ is normalized to be 1. We set $s=10$, $n=500$, $m=100$, and $\b_i=\A_i^Tx$ with some unknown $x$. We test the performance of the proposed algorithms on both the strongly convex problem and nonstrongly convex one. For the strongly convex case, we test on $\mu=10^{-6}$ and $\mu=10^{-8}$, respectively. In general, the accelerated algorithms apply to ill-conditioned problems with large condition numbers. For the nonstrongly convex one, we let $\mu=0$.

We test the performance on two kinds of networks: (1) The first is Erd\H{o}s$-$R\'{e}nyi random graph, where each pair of nodes has a connection with the ratio of $p$. We test two different settings of $p$: $p=0.5$ and $p=0.1$, which results in $\frac{1}{1-\sigma_2(\W)}=2.87$ and $\frac{1}{1-\sigma_2(\W)}=7.74$, respectively. (2) The second is the geometric graph, where $m$ nodes are placed uniformly and independently in the unit square $[0,1]$ and two nodes are connected if their distance is at most $d$. We test on $d=0.5$ and $d=0.3$, which leads to $\frac{1}{1-\sigma_2(\W)}=8.13$ and $\frac{1}{1-\sigma_2(\W)}=30.02$, respectively. We set the weight matrix as $\W=\frac{\I+\M}{2}$ for both graphes, where $\M$ is the Metropolis weight matrix \cite{boyd2004}: $\M_{i,j}=\left\{
  \begin{array}{ll}
    1/(1+\max\{d_i,d_j\}),& \mbox{if }(i,j)\in\mathcal{E},\\
    0,& \mbox{if }(i,j)\notin\mathcal{E}\mbox{ and }i\neq j,\\
    1-\sum_{l\in\N_i}\W_{i,l}, &\mbox{if }i=j,
  \end{array}
\hspace*{-0.2cm}\right.$ and $d_i$ is the number of the $i$-th agent's neighbors.

\begin{figure}
\centering
\begin{tabular}{@{\extracolsep{0.001em}}c@{\extracolsep{0.001em}}c@{\extracolsep{0.001em}}c@{\extracolsep{0.001em}}c@{\extracolsep{0.001em}}c}
\hspace*{-0.8cm}\includegraphics[width=0.32\textwidth,keepaspectratio]{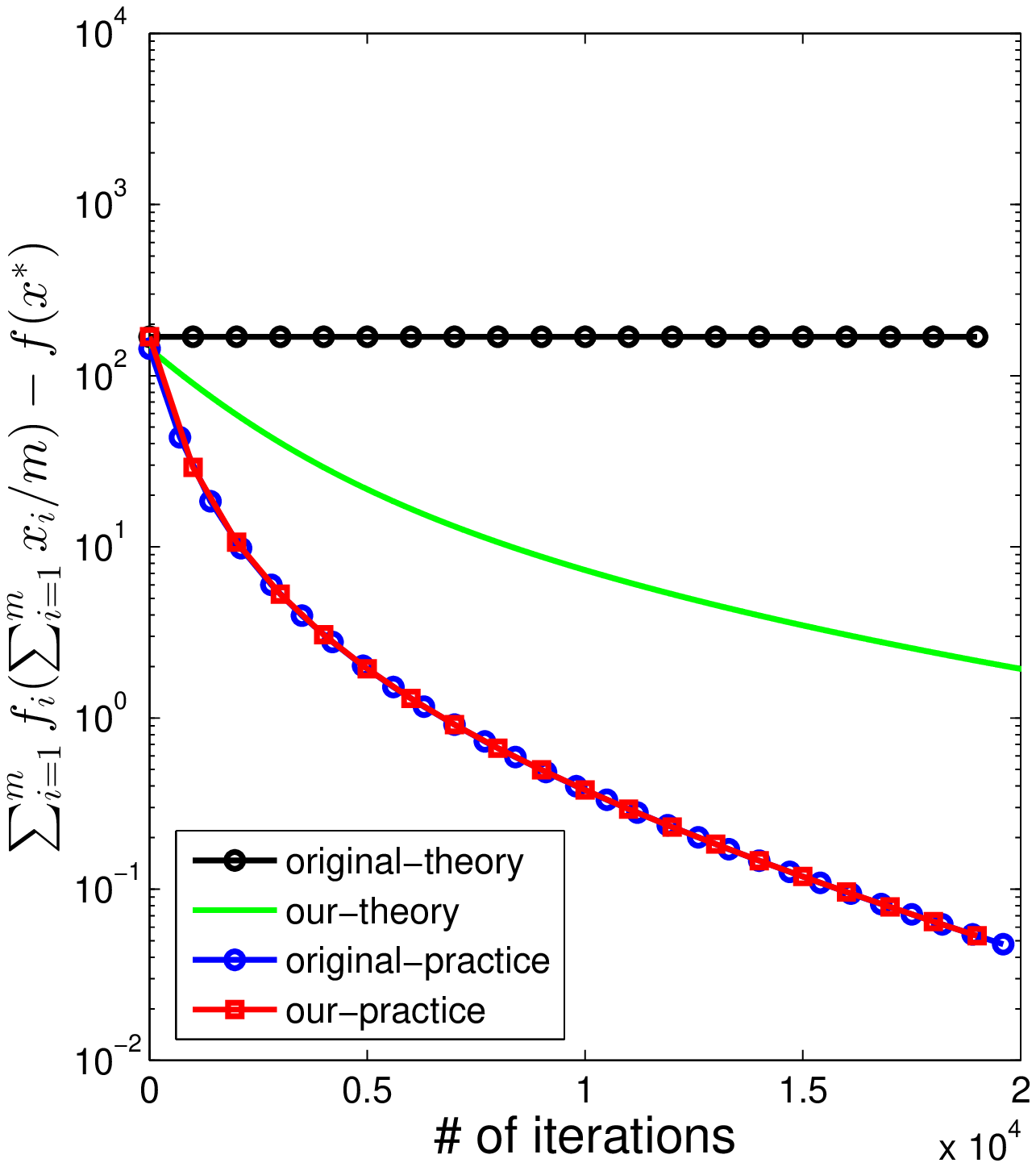}
&\hspace*{-0.28cm}\includegraphics[width=0.26\textwidth,keepaspectratio]{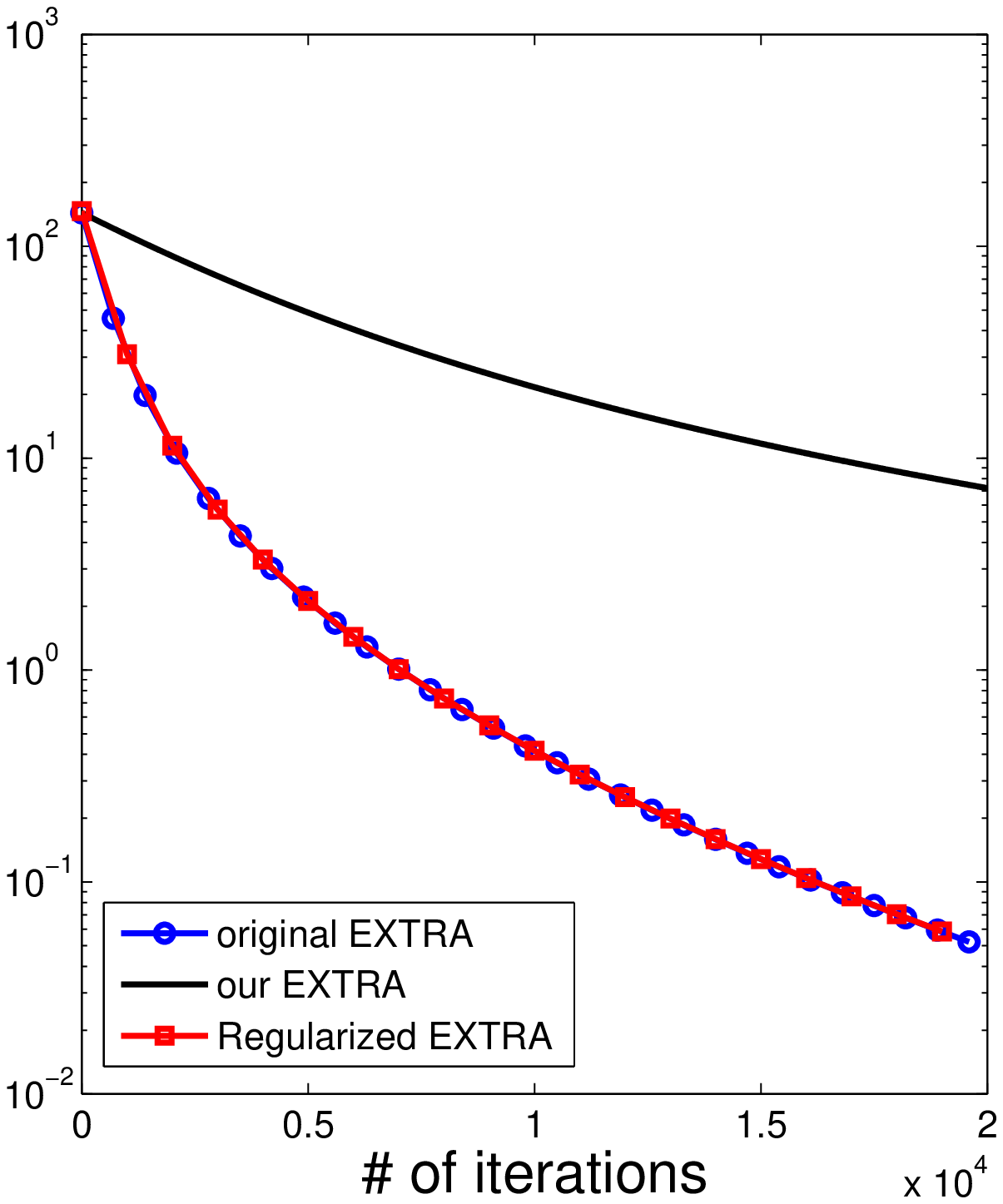}
&\hspace*{-0.28cm}\includegraphics[width=0.26\textwidth,keepaspectratio]{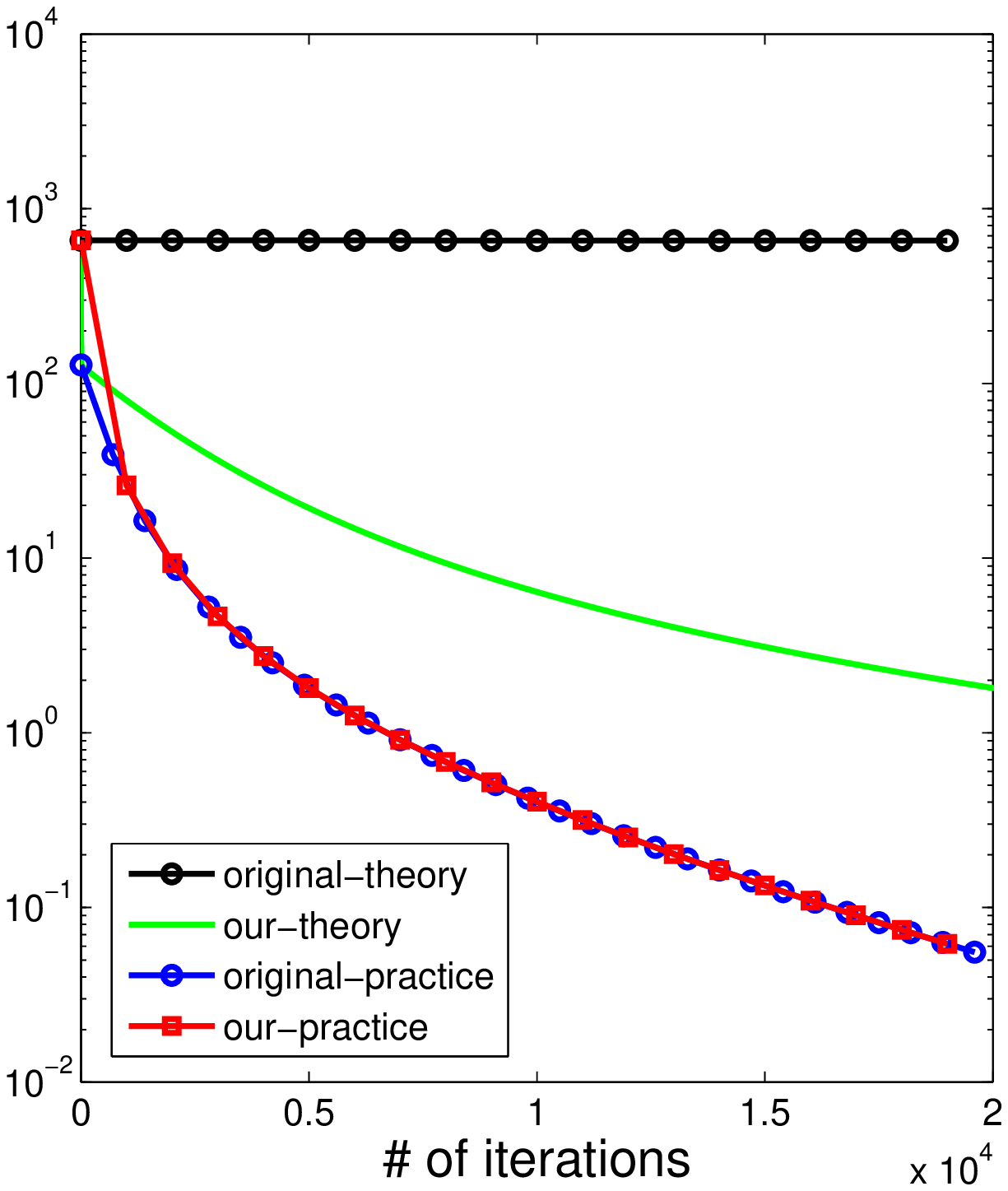}
&\hspace*{-0.28cm}\includegraphics[width=0.26\textwidth,keepaspectratio]{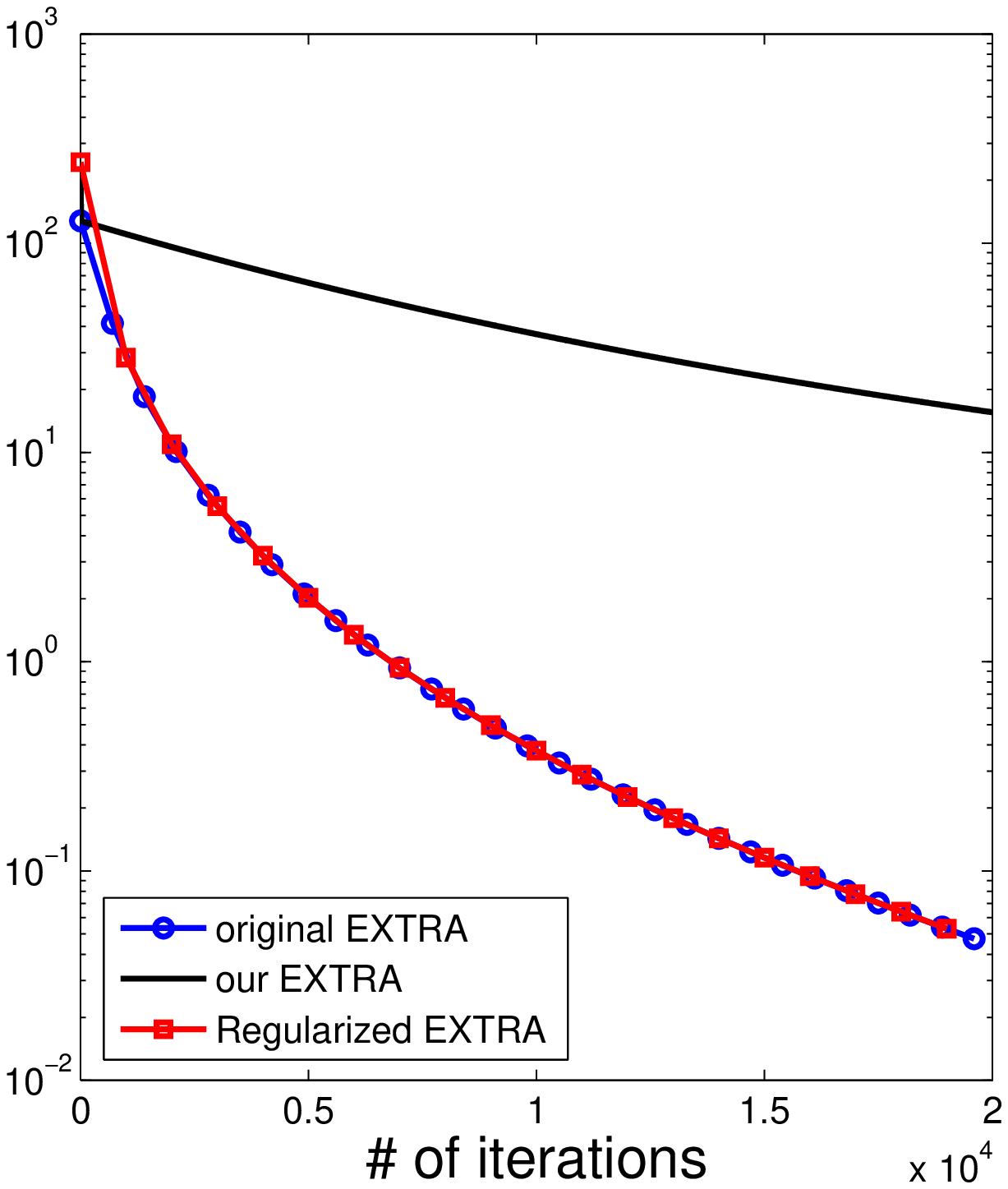}\\
\scriptsize (a). SC, Erd\H{o}s$-$R\'{e}nyi graph&\scriptsize (b). NS, Erd\H{o}s$-$R\'{e}nyi graph&\scriptsize (c). SC, Geometric graph&\scriptsize (d). NS, Geometric graph\normalsize
\end{tabular}
\caption{Comparisons between different EXTRA on the Erd\H{o}s$-$R\'{e}nyi random graph (p=0.1) and the geometric graph (d=0.3). SC means the strongly convex problem ($\mu=10^{-6}$), and NS means the nonstrongly convex one.}\label{fig0}
\end{figure}

We first compare EXTRA analyzed in this paper with the original EXTRA \cite{Shi-2015}. For the strongly convex problem, the authors of [Remark 4]\cite{Shi-2015} analyzed the algorithm with $\alpha=\frac{1}{\beta}=\frac{\mu^2}{L}$ and $\alpha=\frac{1}{\beta}=\frac{1}{L}$ being suggested in practice. In our theory, we use $\beta=L$ and $\alpha=\frac{1}{4L}$. In practice, we observe that $\beta=L$ and $\alpha=\frac{1}{L}$ performs the best. Figures \ref{fig0}(a) and \ref{fig0}(c) plot the results. We can see that the theoretical setting in the original EXTRA makes almost no decreasing in the objective function values due to small step-size and that our theoretical setting works much better. On the other hand, both the original EXTRA and our analyzed one work best for $\beta=L$ and $\alpha=\frac{1}{L}$. For the nonstrongly convex problems, \cite{Shi-2015} suggests $\alpha=\frac{1}{\beta}=\frac{1}{L}$ in both theory and practice. In our theory, Lemma \ref{nsc_lemma} suggests $\beta=\frac{L}{\sqrt{1-\sigma_2(\W)}}$ and $\alpha=\frac{1}{2(L+\beta)}$. From Figure \ref{fig0}.b and Figure \ref{fig0}.d, we observe that a larger $\beta$ (i.e., a smaller step-size) makes the algorithm slow. On the other hand, our regularized EXTRA performs as well as the original EXTRA.

Then, we compare the proposed accelerated EXTRA (Acc-EXTRA) with the original EXTRA \cite{Shi-2015}, accelerated distributed Nesterov gradient descent (Acc-DNGD) \cite{qu2017-2}, accelerated dual ascent (ADA) \cite{Uribe-2017} and the accelerated penalty method with consensus (APM-C) \cite{li-2018-pm}. For the strongly convex problem, we set $\tau=L(1-\sigma_2(\W))-\mu$ and $T_k=\lceil\frac{1}{5(1-\sigma_2(\W))}\log\frac{L}{\mu(1-\sigma_2(\W))}\rceil$ for Acc-EXTRA, $T_k=\lceil\frac{k\sqrt{\mu/L}}{4\sqrt{1-\sigma_2(\W)}}\rceil$ and the step-size as $\frac{1}{L}$ for APM-C, $T_k=\lceil\sqrt{\frac{L}{\mu}}\log\frac{L}{\mu}\rceil$ and the step-size as $\mu$ for ADA, where $T_k$ means the number of inner iterations at the $k$th outer iteration and $\lceil\cdot\rceil$ is the top integral function. We set the step-size as $\frac{1}{L}$ for EXTRA and tune the best step-size for Acc-DNGD with different graphs and different $\mu$. All the compared algorithms start from $x_{(i)}=\0$ for all $i$.

\begin{figure}
\centering
\begin{tabular}{@{\extracolsep{0.001em}}c@{\extracolsep{0.001em}}c@{\extracolsep{0.001em}}c@{\extracolsep{0.001em}}c@{\extracolsep{0.001em}}c}
\hspace*{-0.8cm}\includegraphics[width=0.325\textwidth,keepaspectratio]{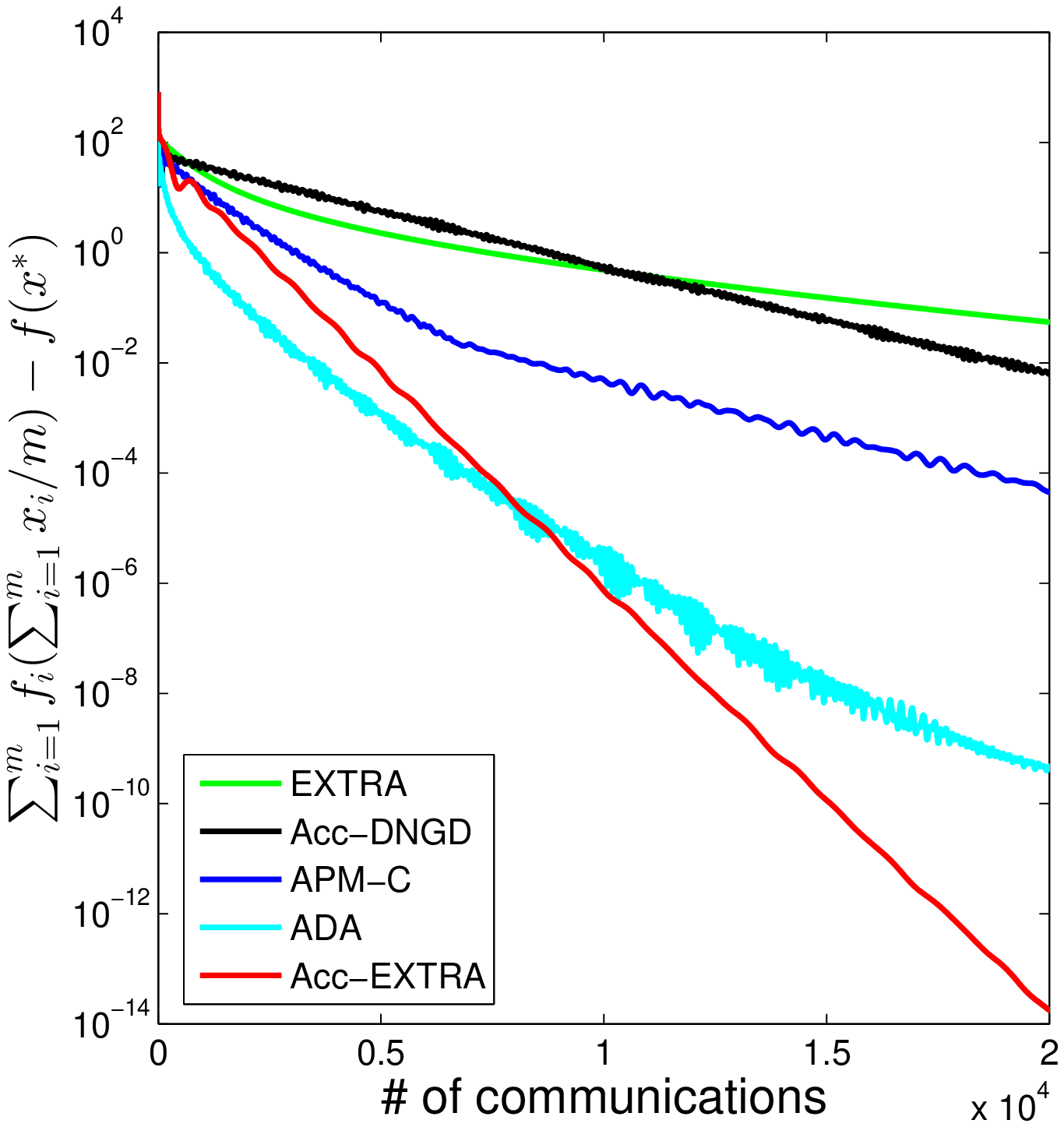}
&\hspace*{-0.28cm}\includegraphics[width=0.26\textwidth,keepaspectratio]{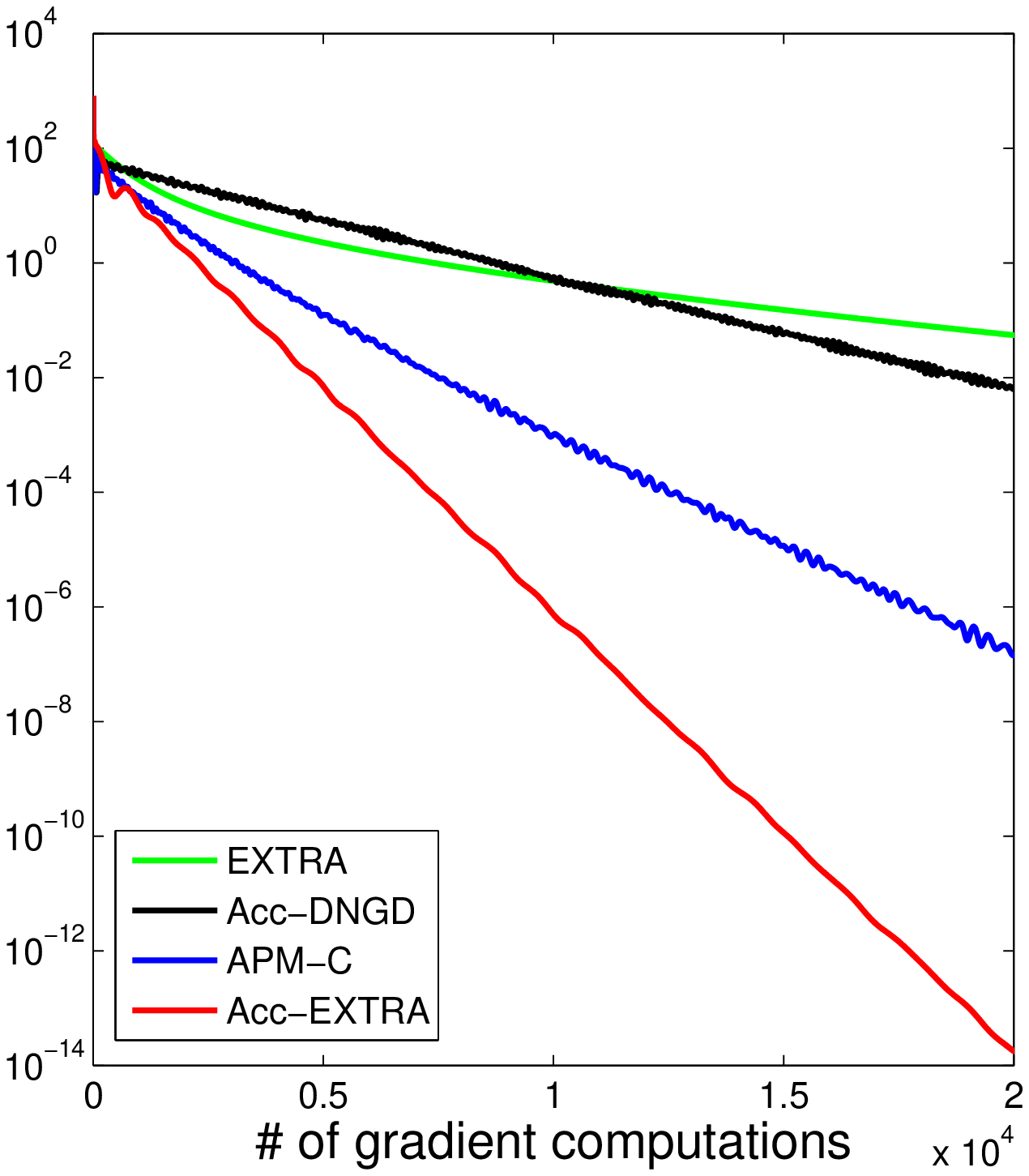}
&\hspace*{-0.28cm}\includegraphics[width=0.26\textwidth,keepaspectratio]{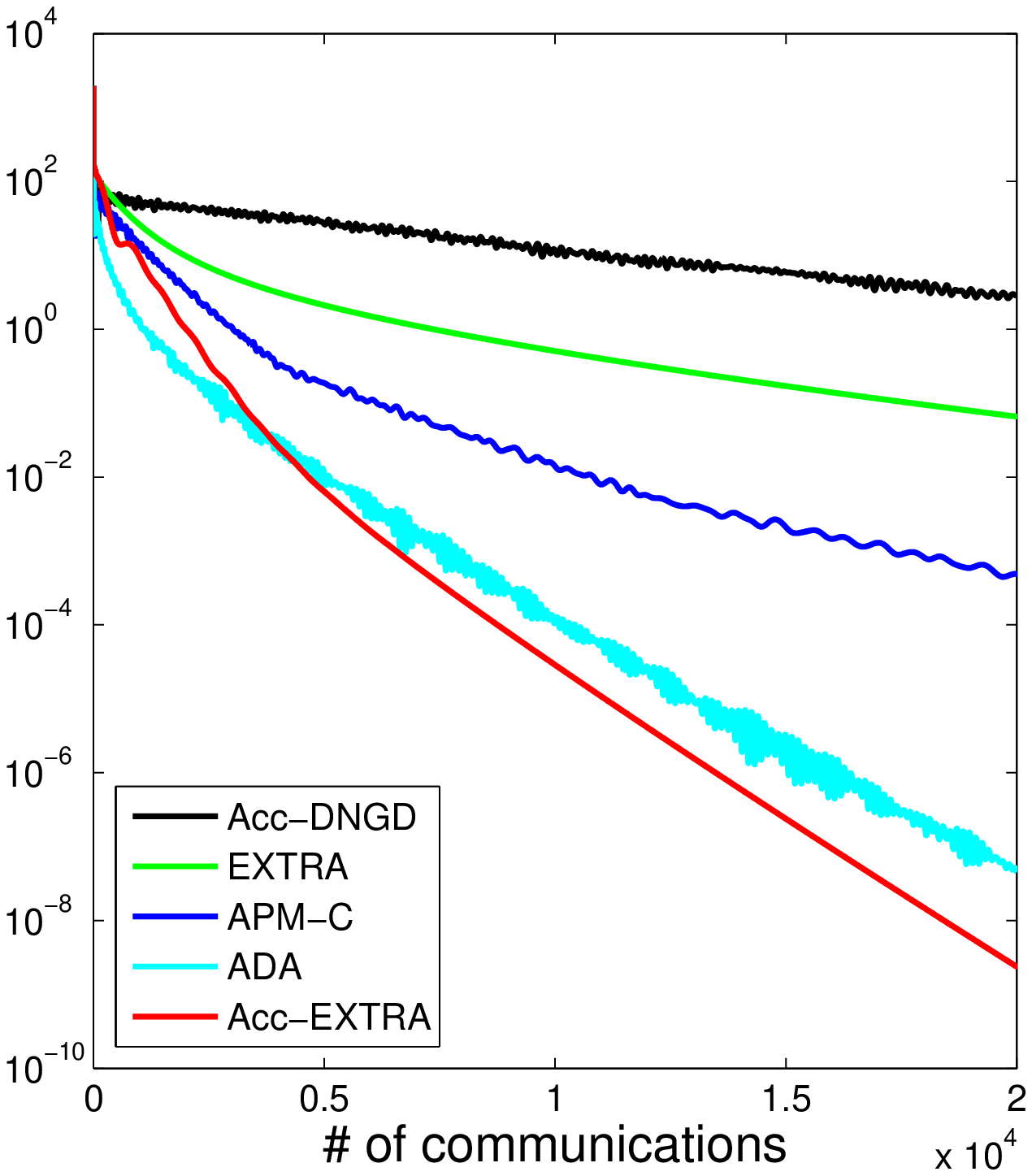}
&\hspace*{-0.28cm}\includegraphics[width=0.26\textwidth,keepaspectratio]{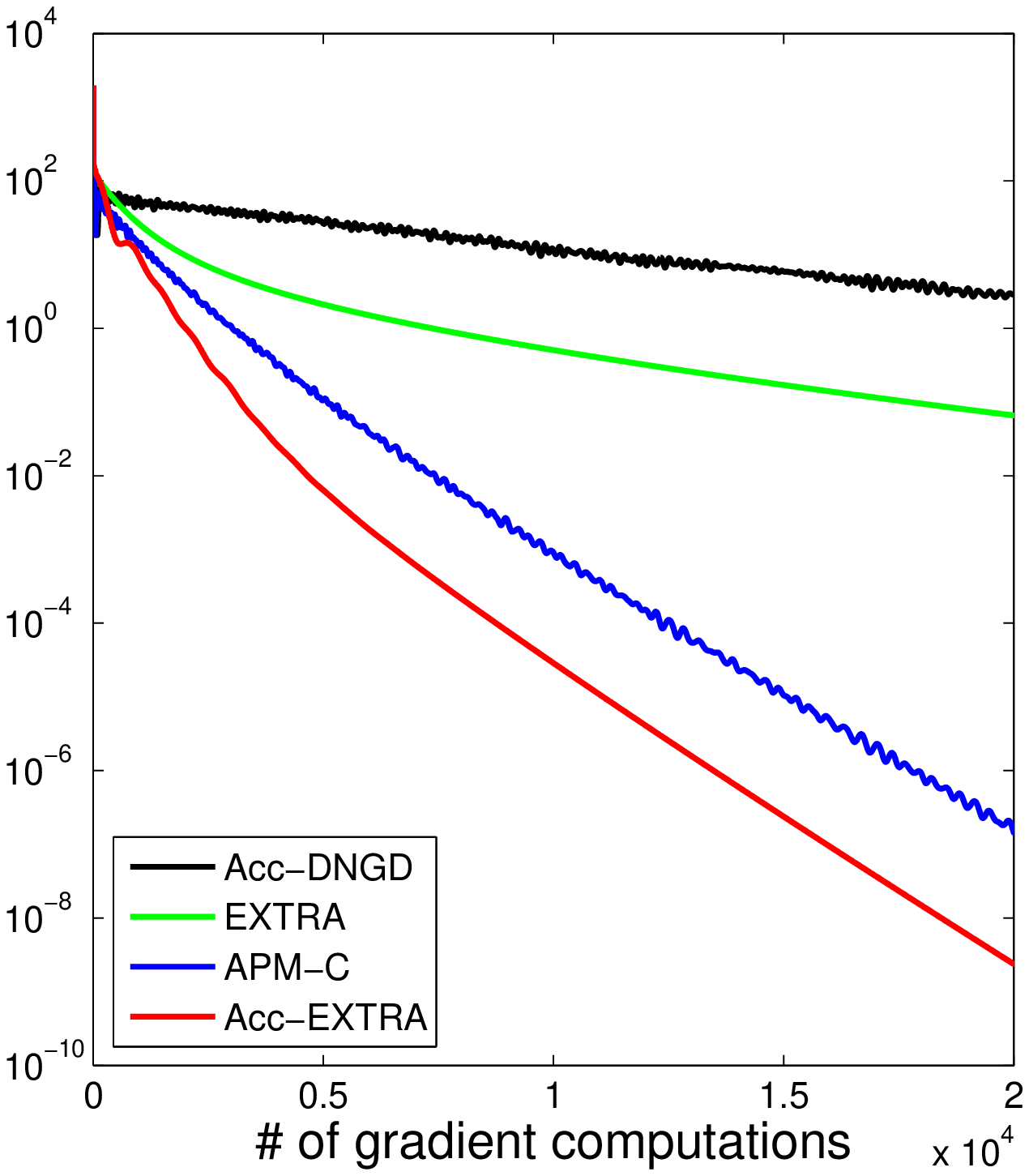}\\
\hspace*{-0.8cm}\includegraphics[width=0.325\textwidth,keepaspectratio]{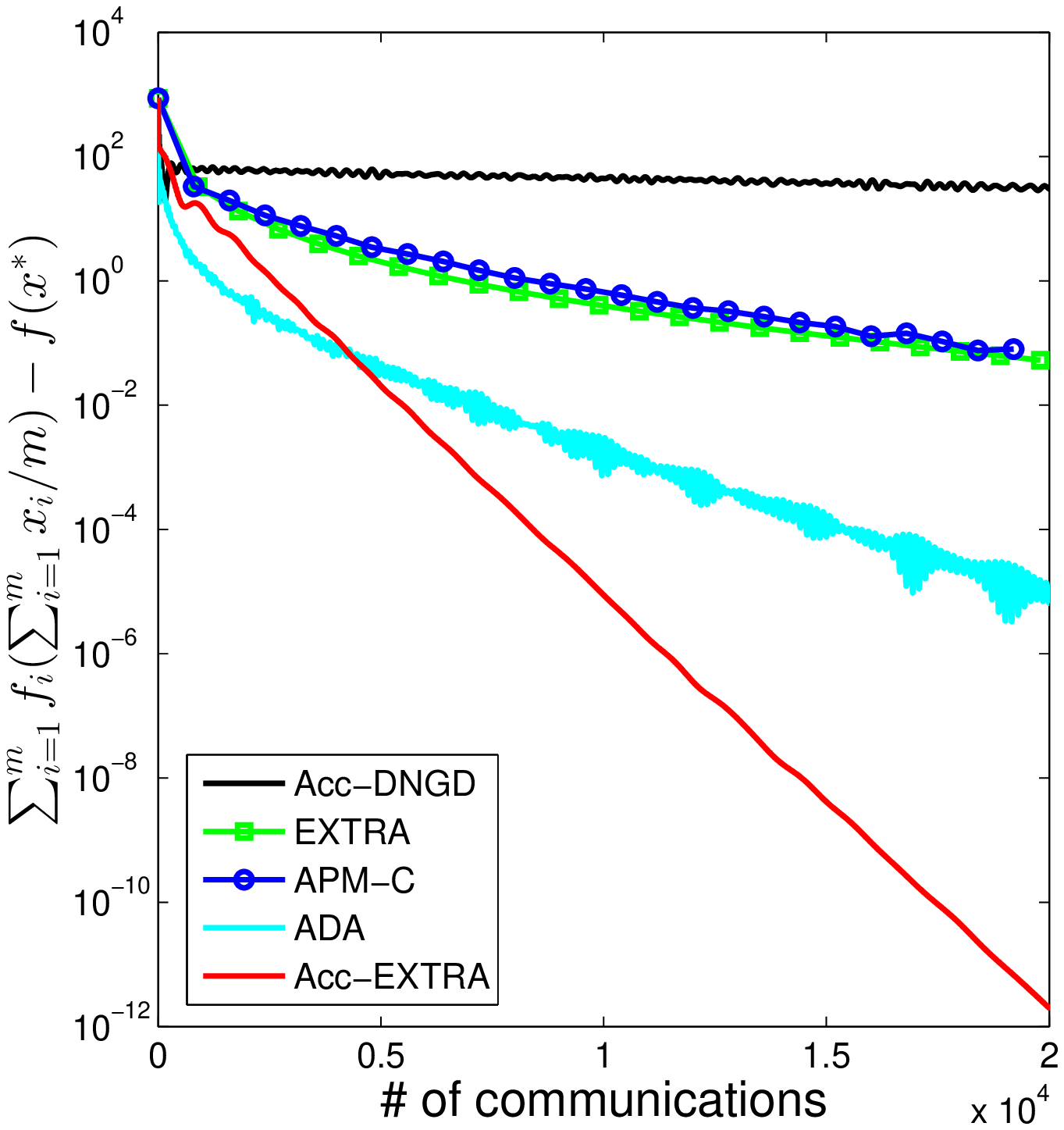}
&\hspace*{-0.28cm}\includegraphics[width=0.26\textwidth,keepaspectratio]{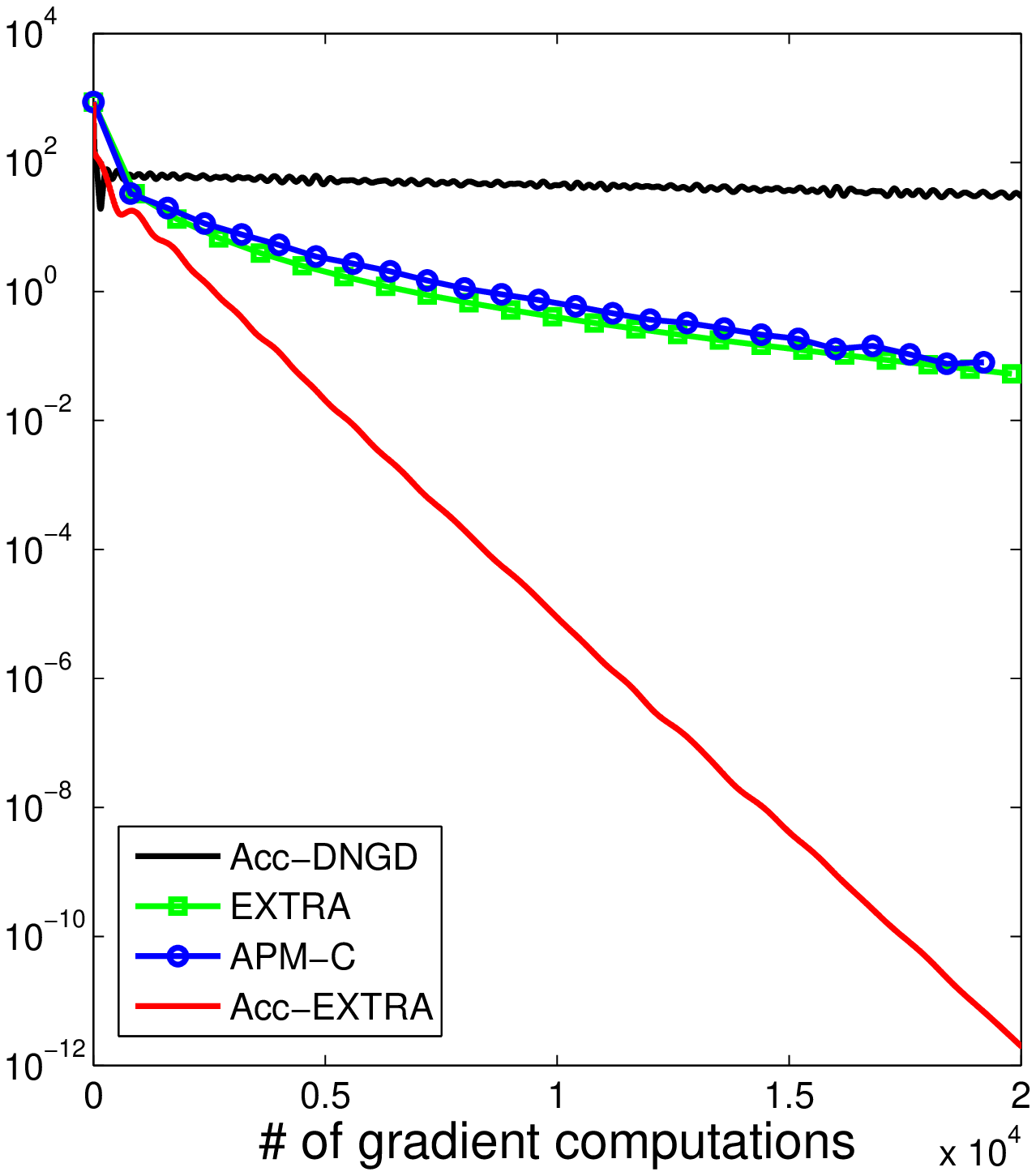}
&\hspace*{-0.28cm}\includegraphics[width=0.26\textwidth,keepaspectratio]{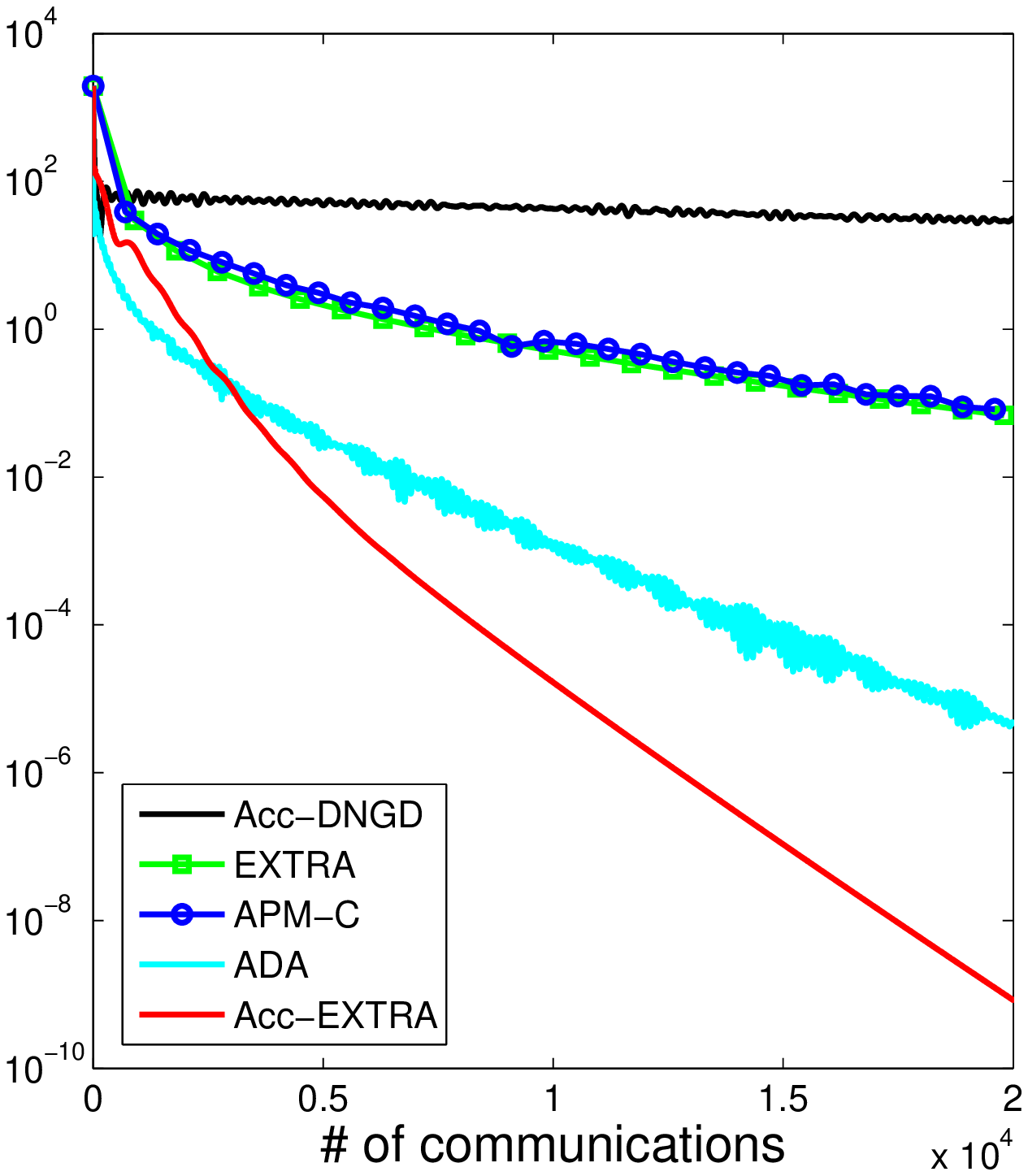}
&\hspace*{-0.28cm}\includegraphics[width=0.26\textwidth,keepaspectratio]{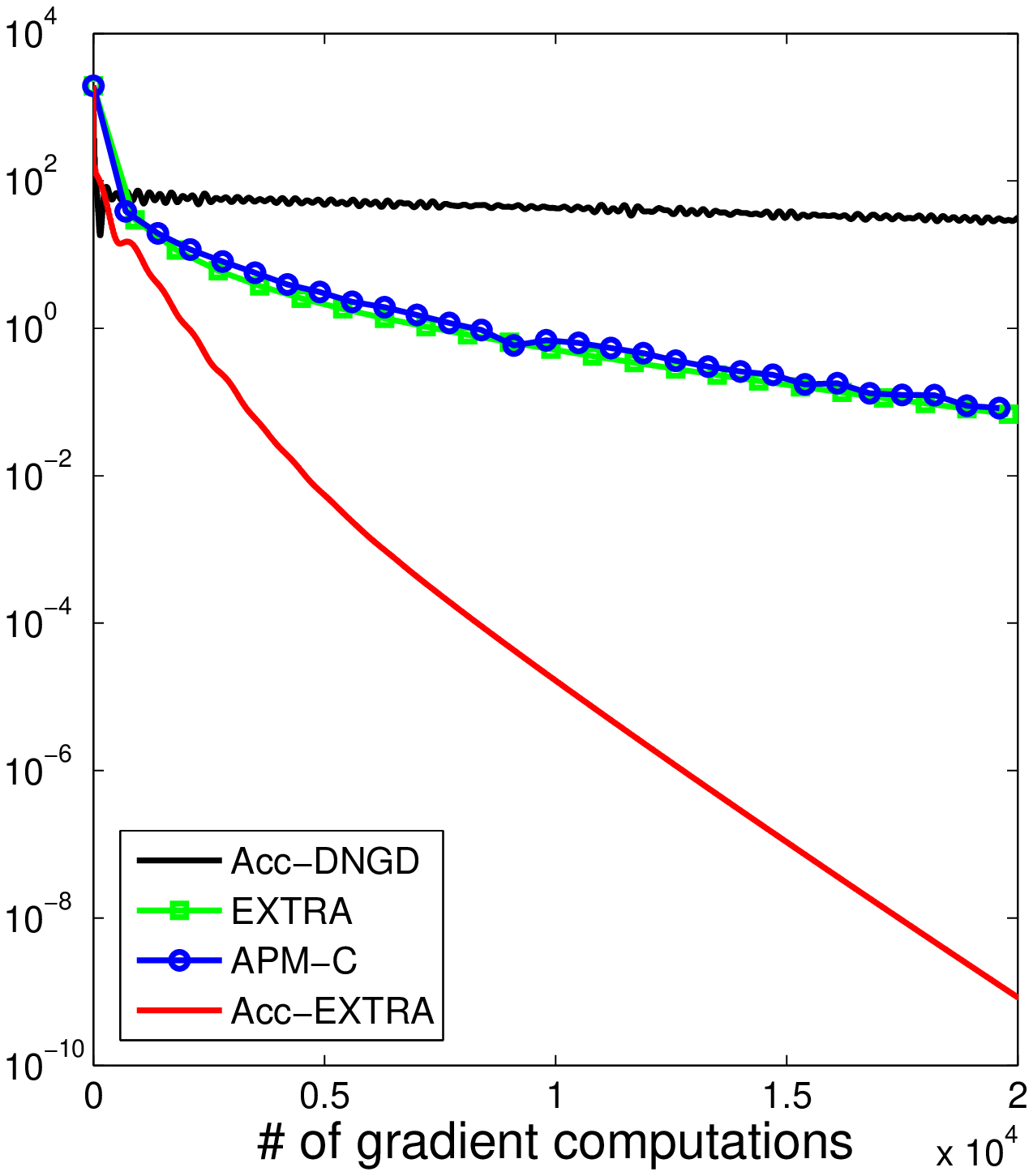}\\
\scriptsize$\frac{1}{1-\sigma_2(\W)}=2.87$&\scriptsize$\frac{1}{1-\sigma_2(\W)}=2.87$&\scriptsize$\frac{1}{1-\sigma_2(\W)}=7.74$&\scriptsize$\frac{1}{1-\sigma_2(\W)}=7.74$\normalsize\\
\end{tabular}
\caption{Comparisons on the strongly convex problem with the Erd\H{o}s$-$R\'{e}nyi random graph. $p=0.5$ for the two left plots, and $p=0.1$ for the two right. $\mu=10^{-6}$ for the top four plots, and $\mu=10^{-8}$ for the bottom four.}\label{fig1}
\end{figure}

\begin{figure}
\centering
\begin{tabular}{@{\extracolsep{0.001em}}c@{\extracolsep{0.001em}}c@{\extracolsep{0.001em}}c@{\extracolsep{0.001em}}c@{\extracolsep{0.001em}}c}
\hspace*{-0.8cm}\includegraphics[width=0.325\textwidth,keepaspectratio]{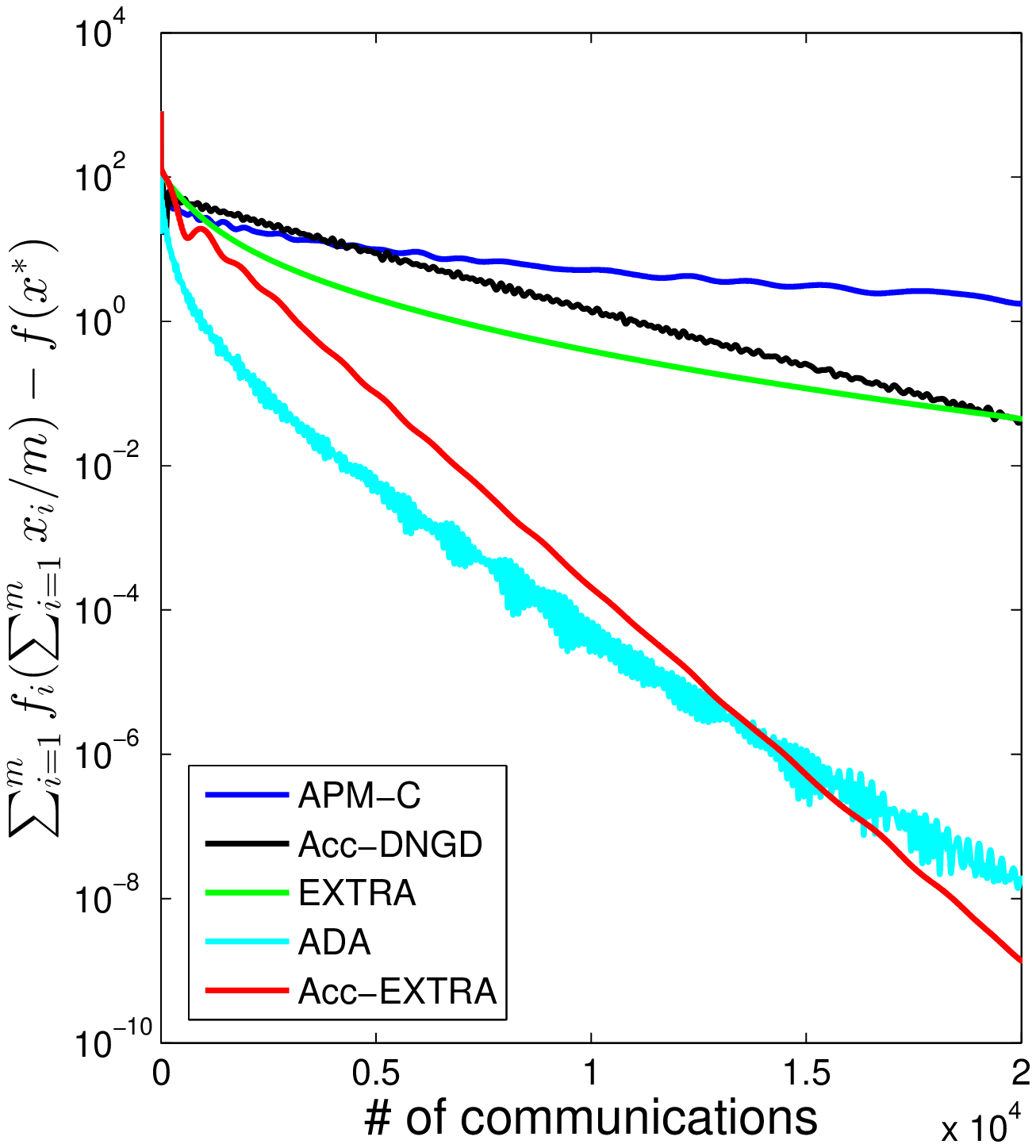}
&\hspace*{-0.28cm}\includegraphics[width=0.26\textwidth,keepaspectratio]{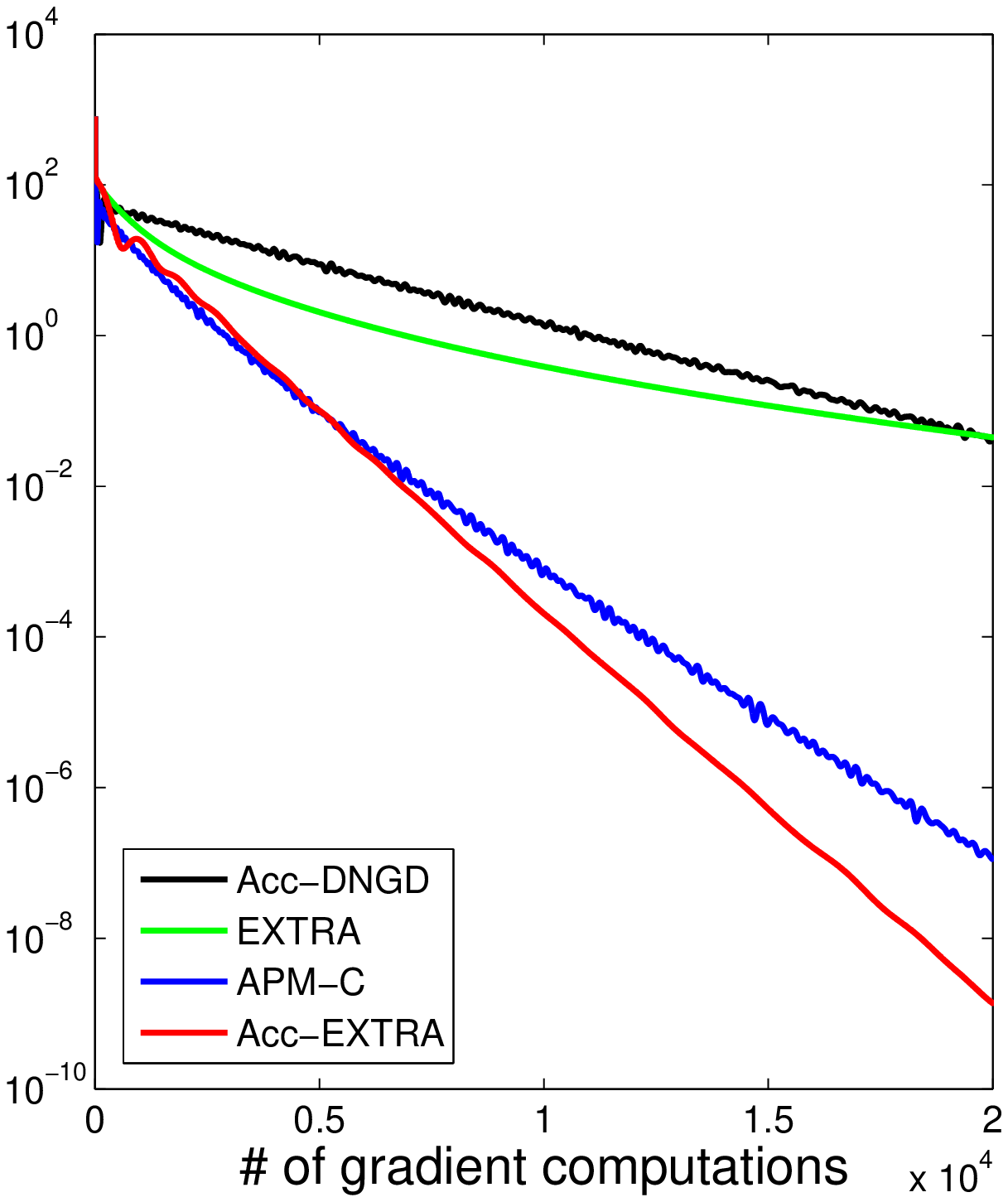}
&\hspace*{-0.28cm}\includegraphics[width=0.26\textwidth,keepaspectratio]{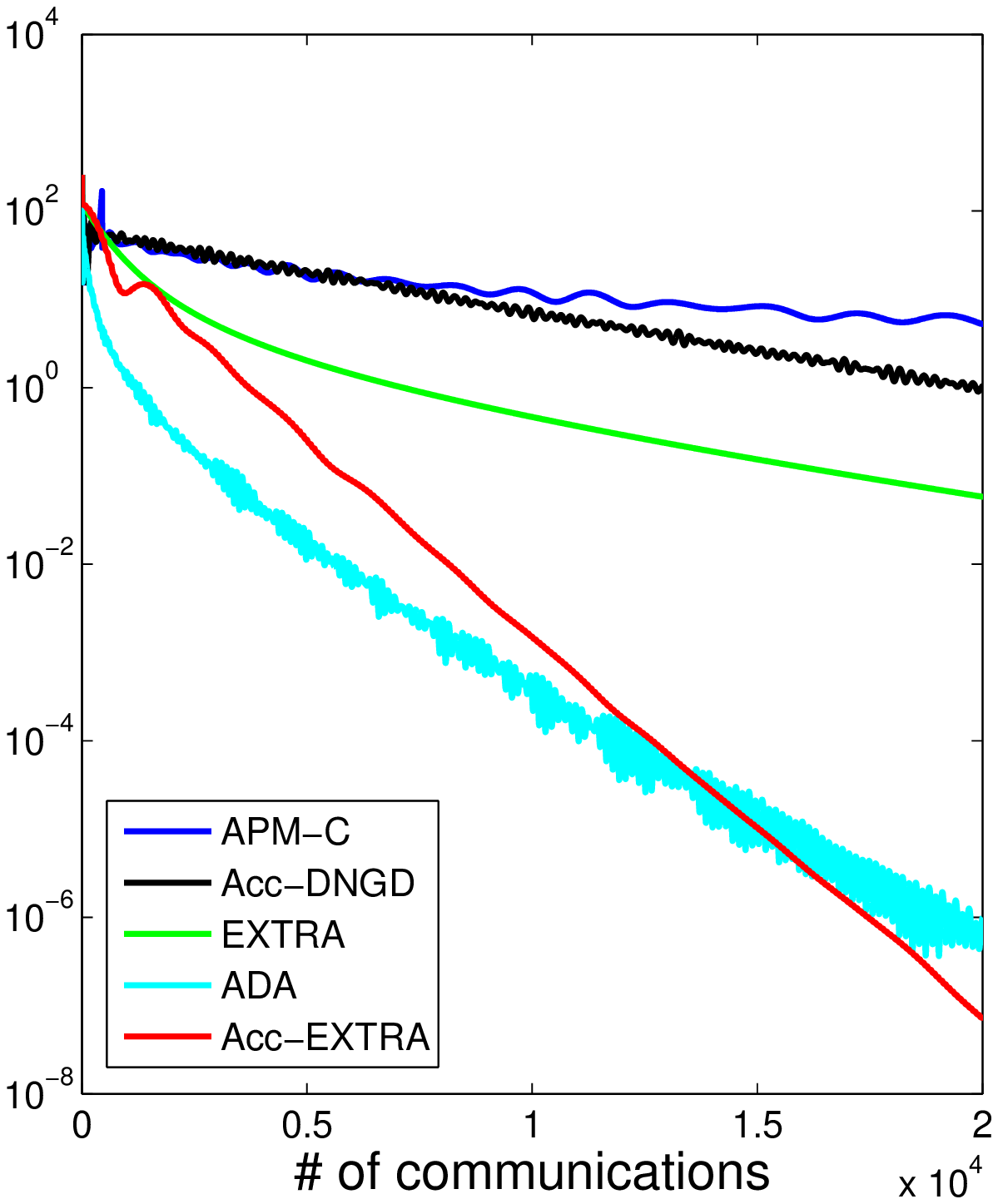}
&\hspace*{-0.28cm}\includegraphics[width=0.26\textwidth,keepaspectratio]{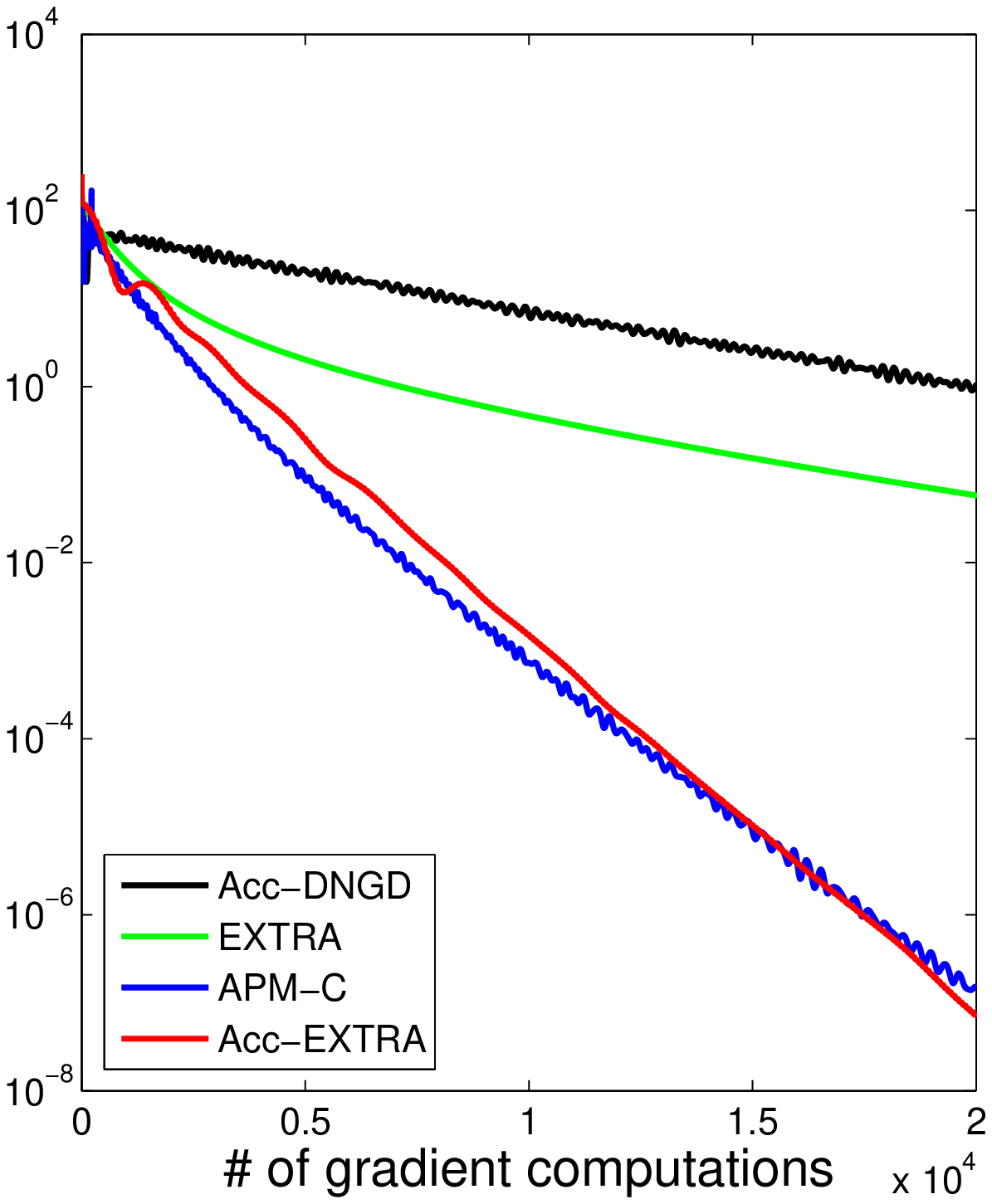}\\
\hspace*{-0.8cm}\includegraphics[width=0.325\textwidth,keepaspectratio]{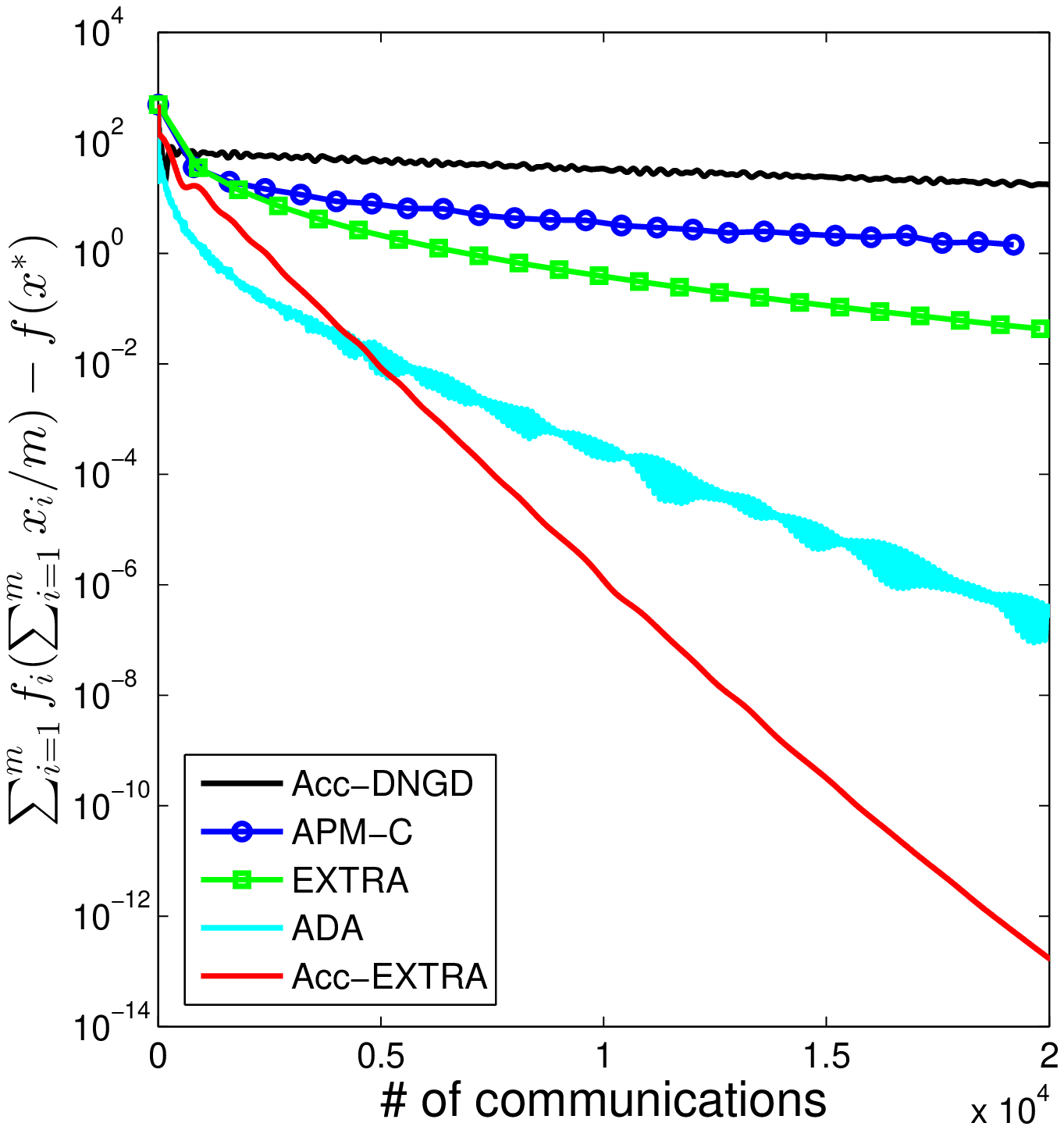}
&\hspace*{-0.28cm}\includegraphics[width=0.26\textwidth,keepaspectratio]{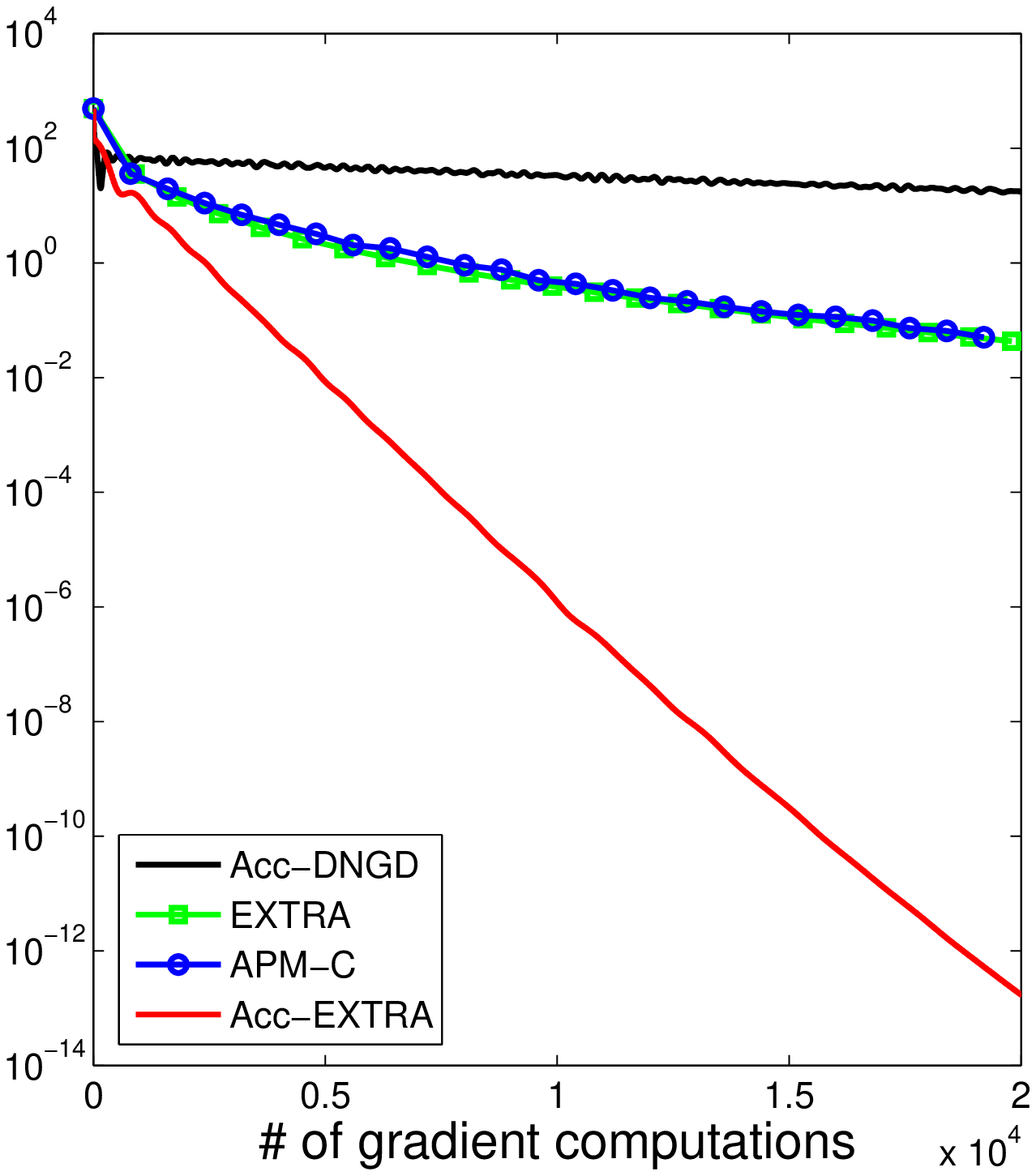}
&\hspace*{-0.28cm}\includegraphics[width=0.26\textwidth,keepaspectratio]{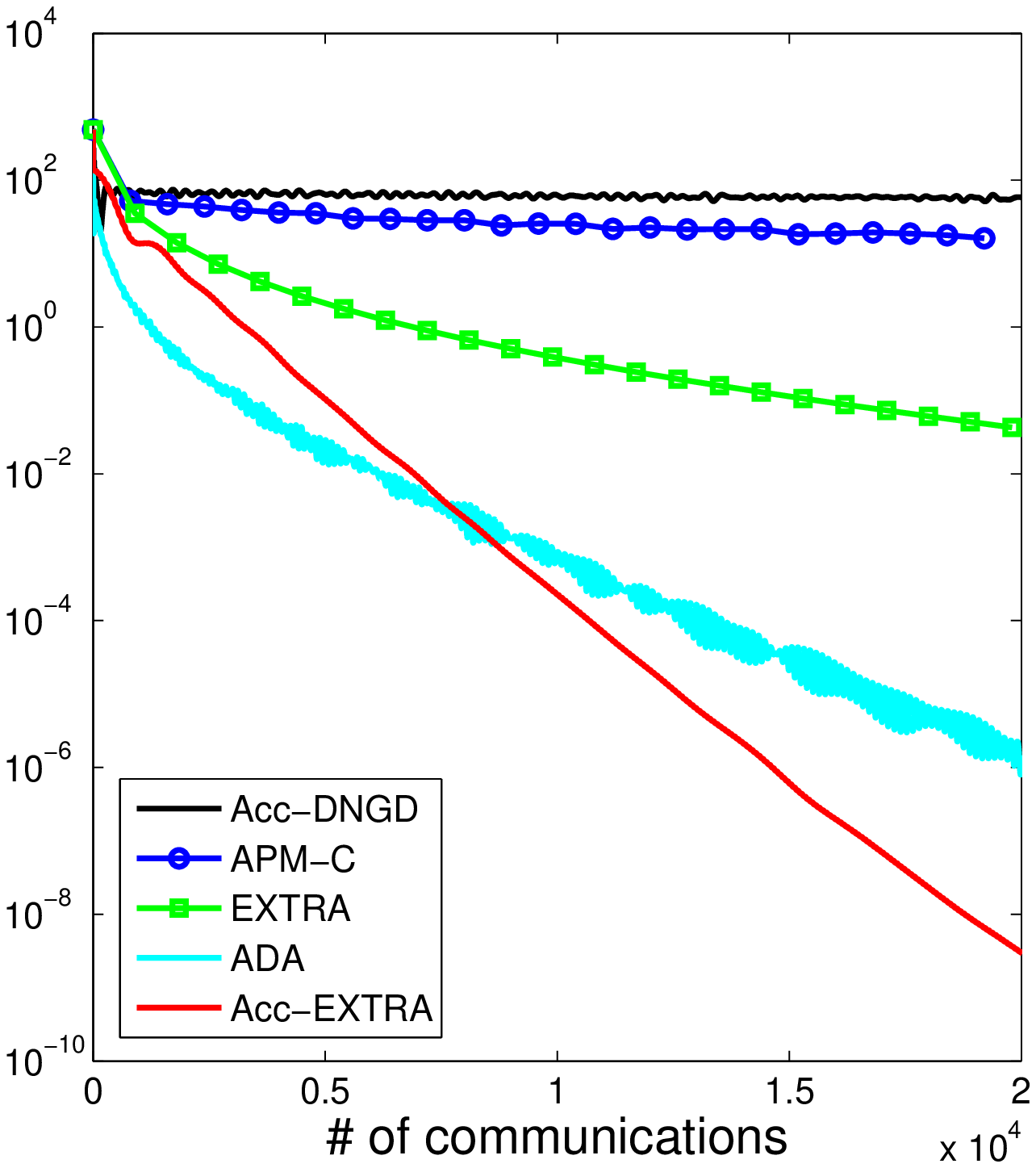}
&\hspace*{-0.28cm}\includegraphics[width=0.26\textwidth,keepaspectratio]{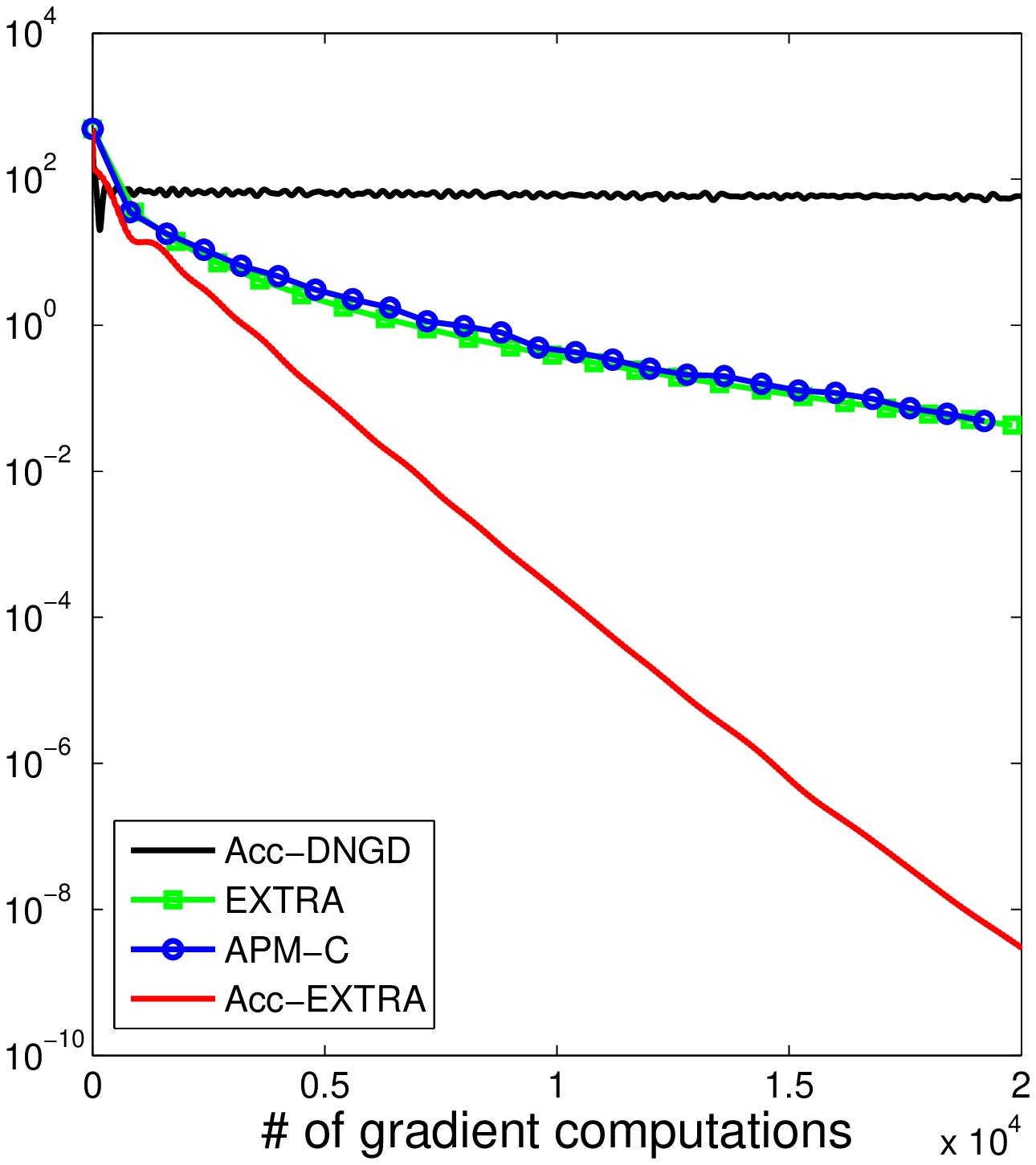}\\
\scriptsize$\frac{1}{1-\sigma_2(\W)}=8.13$&\scriptsize$\frac{1}{1-\sigma_2(\W)}=8.13$&\scriptsize$\frac{1}{1-\sigma_2(\W)}=30.02$&\scriptsize$\frac{1}{1-\sigma_2(\W)}=30.02$\normalsize
\end{tabular}
\caption{Comparisons on the strongly convex problem with the geometric graph. $d=0.5$ for the two left plots, and $d=0.3$ for the two right. $\mu=10^{-6}$ for the top four plots, and $\mu=10^{-8}$ for the bottom four.}\label{fig3}
\end{figure}

\begin{figure}
\centering
\begin{tabular}{@{\extracolsep{0.001em}}c@{\extracolsep{0.001em}}c@{\extracolsep{0.001em}}c@{\extracolsep{0.001em}}c@{\extracolsep{0.001em}}c}
\hspace*{-0.8cm}\includegraphics[width=0.33\textwidth,keepaspectratio]{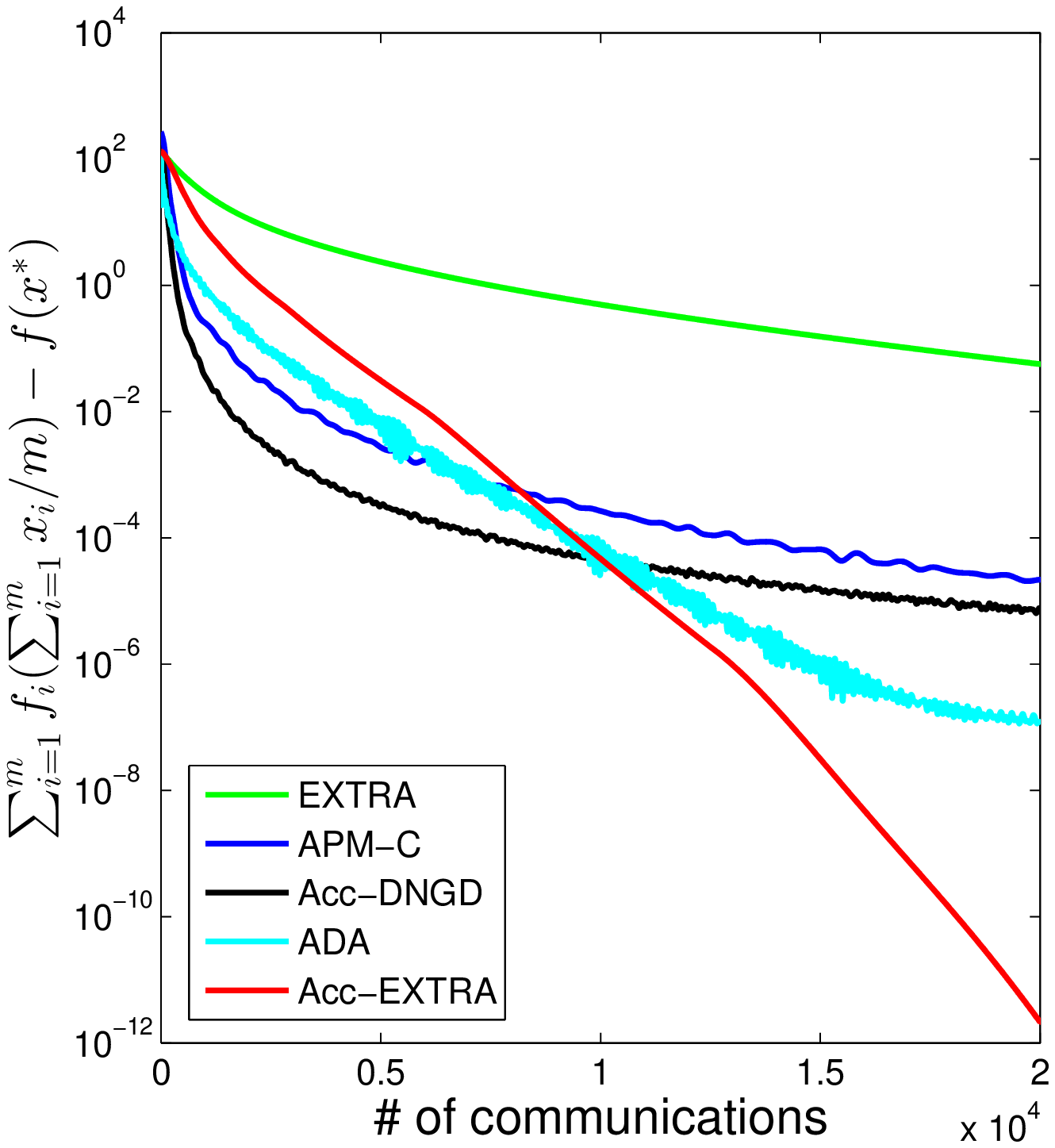}
&\hspace*{-0.28cm}\includegraphics[width=0.26\textwidth,keepaspectratio]{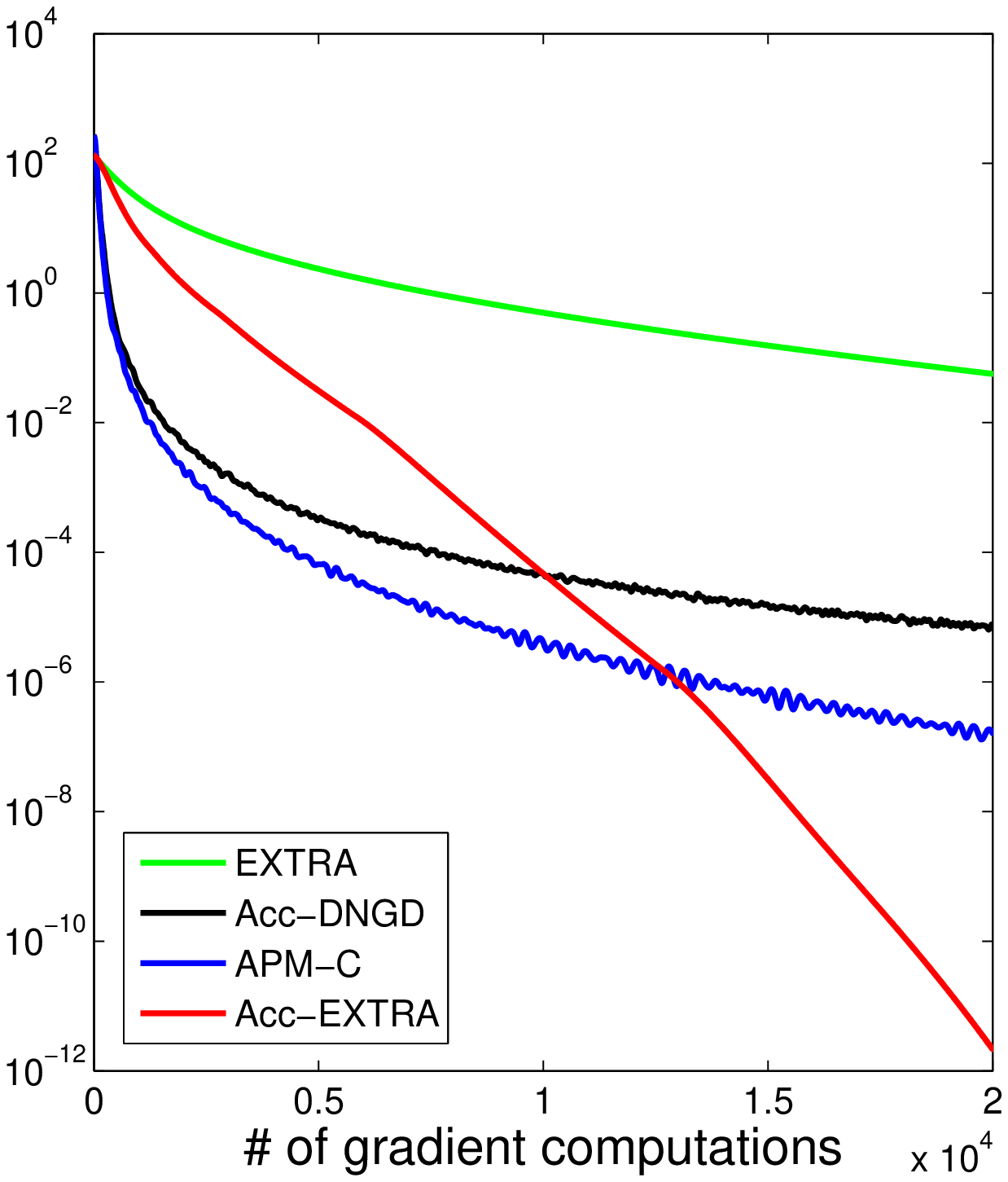}
&\hspace*{-0.28cm}\includegraphics[width=0.26\textwidth,keepaspectratio]{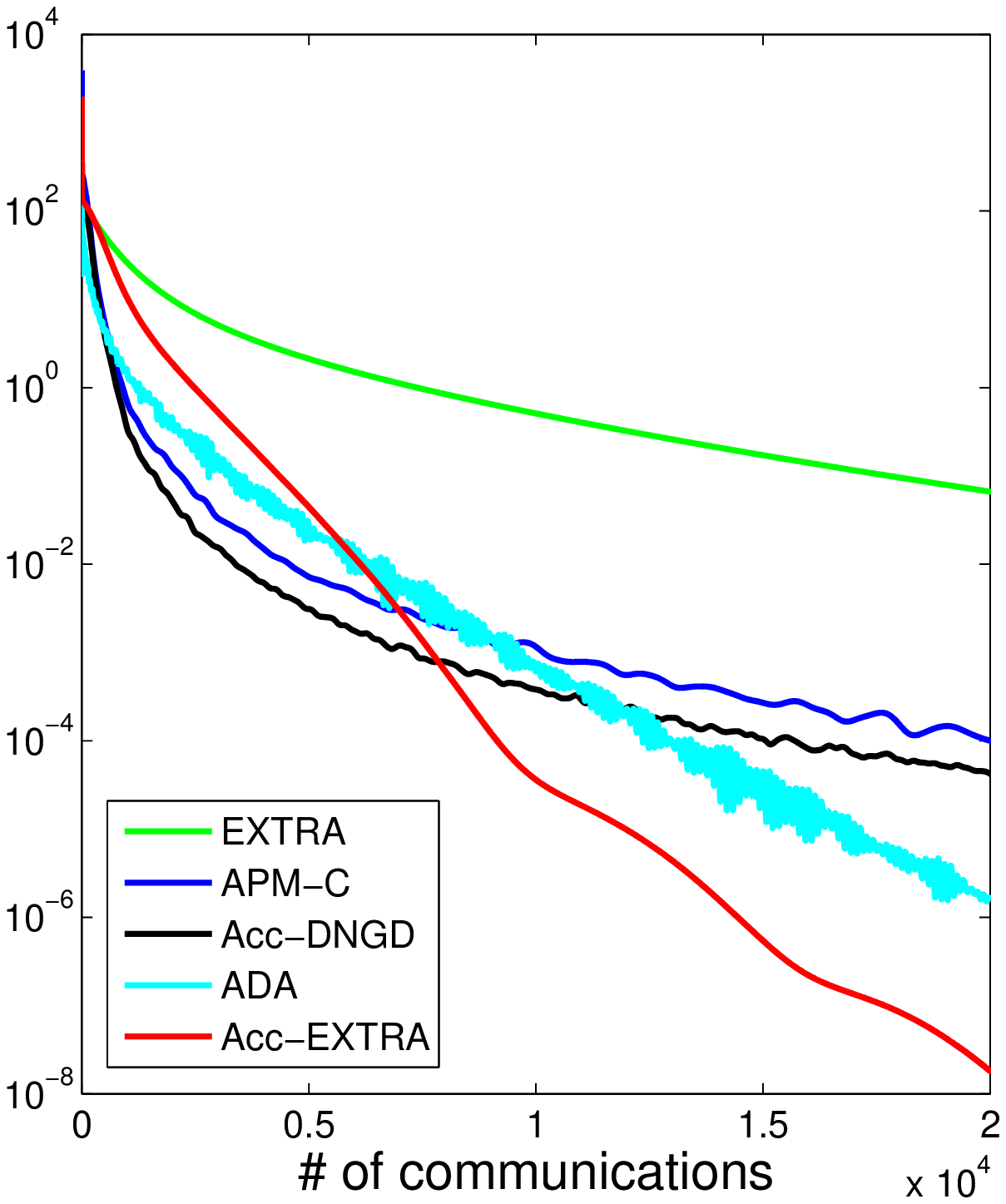}
&\hspace*{-0.28cm}\includegraphics[width=0.26\textwidth,keepaspectratio]{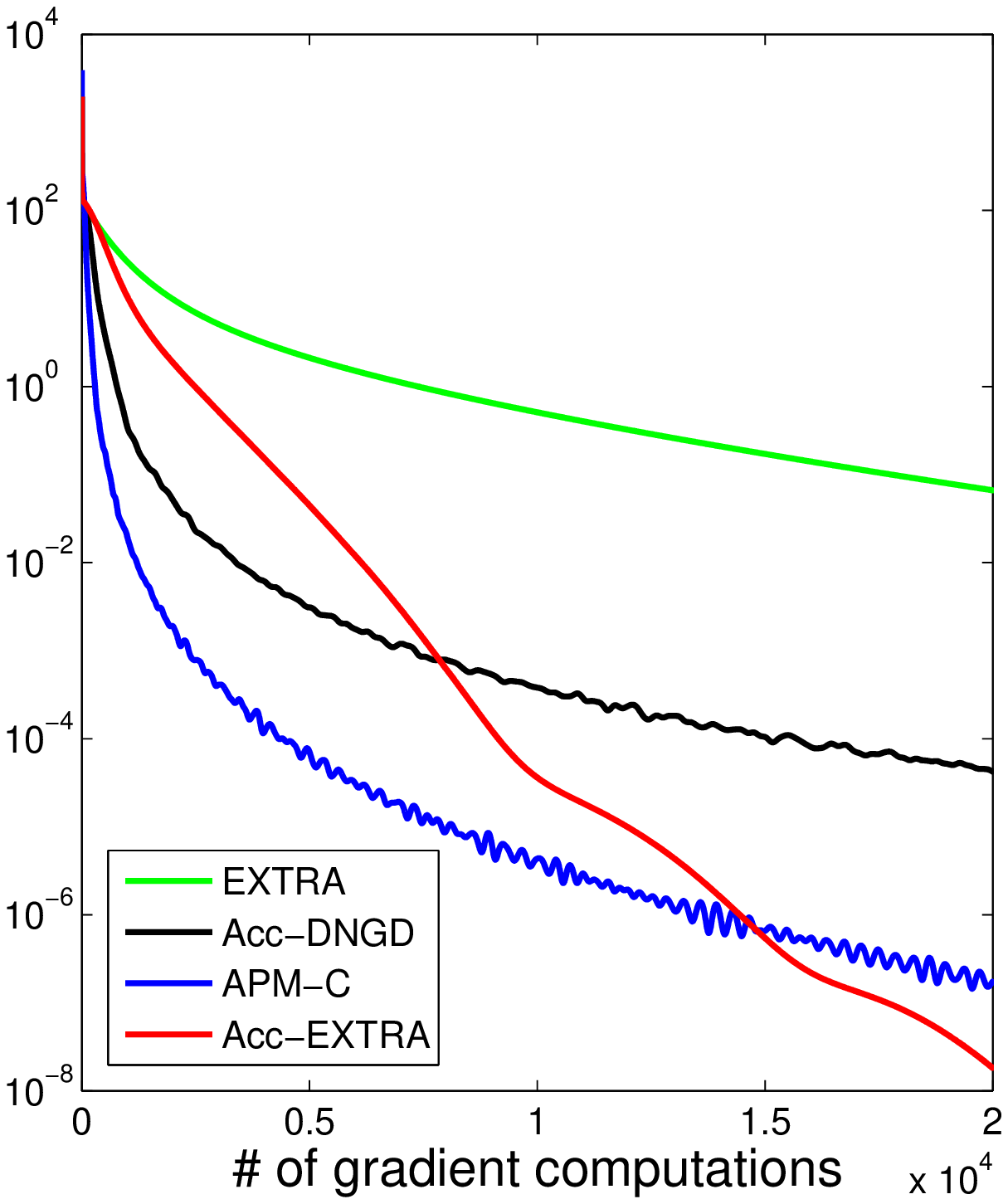}\\
\scriptsize$\frac{1}{1-\sigma_2(\W)}=2.87$&\scriptsize$\frac{1}{1-\sigma_2(\W)}=2.87$&\scriptsize$\frac{1}{1-\sigma_2(\W)}=7.74$&\scriptsize$\frac{1}{1-\sigma_2(\W)}=7.74$\normalsize
\end{tabular}
\caption{Comparisons on the nonstrongly convex problem with the Erd\H{o}s$-$R\'{e}nyi random graph. $p=0.5$ for the two left plots, and $p=0.1$ for the two right plots.}\label{fig2}
\end{figure}

\begin{figure}
\centering
\begin{tabular}{@{\extracolsep{0.001em}}c@{\extracolsep{0.001em}}c@{\extracolsep{0.001em}}c@{\extracolsep{0.001em}}c@{\extracolsep{0.001em}}c}
\hspace*{-0.8cm}\includegraphics[width=0.325\textwidth,keepaspectratio]{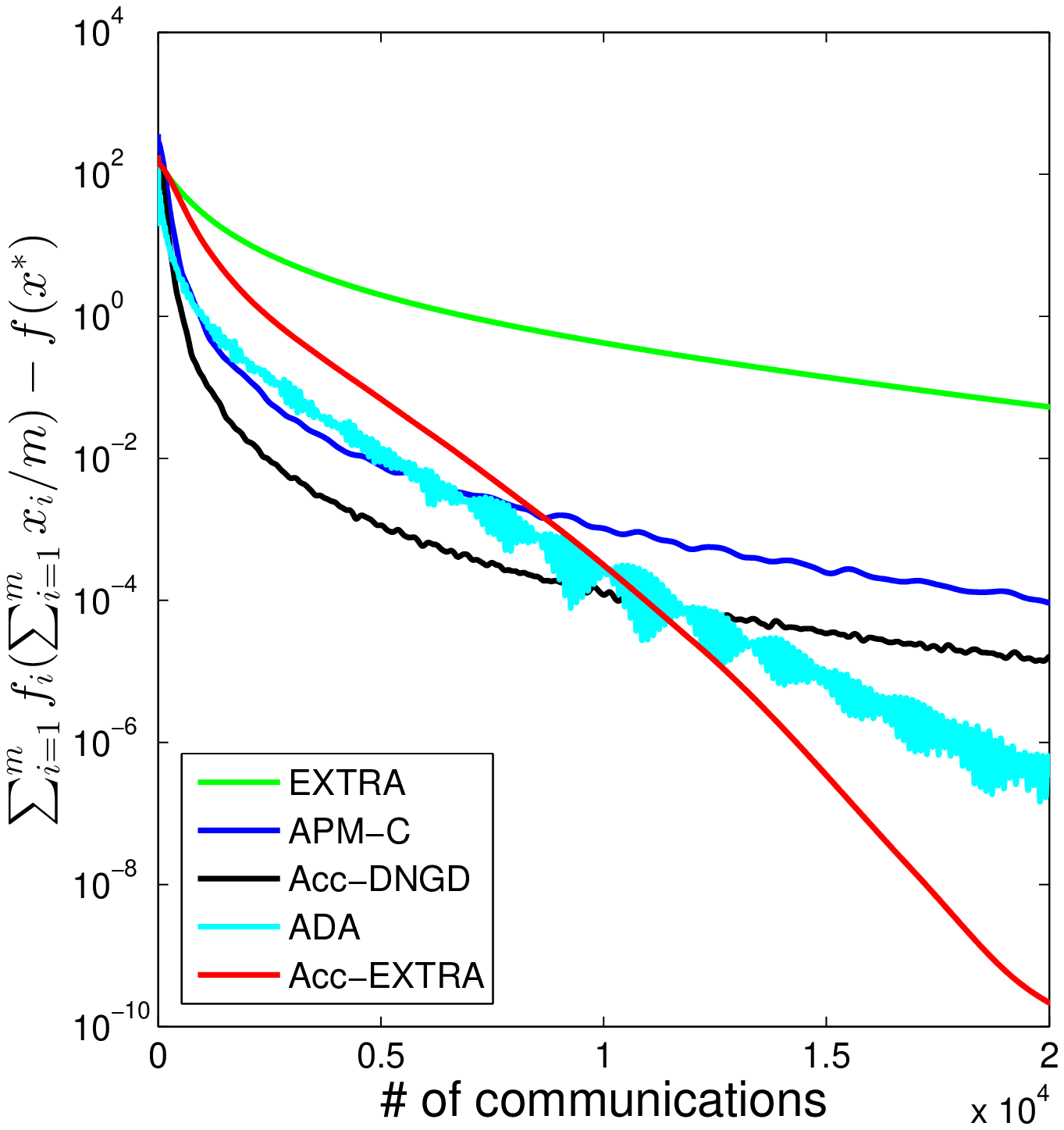}
&\hspace*{-0.28cm}\includegraphics[width=0.26\textwidth,keepaspectratio]{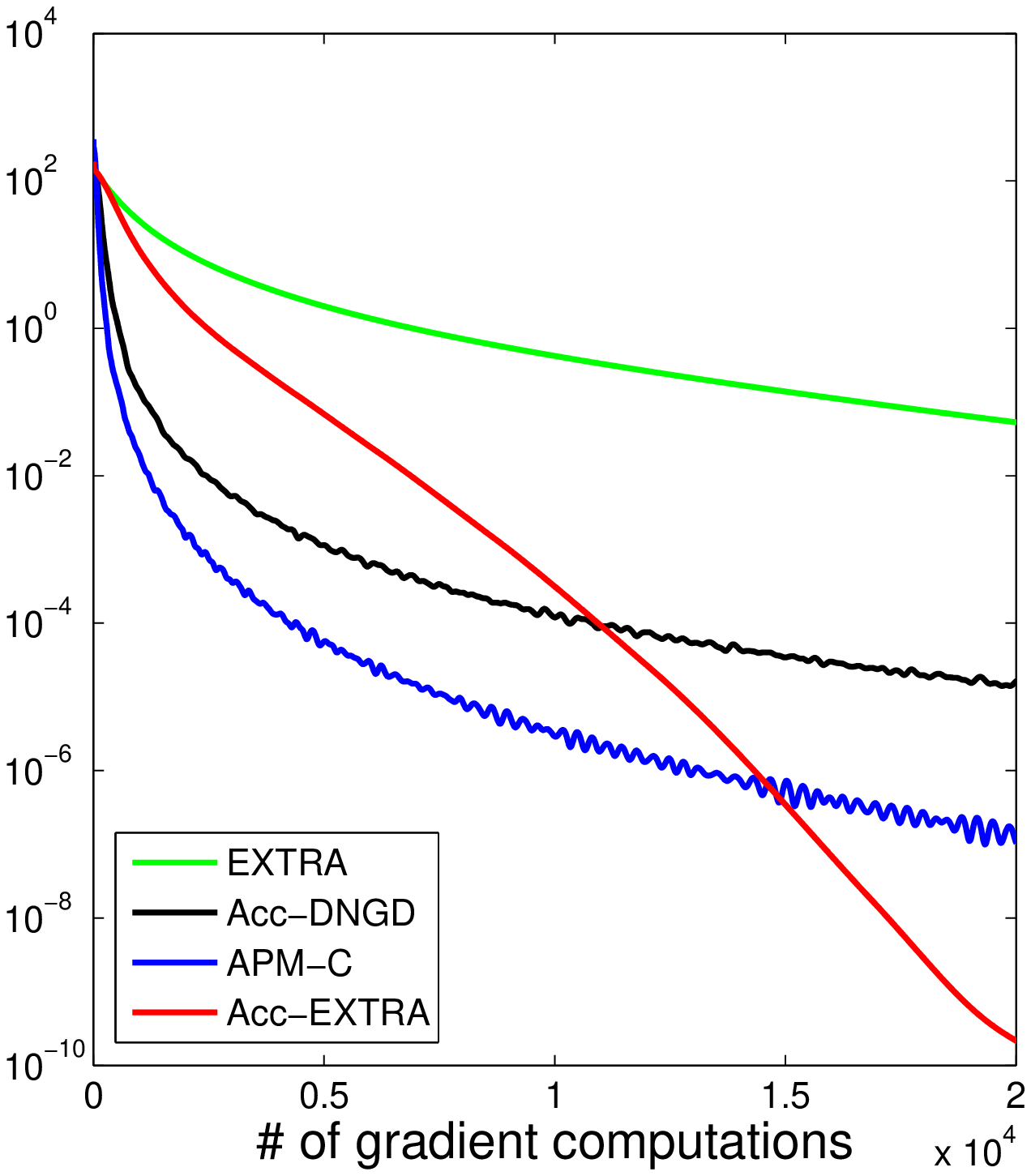}
&\hspace*{-0.28cm}\includegraphics[width=0.26\textwidth,keepaspectratio]{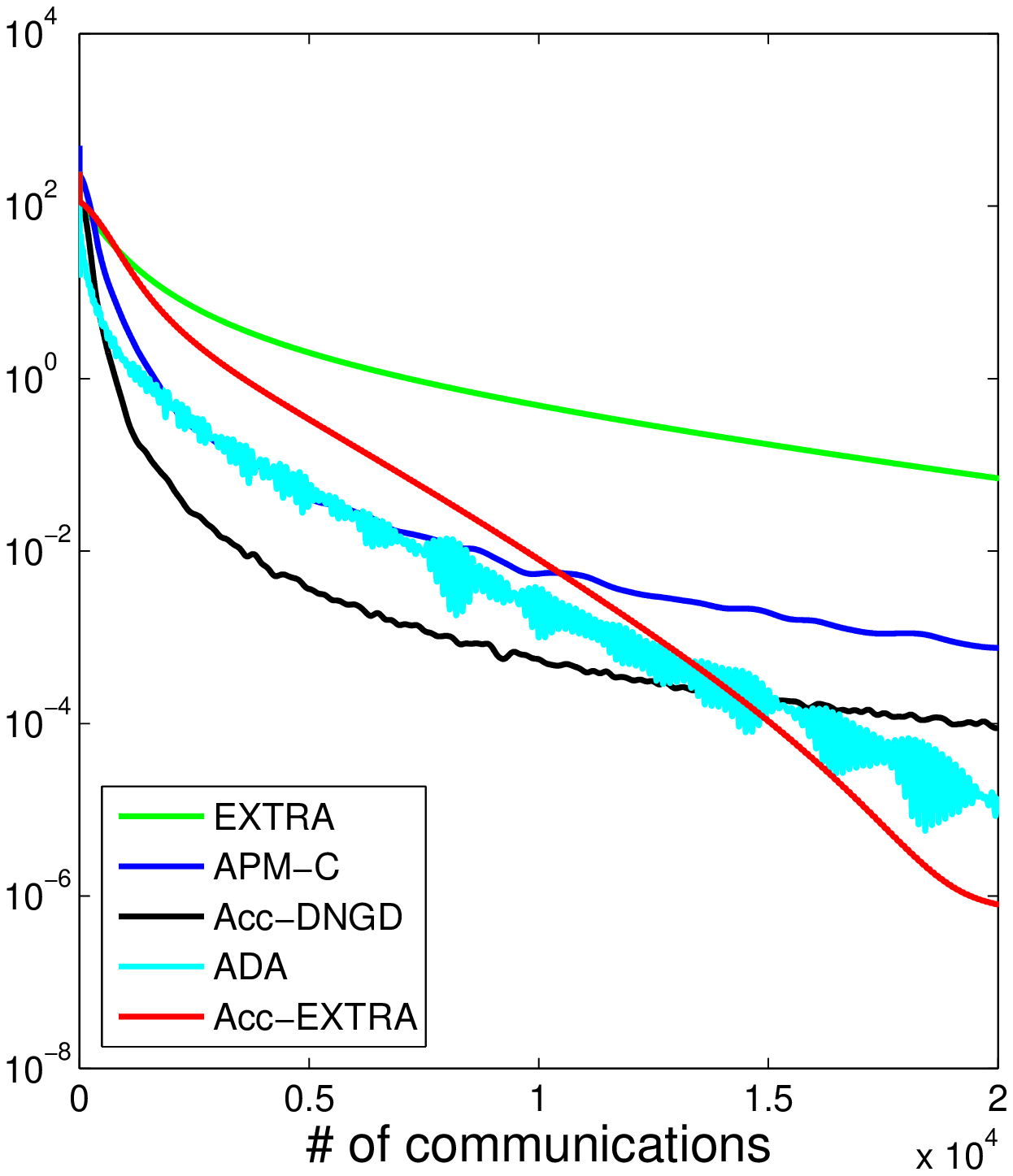}
&\hspace*{-0.28cm}\includegraphics[width=0.26\textwidth,keepaspectratio]{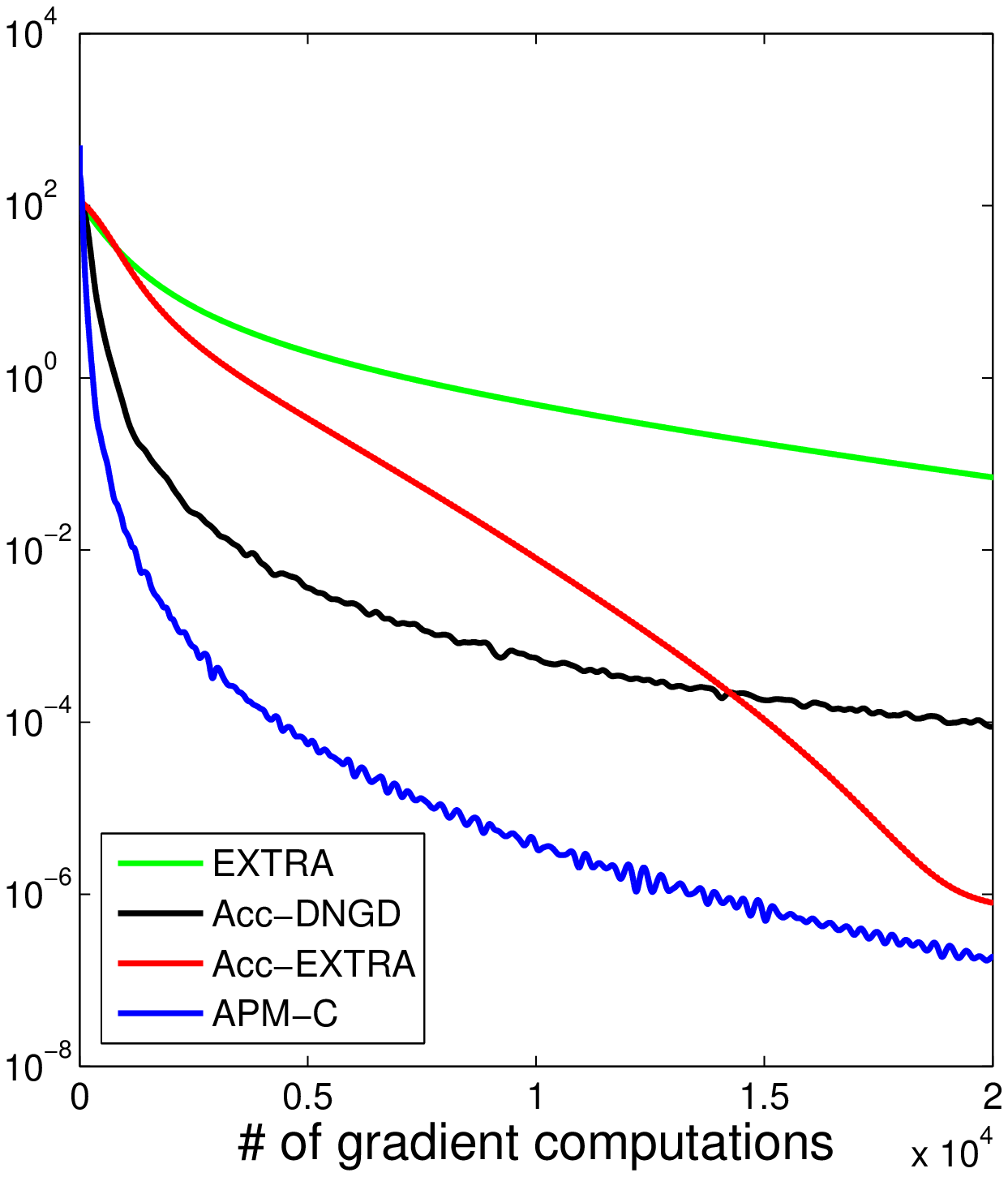}\\
\scriptsize$\frac{1}{1-\sigma_2(\W)}=8.13$&\scriptsize$\frac{1}{1-\sigma_2(\W)}=8.13$&\scriptsize$\frac{1}{1-\sigma_2(\W)}=30.02$&\scriptsize$\frac{1}{1-\sigma_2(\W)}=30.02$\normalsize
\end{tabular}
\caption{Comparisons on the nonstrongly convex problem with the geometric graph. $d=0.5$ for the two left plots, and $d=0.3$ for the two right plots.}\label{fig4}
\end{figure}

\begin{figure}
\centering
\begin{tabular}{@{\extracolsep{0.001em}}c@{\extracolsep{0.001em}}c@{\extracolsep{0.001em}}c@{\extracolsep{0.001em}}c@{\extracolsep{0.001em}}c}
\hspace*{-0.8cm}\includegraphics[width=0.325\textwidth,keepaspectratio]{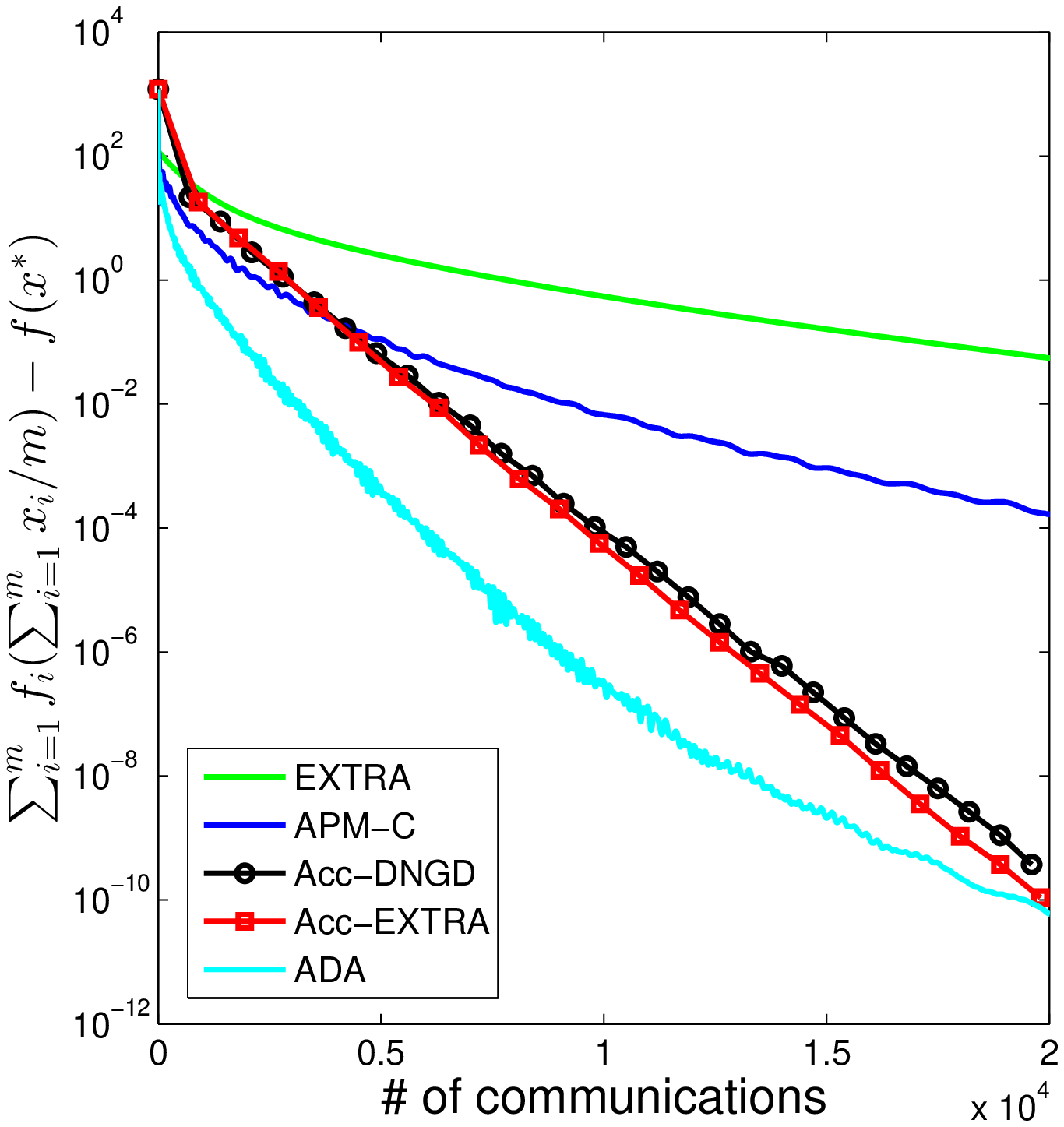}
&\hspace*{-0.28cm}\includegraphics[width=0.26\textwidth,keepaspectratio]{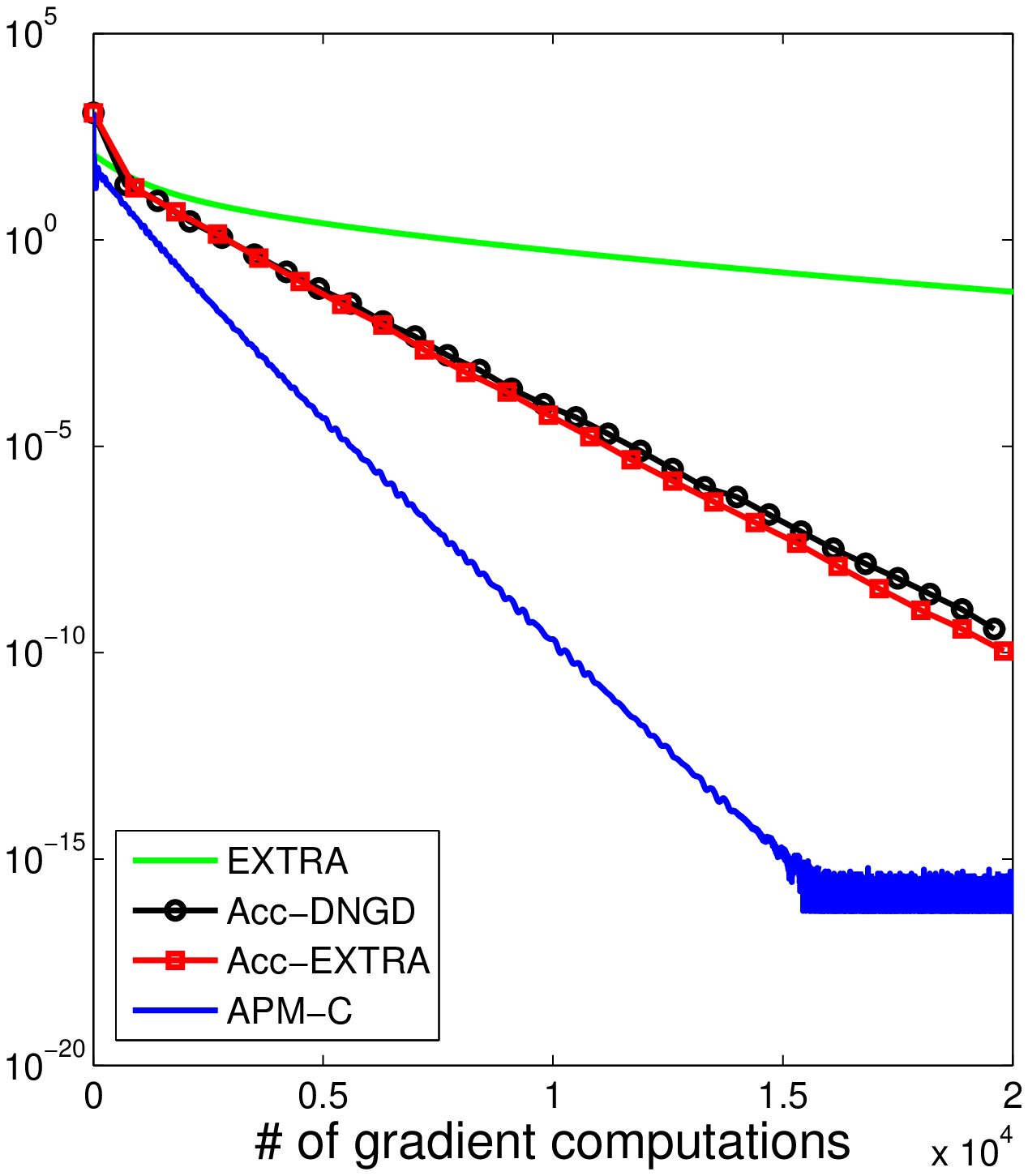}
&\hspace*{-0.28cm}\includegraphics[width=0.26\textwidth,keepaspectratio]{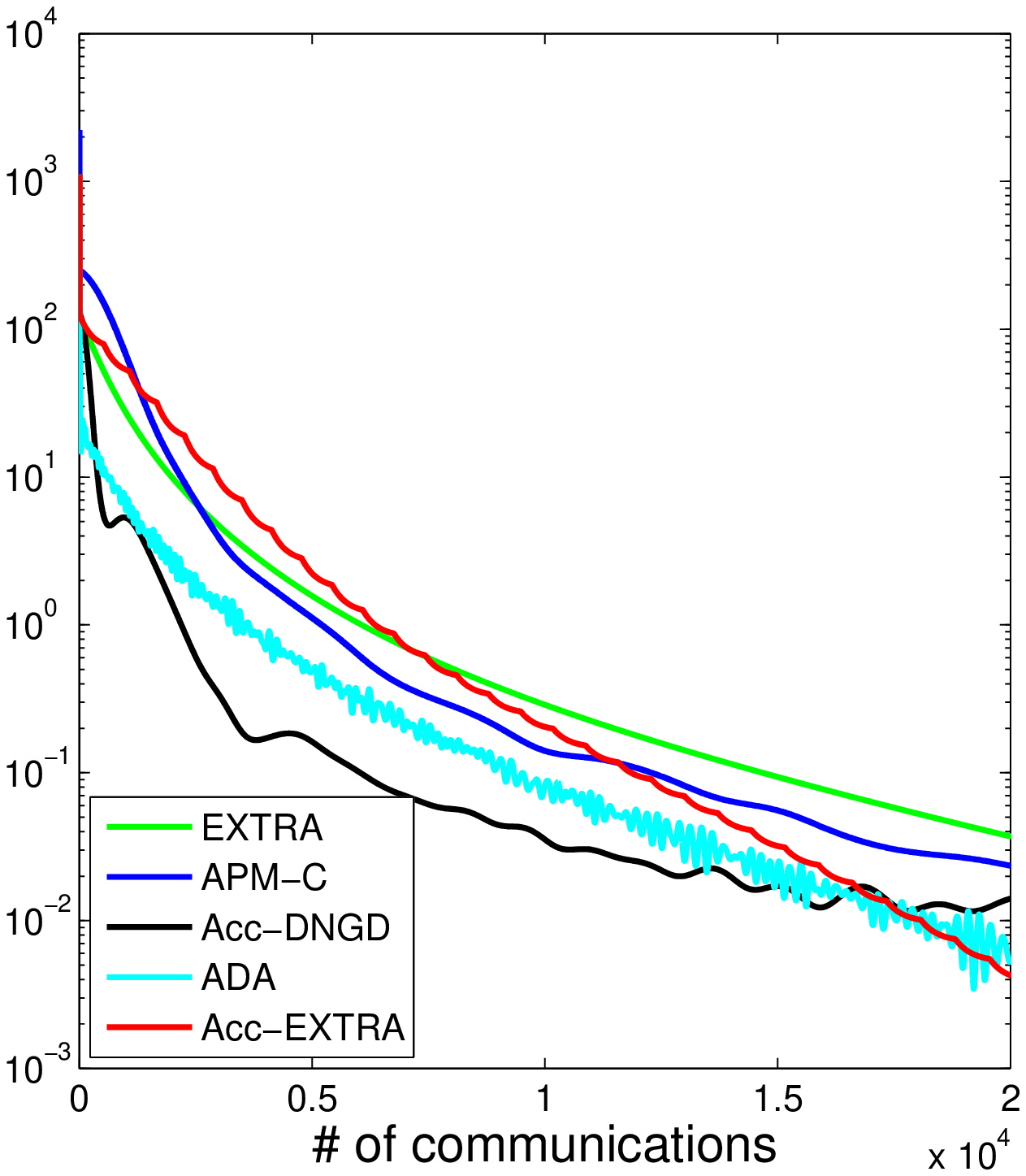}
&\hspace*{-0.28cm}\includegraphics[width=0.26\textwidth,keepaspectratio]{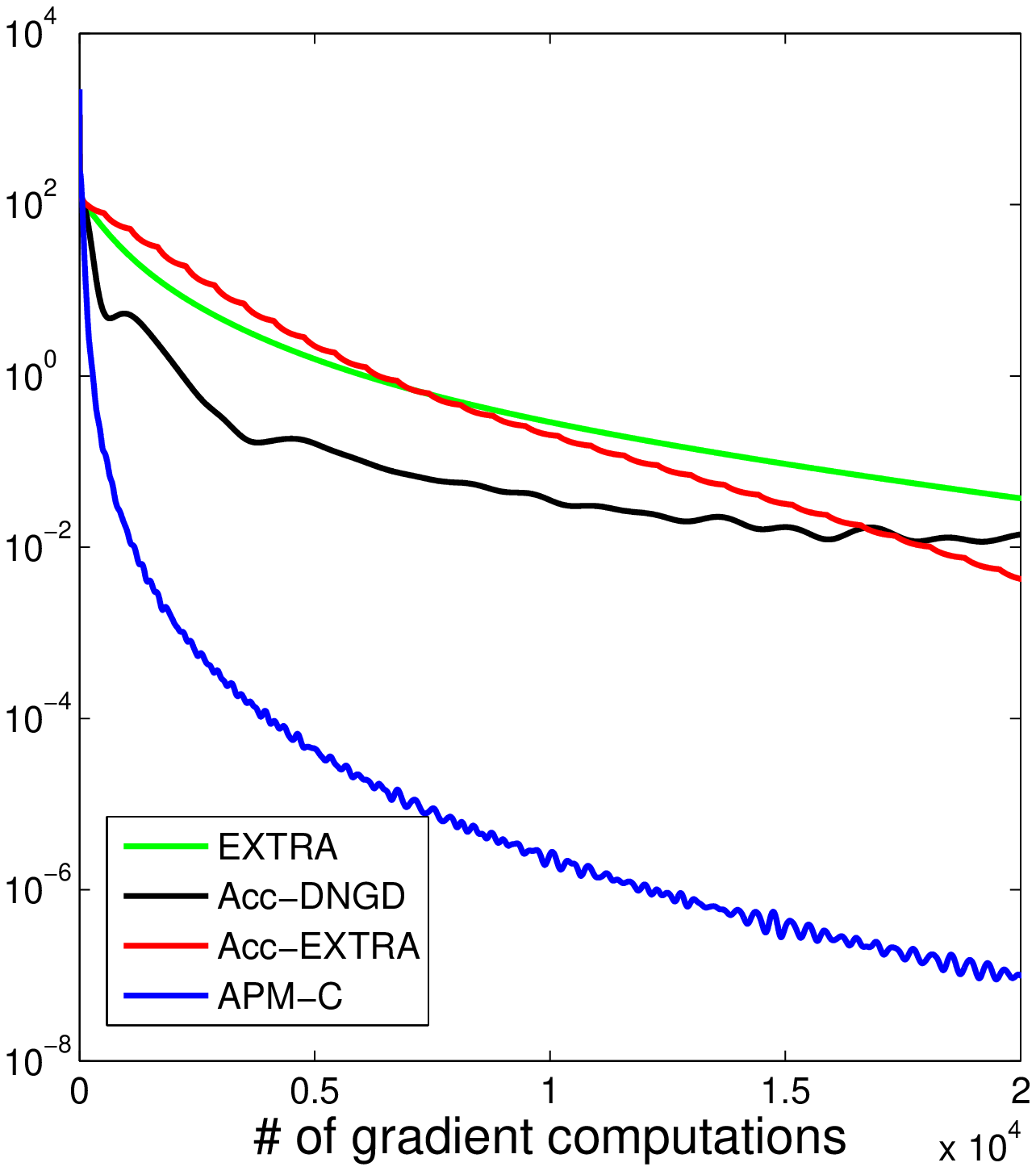}\\
\end{tabular}
\caption{Further comparisons on the geometric graph. $d=0.5$ and $\mu=10^{-5}$ for the two left plots, and $d=0.15$ and $\mu=0$ for the two right plots.}\label{fig5}
\end{figure}


The numerical results are illustrated in Figures \ref{fig1} and \ref{fig3}. The computation cost of ADA is high, and it has almost no visible decreasing in the first $20,000$ gradient computations [Figure 2]\cite{li-2018-pm}. Thus, we do not paint it in the second and fourth plots of Figures \ref{fig1}-\ref{fig4}. We can see that Acc-EXTRA performs better than the original EXTRA on both the Erd\H{o}s$-$R\'{e}nyi random graph and the geometric graph. We also observe that Acc-EXTRA is superior to ADA and APM-C on the graphs with small $\mu$ and $\frac{1}{1-\sigma_2(\W)}$. The performance of Acc-EXTRA degenerates when $\mu$ and $\frac{1}{1-\sigma_2(\W)}$ become larger. When preparing the experiments, we observe that Acc-EXTRA applies to ill-conditioned problems with large condition numbers for strongly convex problems. In this case, Acc-EXTRA runs with a certain number of outer iterations and the acceleration takes effect.

For the nonstrongly convex problem ($\mu=0$), we set $\tau=L(1-\sigma_2(\W))$ and $T_k=\lceil\frac{1}{2(1-\sigma_2(\W))}\log\frac{k+1}{1-\sigma_2(\W)}\rceil$ for Acc-EXTRA, $T_k=\lceil\frac{\log(k+1)}{5\sqrt{1-\sigma_2(\W)}}\rceil$ and the step-size as $\frac{1}{L}$ for APM-C. We tune the best step-size as $\frac{1}{L}$ and $\frac{0.2}{L}$ for EXTRA and Acc-DNGD, respectively. For ADA, we add a small regularizer of $\frac{\epsilon}{2}\|x\|^2$ to each $f_i(x)$ and solve a regularized strongly convex problem with $\epsilon=10^{-7}$. The numerical results are illustrated in Figures \ref{fig2} and \ref{fig4}. We observe that Acc-EXTRA also outperforms the original EXTRA and Acc-EXTRA is superior with small $\frac{1}{1-\sigma_2(\W)}$. Moreover, at the first 10000 iterations, the advantage of Acc-EXTRA is not obvious and it performs better at the last 5000 iterations. Thus, Acc-EXTRA is suited for the applications requiring high precision and the well-connected networks with small $\frac{1}{1-\sigma_2(\W)}$.

Finally, we report two results in Figure \ref{fig5} that Acc-EXTRA does not perform well, where the two left plots are for the strongly convex problem and the two right ones are for the nonstrongly convex one. Comparing the two left plots in Figure \ref{fig5} with the left and top two plots in Figure \ref{fig3}, we can see that Acc-EXTRA is inferior to ADA and APM-C in cases with a larger $\mu$, i.e., a smaller condition number for strongly convex problems. On the other hand, comparing the two right ones in Figure \ref{fig5} with the four plots in Figure \ref{fig4}, we observe that ADA and APM-C outperform Acc-EXTRA in cases with a larger $\frac{1}{1-\sigma_2(\W)}$ (it equals 268.67 when $d=0.15$) for nonstrongly convex problems. These observations further support the above conclusions.


\section{Conclusion}\label{sec:conclude}

In this paper, we first give a sharp analysis on the original EXTRA with improved complexities, which depends on the sum of $\frac{L}{\mu}$ (or $\frac{L}{\epsilon}$) and $\frac{1}{1-\sigma_2(\W)}$, rather than their product. Then, we use the Catalyst framework to accelerate it and obtain the near optimal communication complexities and competitive computation complexities. Our communication complexities of the proposed accelerated EXTRA are only worse by the factors of $(\log\frac{L}{\mu(1-\sigma_2(\W))})$ and $(\log\frac{1}{\epsilon})$ form the lower bounds for strongly convex and nonstrongly convex problems, respectively. 

\bibliography{arxiv}
\bibliographystyle{icml2020}
\end{document}